\documentclass[11pt]{article}

\usepackage{epsfig,epsf,fancybox}
\usepackage{amsmath}
\usepackage{mathrsfs}
\usepackage{amssymb}
\usepackage{graphicx}
\usepackage{color}
\usepackage{multirow}
\usepackage{paralist}
\usepackage{verbatim}
\usepackage{galois}
\usepackage{algorithm}
\usepackage{algorithmic}
\usepackage{boxedminipage}
\usepackage{booktabs}
\usepackage{accents}
\usepackage{stmaryrd}
\usepackage{subcaption}
\usepackage{wrapfig}
\usepackage[bottom]{footmisc}
\usepackage{natbib}

\usepackage[colorlinks,linkcolor=magenta,citecolor=blue, pagebackref=true,backref=true]{hyperref}
\renewcommand*{\backrefalt}[4]{%
    \ifcase #1 \footnotesize{(Not cited.)}%
    \or        \footnotesize{(Cited on page~#2.)}%
    \else      \footnotesize{(Cited on pages~#2.)}%
    \fi}

\long\def\comment#1{}
\usepackage{nicefrac}

\newtheorem{theorem}{Theorem}[section]

\newtheorem{lemma}[theorem]{Lemma}
\newtheorem{proposition}[theorem]{Proposition}

\newtheorem{definition}{Definition}[section]

\newtheorem{remark}{Remark}[section]
\newtheorem{assumption}{Assumption}[section]

\numberwithin{equation}{section}

\newcommand{\BB}{\mathbb{B}}

\newcommand{\Tg}{\textnormal{T}}

\newcommand{\st}{\textnormal{s.t.}}

\newcommand{\subdiff}{\textnormal{subdiff}\,}

\newcommand{\retr}{\textnormal{Retr}}

\newcommand{\op}{\textnormal{op}\,}

\newcommand{\Var}{\textnormal{Var}\,}

\newcommand{\argmin}{\mathop{\rm argmin}}

\newcommand{\ECal}{\mathcal{E}}

\newcommand{\MCal}{\mathcal{M}}
\newcommand{\SCal}{\mathcal{S}}
\newcommand{\XCal}{\mathcal{X}}
\newcommand{\YCal}{\mathcal{Y}}

\newcommand{\card}{\textnormal{card}}

\newcommand{\proj}{\textnormal{proj}}

\newcommand{\St}{\textnormal{St}}

\newcommand{\br}{\mathbb{R}}
\newcommand{\bs}{\mathbb{S}}
\newcommand{\bn}{\mathbb{N}}
\newcommand{\ba}{\begin{array}}
\newcommand{\ea}{\end{array}}

\newcommand{\ZCal}{\mathcal{Z}}
\newcommand{\WCal}{\mathcal{W}}

\newcommand{\FCal}{\mathcal{F}}

\newcommand{\PCal}{\mathcal{P}}
\newcommand{\PWCal}{\mathcal{PW}}
\newcommand{\PScr}{\mathscr{P}}
\newcommand{\EE}{{\mathbb{E}}}
\newcommand{\PP}{\mathbb{P}}

\newcommand{\NCal}{\mathcal{N}}

\newcommand{\UCal}{\mathcal{U}}
\newcommand{\one}{\textbf{1}}
\newcommand{\zero}{\textbf{0}}

\newcommand{\mydefn}{:=}

\begin{document}

\begin{center}

{\bf{\LARGE{On Projection Robust Optimal Transport: Sample Complexity and Model Misspecification}}}

\vspace*{.2in}
{\large{
\begin{tabular}{ccccc}
Tianyi Lin$^\diamond$ & Zeyu Zheng$^\star$ & Elynn Y. Chen$^\diamond$ & Marco Cuturi$^{\triangleleft, \triangleright}$ & Michael I. Jordan$^{\diamond, \dagger}$ \\
 \end{tabular}
}}

\vspace*{.2in}

\begin{tabular}{c}
Department of Electrical Engineering and Computer Sciences$^\diamond$ \\
Department of Industrial Engineering and Operations Research$^\star$ \\
Department of Statistics$^\dagger$ \\ 
University of California, Berkeley \\
CREST - ENSAE$^\triangleleft$, Google Brain$^\triangleright$
\end{tabular}

\vspace*{.2in}

\today

\vspace*{.2in}

\begin{abstract}
Optimal transport (OT) distances are increasingly used as loss functions for statistical inference, notably in the learning of generative models or supervised learning. Yet, the behavior of minimum Wasserstein estimators is poorly understood, notably in high-dimensional regimes or under model misspecification. In this work we adopt the viewpoint of projection robust (PR) OT, which seeks to maximize the OT cost between two measures by choosing a $k$-dimensional subspace onto which they can be projected. Our first contribution is to establish several fundamental statistical properties of PR Wasserstein distances, complementing and improving previous literature that has been restricted to one-dimensional and well-specified cases. Next, we propose the integral PR Wasserstein (IPRW) distance as an alternative to the PRW distance, by averaging rather than optimizing on subspaces. Our complexity bounds can help explain why both PRW and IPRW distances outperform Wasserstein distances empirically in high-dimensional inference tasks. Finally, we consider parametric inference using the PRW distance. We provide an asymptotic guarantee of two types of minimum PRW estimators and formulate a central limit theorem for max-sliced Wasserstein estimator under model misspecification. To enable our analysis on PRW with projection dimension larger than one, we devise a novel combination of variational analysis and statistical theory. 
\end{abstract}
\end{center}

\section{Introduction}
Recent years have witnessed an ever-increasing role for ideas from optimal transport (OT)~\citep{Villani-2008-Optimal} in machine learning. Combining OT distances with the general principles of minimal distance estimation (MDE)~\citep{Wolfowitz-1957-Minimum, Basu-2011-Statistical} yields a powerful basis for various statistical inference problems, such as density estimation~\cite{Bassetti-2006-Minimum}, training of generative model~\citep{Arjovsky-2017-Wasserstein, Gulrajani-2017-Improved, Montavon-2016-Wasserstein, Adler-2018-Banach, Cao-2019-Multi}, auto-encoders~\citep{Tolstikhin-2018-Wasserstein}, clustering~\citep{Cuturi-2014-Fast, Bonneel-2016-Barycoord, Ho-2017-Multilevel, Ye-2017-Fast}, multitask regression~\citep{Janati-2020-Multi}, trajectory inference~\citep{Hashimoto-2016-Learning, Schiebinger-2017-Reconstruction, Yang-2020-Predicting, Tong-2020-Trajectorynet} or nonparametric testing~\citep{Ramdas-2017-Wasserstein}; see~\citet{Peyre-2019-Computational} and~\citet{Panaretos-2019-Statistical} for reviews on these topics.  

For OT ideas to continue to bear fruit in machine learning, it will be necessary to tackle two characteristic challenges: (1) high dimensionality and (2) model misspecification.  Initial progress has been made on the latter problem by~\citet{Bernton-2019-Parameter}, who showed that in the misspecified case the minimum Wasserstein estimator (MWE) outputs the Wasserstein projection of the data-generating distribution onto the fitted model class.  These authors also obtained results on robustness and the asymptotic distribution of the projection, while these results only apply to the one-dimensional setting. High-dimensional settings are challenging; indeed, it is known that the sample complexity of estimating the Wasserstein distance can grow exponentially in dimension~\citep{Dudley-1969-Speed, Fournier-2015-Rate, Singh-2018-Minimax, Weed-2019-Sharp, Lei-2020-Convergence}.

We focus on a promising approach to treating high-dimensional problems: Compute the OT distance between low-dimensional projections of high-dimensional input measures. The simplest and most representative example of this approach is the sliced Wasserstein distance~\citep{Rabin-2011-Wasserstein, Bonnotte-2013-Unidimensional, Bonneel-2015-Sliced, Deshpande-2019-Max, Kolouri-2019-Generalized, Nadjahi-2020-Statistical, Manole-2019-Minimax}, which is defined as the average OT distance obtained between random 1-dimensional projections, and which is shown practical in real applications~\citep{Deshpande-2018-Generative, Deshpande-2019-Max, Kolouri-2016-Sliced, Kolouri-2018-Sliced, Carriere-2017-Sliced, Wu-2019-Sliced, Liutkus-2019-Sliced}. In an important extension,~\citet{Paty-2019-Subspace} and~\citet{Weed-2019-Estimation} proposed very recently to seek the $k$-dimensional subspace ($k>1$) that would maximize the OT distance between two measures after projection. The quantity is named as~\textit{projection robust Wasserstein} (PRW) distance\footnote{This quantity is also named as~\textit{Wasserstein Projection Pursuit} (WPP)~\citep{Weed-2019-Estimation}. For simplicity, we refer from now on to PRW/WPP as PRW.}, which is conceptually simple and does solve the curse of dimensionality in the so-called spiked model as proved in~\citep[Theorem 1]{Weed-2019-Estimation} by recovering the $n^{-1/k}$ rate under the Talagrand transport inequality. This result suggests that PRW can be significantly more useful than the OT distance for inference tasks when the dimension is large. From a computational point of view, PRW becomes the max-sliced Wasserstein distance when the projection dimension is $k=1$ and has an efficient implementation~\citep{Deshpande-2019-Max}. For general $k \geq 1$,~\citet{Lin-2020-Projection} proposed to compute PRW using Riemannian optimization toolbox and provided theoretical guarantee and encouraging empirical results. However, it is desirable to understand its statistical behavior which mostly determines the practical performance of PRW.

\paragraph{Contributions.} In this paper, we study the statistical properties of PRW and another so-called integrated PRW (IPRW), which replaces the maximum in the original PRW with an average of OT distance over $k$-dimensional projections. Our contributions can be summarized as follows. 
\begin{enumerate}
\item We prove that the empirical measure $\widehat{\mu}_n$ converges to true measure $\mu_\star$ under both PRW and IPRW with different rates. These rates are new to our knowledge. For example, when the order $p=3/2$ and the projected dimension $k \geq 3$, the rate is $n^{-1/k}$ for IPRW. For PRW, the rate is $(n^{-1/k} + n^{-1/6}\sqrt{dk\log(n)} + n^{-2/3}dk\log(n))$ when $\mu_\star$ satisfies a projection Bernstein tail condition and $(n^{-1/k} + n^{-1/2}\sqrt{dk\log(n)} + n^{-2/3}dk\log(n))$ when $\mu_\star$ satisfies a projection Poincar\'{e} inequality.

\item We derive the concentration results when $\mu_\star$ satisfies a Bernstein tail condition or a projection Poincar\'{e} inequality. In terms of tail conditions, our Bernstein condition and Poincare inequality handle subexponential tail while Talagrand inequality in~\citet{Weed-2019-Estimation} addresses subgaussian tail. Our assumptions are thus weaker than~\citet{Weed-2019-Estimation}.
 
\item We establish asymptotic guarantees for the minimal PRW and expected minimal PRW estimators under model misspecification. For minimal PRW estimator with the order $p=1$ and the projected dimension $k=1$, we derive an asymptotic distribution for arbitrary dimension $d$ with the $n^{-1/2}$ rate in the Hausdorff metric. Our assumptions are weaker than those used in~\citet{Bernton-2019-Parameter}, not requiring the nonsingularity of the Jacobian or the separability of the parameters. Our techniques for CLT in misspecified settings did not appear in~\citet{Nadjahi-2019-Asymptotic} and complete the analysis in~\citet{Bernton-2019-Parameter}.

\item We conduct experiments on synthetic data and neural networks to validate our theory. We also present a simple optimization algorithm that can efficiently compute the PRW distance in practice even when $k \geq 2$; see Appendix~\ref{appendix:computational} or Appendix B of the concurrent work by\citet{Lin-2020-Projection}. 
\end{enumerate}

\section{Preliminaries on Projected Optimal Transport} \label{sec:preliminaries}
In this section, we provide some technical background materials on projection optimal transport. Throughout the paper, we denote $\|\cdot\|$ as the Euclidean norm (in the corresponding vector space) and $\Rightarrow$ as the convergence in the weak sense.  
 
\paragraph{Wasserstein and sliced Wasserstein.} Let $p \geq 1$ and define $\PScr(\br^d)$ and $\PScr_p(\br^d)$ as the set of all Borel measures on $\br^d$ and the subset that satisfies $M_p(\mu) := \int_{\br^d} \|x\|^p d\mu(x) < +\infty$. For two probability measures $\mu, \nu \in \PScr_p(\br^d)$, their Wasserstein distance of order $p$ is defined as follows: 
\begin{equation}\label{def:Wasserstein}
\WCal_p(\mu, \nu) \mydefn \left(\inf_{\pi \in \Pi(\mu, \nu)} \int_{\br^d \times \br^d} \|x-y\|^p d\pi(x, y)\right)^{1/p},
\end{equation}
where the infimum is taken over $\Pi(\mu, \nu) \subseteq \PScr(\br^d \times \br^d)$---the set of probability measures with marginals $\mu$ and $\nu$. In the 1D case,~\citet[Theorem~3.1.2.(a)]{Rachev-1998-Mass} have shown that $\WCal_p(\mu, \nu) = (\int_0^1 |F_\mu^{-1}(t) - F_\nu^{-1}(t)|^p dt)^{1/p}$, where $F_\mu^{-1}$ and $F_\nu^{-1}$ are the quantile functions of $\mu$ and $\nu$. This 1D formula motivates the \textit{sliced Wasserstein (SW) and max-sliced Wasserstein (max-SW)} distances~\citep{Bonnotte-2013-Unidimensional, Bonneel-2015-Sliced, Deshpande-2019-Max}. In particular, the idea is to use as a proxy of \eqref{def:Wasserstein} the average or maximum of a set of 1D Wasserstein distances constructed by projecting $d$-dimensional measures to a random collection of 1D spaces. Computationally appealing, both SW and max-SW distances are widely used in practice, especially in generative modeling~\citep{Kolouri-2018-Sliced, Deshpande-2019-Max, Liutkus-2019-Sliced}. Practitioners observe that the SW distance only outputs a good Monte-Carlo approximation with a large number of projections, while the max-SW distance achieves similar results with fewer projections~\citep{Kolouri-2019-Generalized, Nguyen-2020-Distributional}. 

Encouraged by the success of SW and max-SW,~\citet{Paty-2019-Subspace} asked whether we can gain more by using a subspace of dimension $k \geq 2$, define the \textit{projection robust Wasserstein} (PRW) distance, and prove that this quantity is well posed if the order is $p \geq 1$. More specifically, let $\bs_{d, k} = \{E \in \br^{d \times k}: E^\top E = I_k\}$ be the set of $d \times k$ orthogonal matrices and $E^\star$ be the linear transformation associated with $E$ for any $x \in \br^d$ by $E^\star(x) = E^\top x$. For any measurable function $f$ and $\mu \in \PScr(\br^d)$, we denote $f_{\#}\mu$ as the push-forward of $\mu$ by $f$, so that $f_{\#}\mu(A) = \mu(f^{-1}(A))$ where $f^{-1}(A) = \{x \in \br^d: f(x) \in A\}$ for any Borel set $A$. For any given subspace dimension $K$, the PRW distance of order $p$ between $\mu$ and $\nu$ is defined by 
\begin{equation}\label{def:PRW}
\overline{\PWCal}_{p, k}(\mu, \nu) \mydefn \sup_{E \in \bs_{d, k}} \WCal_p(E_{\#}^\star\mu, E_{\#}^\star\nu).
\end{equation}
The PRW distance has better discriminative power than the SW or max-SW distances since it can extract more geometric information from high-dimensional projections than that from 1-dimensional projections; see~\citet{Paty-2019-Subspace} for more details.  

As an alternative, we define the IPRW distance, which replaces the supremum in Eq.~\eqref{def:PRW} with an average. The IPRW distance of order $p$ between $\mu$ and $\nu$ is
\begin{equation}\label{def:IPRW}
\underline{\PWCal}_{p, k}(\mu, \nu) \mydefn \left(\int_{\bs_{d, k}} \WCal_p^p(E_{\#}^\star\mu, E_{\#}^\star\nu) d\sigma(E)\right)^{1/p},   
\end{equation}
where $\sigma$ is the uniform distribution on $\bs_{d, k}$. Note that IPRW is well defined for comparing two measures and match our intuition. For example, given three Gaussian distributions $\mu_i = \NCal(u_i, I_d)$ for $i=1,2,3$, we have $\underline{\mathcal{PW}}_{p,2}(\mu_i, \mu_j)=c\|u_i-u_j\|$ where $c>0$ only depend on $p$ and the dimension $d$.    

The IPRW and PRW distances generalize the SW and max-SW distances to the high-dimensional projection setting. Both PRW and and IPRW are distances and satisfy the triangle inequality: the proof for PRW is in~\citet[Proposition~1]{Paty-2019-Subspace}, while that for IPRW is the same as that for SW in~\citet{Bonnotte-2013-Unidimensional}. Compared to the PRW distance, the IPRW distance performs better statistically but remains unfavorable in computational sense. Indeed, a large amount of projections from $\bs_{d, k}$ are necessary to approximate the IPRW distance. However, if the intrinsic dimension of data distribution is small, the required number of random projections is small; see~\citet{Nadjahi-2019-Asymptotic}.

Let $X_{1:n} = (X_1, \ldots, X_n)$ be independent and identically distributed samples according to the true measure $\mu_\star \in \PScr_q(\br^d)$. The empirical measure of $X_{1:n}$ is defined by $\widehat{\mu}_n \mydefn (1/n)\sum_{i=1}^n \delta_{X_i}$. It is known that $\widehat{\mu}_n \Rightarrow \mu_\star$ almost surely, and $\WCal_p(\widehat{\mu}_n, \mu_\star) \rightarrow 0$ almost surely since Wasserstein distances metrizes weak convergence~\citep[Theorem~6.9]{Villani-2008-Optimal} (note that $q \geq p \geq 1$). However, $\EE[\WCal_p(\widehat{\mu}_n, \mu_\star)] \simeq n^{-1/d}$ whenever $\mu$ is absolutely continuous with respect to Lebesgue measure and $d > 2p$~\citep{Dudley-1969-Speed, Fournier-2015-Rate, Weed-2019-Sharp} ($\simeq$ means ``equal to" with a constant independent of $n$). The convergence is slow when the dimension is high --- an instance of the well-known curse-of-dimensionality phenomenon. 

Due to the low-dimensional structure of the IPRW and PRW distances, the rate of IPRW and PRW distances is expected to be of $n^{-1/k}$ in the large-$n$ limit. Similar rates have been derived for $\EE[|\overline{\PWCal}_{k, p}(\widehat{\mu}_n, \widehat{\nu}_n) - \WCal_p(\mu, \nu)|]$ as a function of $n$ under a spiked transport model for both $\mu$ and $\nu$; see~\citet[Theorem~8]{Weed-2019-Estimation}. Their bound depends on problem dimension $d$ and requires $\mu$ and $\nu$ to satisfy the Talagrand transport inequality~\citep{Talagrand-1996-Transportation}. For the special case when $k=1$, the rate for the IPRW distance was studied in~\citep{Nadjahi-2020-Statistical} and the minimax confidence intervals were established in~\citet{Manole-2019-Minimax}. To our knowledge, there has been no other paper on the statistical properties of IPRW and PRW distances for $k \geq 2$.  

\paragraph{Parametric modeling and inference.} A statistical model is a family of distributions, $\MCal = \{\mu_\theta \in \PScr(\br^d) \mid \theta \in \Theta\}$, where $\Theta$ is the parameter space. A minimal set of the conditions of a proper family of distribution are: (i) $(\Theta, \|\cdot\|_\Theta)$ is a Polish space, (ii) $\Theta$ is $\sigma$-compact, i.e., it is the union of countably many compact subspaces, and (iii) parameters are identifiable, i.e., $\mu_\theta = \mu_{\theta'}$ implies $\theta = \theta'$ for all $\theta, \theta' \in \Theta$. Since the space $\PScr_p(\br^d)$ endowed with the distance $\WCal_p$ is a Polish space, we estimate model coefficients using minimum distance estimation (MDE)~\citep{Wolfowitz-1957-Minimum, Basu-2011-Statistical}, where the distance we consider here is PRW. The main reason why we do not choose IPRW in this setting is computational. The \textit{minimum project robust Wasserstein} (MPRW) estimator is defined as follows: 
\begin{equation}\label{def:MPRW}
\widehat{\theta}_n \mydefn \argmin_{\theta \in \Theta} \overline{\PWCal}_{p, k}(\widehat{\mu}_n, \mu_\theta). 
\end{equation}
Note that the probability density function of $\mu_\theta$ can be difficult to evaluate in practice, especially when $\mu_\theta$ is a generative model. Nevertheless, in various settings, even if the density is not available, one can generate samples $Z_{1:m}$ from $\mu_\theta$ and use them to approximate $\mu_\theta$. With this approximation, a natural alternative is the \textit{minimum expected projection robust Wasserstein} (MEPRW) estimator, which is defined as follows~\citep{Bernton-2019-Parameter, Nadjahi-2019-Asymptotic}: 
\begin{equation}\label{def:MEPRW}
\widehat{\theta}_{n, m} \mydefn \argmin_{\theta \in \Theta} \EE[\overline{\PWCal}_{p, k}(\widehat{\mu}_n, \widehat{\mu}_{\theta, m}) \mid X_{1:n}], 
\end{equation}
where $n$ is the number of samples from the data distribution $\mu_\star$, $m$ is the number of samples from the parametric distribution $\mu_\theta$, and $\widehat{\mu}_{\theta, m}$ is an empirical version of $\mu_\theta$ based on samples $Z_{1:m}$. 

Existing works have established asymptotic guarantees for minimal Wasserstein and sliced Wasserstein estimators~\citep{Bernton-2019-Parameter, Nadjahi-2019-Asymptotic}. Despite the similar proof paths, our results for the MPRW and MEPRW estimators are new and derived under weaker assumptions and more general settings than previous work; see Sections~\ref{subsec:estimators} and~\ref{subsec:distribution}.  

\section{Main Results on Projection Robust Optimal Transport Estimation}\label{sec:results}
Throughout this section, we assume $p \geq 1$ and $k \in [d] \triangleq \{1,2,\ldots,d\}$ unless stated otherwise. Focusing on the IPRW and PRW distances, we prove that they are lower semi-continuous and metrize weak convergence. Through a new sample complexity analysis, we derive the convergence rate of empirical measures under both distances as well as an improved rate for the PRW distance when $\mu_\star$ satisfies either a Bernstein tail condition or the Poincar\'{e} inequality. For the generative models with the PRW distance, we study the misspecified setting where the limit $\theta_\star$ is not necessarily the limit of the maximum likelihood estimator. We establish the asymptotic properties of the MPRW and MEPRW estimators and formulate a central limit theorem when $p=1$ and $k=1$.

\subsection{Topological properties}\label{subsec:topology}
We begin with the results on the relationship between the IPRW, PRW and Wasserstein distances. The following lemma demonstrates their equivalence in a topological sense. 
\begin{lemma}\label{lemma:PRW-equivalence}
The IPRW, PRW and Wasserstein distances are equivalent. That is, for any sequence of probability measures $\{\mu_i\}_{i \in \bn}$ and probability measure $\mu$ in $\PScr_p(\br^d)$, we have $\underline{\PWCal}_{p, k}(\mu_i, \mu) \rightarrow 0$ if and only if $\overline{\PWCal}_{p, k}(\mu_i, \mu) \rightarrow 0$ if and only if $\WCal_p(\mu_i, \mu) \rightarrow 0$. 
\end{lemma}
Lemma~\ref{lemma:PRW-equivalence} is a generalization of~\citet[Theorem~1]{Bayraktar-2019-Strong} where the projection dimension is $k=1$. By Lemma~\ref{lemma:PRW-equivalence} and~\citet[Theorem~6.9]{Villani-2008-Optimal}, we obtain the following result regarding the topology induced by the IPRW and PRW distances of order $p$. 
\begin{theorem}\label{Theorem:PRW-convergence}
The IPRW and PRW distances both metrize weak convergence. In other words, for any sequence of probability measures $\{\mu_i\}_{i \in \bn}$ and probability measure $\mu$ in $\PScr_p(\br^d)$, we have $\underline{\PWCal}_{p, k}(\mu_i, \mu) \rightarrow 0$ if and only if $\overline{\PWCal}_{p, k}(\mu_i, \mu) \rightarrow 0$ if and only if $\mu_i \Rightarrow \mu$.
\end{theorem}
Theorem~\ref{Theorem:PRW-convergence} generalizes~\citet[Theorem~6.9]{Villani-2008-Optimal} since the PRW distance is the Wasserstein distance when the projection dimension $k=d$. When $k=1$, Theorem~\ref{Theorem:PRW-convergence} recovers the results presented by~\citet{Bayraktar-2019-Strong} which implies that the SW and max-SW distances metrize weak convergence. It is worthy noting that this implication is stronger than~\citet[Theorem~1]{Nadjahi-2019-Asymptotic}, which only provides a one-sided argument. 
\begin{theorem}\label{Theorem:PRW-lsc}
The IPRW and PRW distances are both lower semi-continuous in the usual weak topology. In other words, if the sequences of probability measures $\{\mu_i\}_{i \in \bn}, \{\nu_i\}_{i \in \bn} \subseteq \PScr(\br^d)$ satisfy  $\mu_i \Rightarrow \mu$ and $\nu_i \Rightarrow \nu$ for probability measures $\mu, \nu \in \PScr(\br^d)$, then we have $\underline{\PWCal}_{p, k}(\mu, \nu) \leq \liminf_{i \rightarrow +\infty} \underline{\PWCal}_{p, k}(\mu_i, \nu_i)$ and $\overline{\PWCal}_{p, k}(\mu, \nu) \leq \liminf_{i \rightarrow +\infty} \overline{\PWCal}_{p, k}(\mu_i, \nu_i)$. 
\end{theorem}
The above theorem generalizes~\citet[Lemma~S6]{Nadjahi-2019-Asymptotic} and is pivotal to our asymptotic analysis for the MPRW and MEPRW estimators.

\subsection{Convergence and concentration of empirical measures}\label{subsec:convergence}
We provide the rate of the empirical measures under the IPRW and PRW distances of order $p$ with the projection dimension $k$. We present our main result on convergence rates in the following theorem. 
\begin{theorem}\label{Theorem:IPRW-main}
Let $\mu_\star \in \PScr_q(\br^d)$ and $M_q(\mu_\star) < +\infty$ for some $q \geq p \geq 1$. Then we have\footnote{$a \vee b = \max\{a, b\}$ and $a \wedge b = \min\{a, b\}$ here.}
\begin{equation*}
\EE[\underline{\PWCal}_{p, k}(\widehat{\mu}_n, \mu_\star)] \lesssim_{p, q} n^{-[\frac{1}{(2p)\vee k}\wedge(\frac{1}{p}-\frac{1}{q})]}(\log(n))^{\frac{\zeta_{p,q,k}}{p}},
\end{equation*}
where $\lesssim_{p, q}$ refers to ``less than" with a constant depending only on $(p, q)$ and 
\begin{equation*}
\zeta_{p, q, k} = \left\{\begin{array}{ll}
2 & \textnormal{if } k = q = 2p, \\
1 & \textnormal{if } (k \neq 2p \text{ and } q = \frac{kp}{k-p}) \textnormal{ or } (q > k = 2p), \\
0 & \textnormal{otherwise.}
\end{array}\right. 
\end{equation*}
\end{theorem}
\begin{remark}
Theorem~\ref{Theorem:IPRW-main} shows that our bound does not depend on $d$, while all bounds for the Wasserstein distance grow exponentially in $d$ when $d \geq 2p$~\citep[Theorem~3.1]{Lei-2020-Convergence}. This improvement shows that the PRW distance does not suffer from the curse of dimensionality while retaining flexibility via the choice of $k$. We are also aware of concurrent work~\citep{Nath-2020-Statistical} in which the sample complexity has no dependence on dimensionality. 
\end{remark}
\begin{definition}\label{def:proj-Bernstein}
We say $\mu \in \PScr(\br^d)$ satisfies a \emph{projection Bernstein tail condition} if there exist $\sigma, V > 0$ for all $E \in \bs_{d, k}$ and $X \sim E_{\#}^\star\mu$ such that $\EE[\|X\|^r] \leq (1/2)\sigma^2 r! V^{r-2}$ for all $r \geq 2$. 
\end{definition}
\begin{theorem}\label{Theorem:PRW-main-Bernstein}
Suppose $\mu_\star \in \PScr_q(\br^d)$ satisfies a projection Bernstein tail condition and assume the same setting as in Theorem~\ref{Theorem:IPRW-main}. For all $n \geq 1$, the following inequality holds true:
\begin{equation*}
\EE[\overline{\PWCal}_{p, k}(\widehat{\mu}_n, \mu_\star)] \ \lesssim_{p, q} \ n^{-[\frac{1}{(2p)\vee k}\wedge(\frac{1}{p}-\frac{1}{q})]}(\log(n))^{\frac{\zeta_{p,q,k}}{p}} + n^{\frac{1}{2}-\frac{1}{p}}\sqrt{dk\log(n)} + n^{-\frac{1}{p}}dk\log(n). 
\end{equation*}
\end{theorem}
\begin{definition}\label{def:proj-Poincare}
We say $\mu \in \PScr(\br^d)$ satisfies a \emph{projection Poincar\'{e} inequality} if there exists $M > 0$ for all $E \in \bs_{d, k}$ and $X \sim E_{\#}^\star\mu$ such that $\Var(f(X)) \leq M\EE[\|\nabla f(X)\|^2]$ for any $f: \br^d \rightarrow \br$ satisfying that $\EE[f(X)^2] < +\infty$ and $\EE[\|\nabla f(X)\|^2] < +\infty$. 
\end{definition}
\begin{theorem}\label{Theorem:PRW-main-Poincare}
Suppose $\mu_\star \in \PScr_q(\br^d)$ satisfies a projection Poincar\'{e} inequality and assume the same setting as in Theorem~\ref{Theorem:IPRW-main}. For all $n \geq 1$, the following inequality holds true:
\begin{equation*}
\EE[\overline{\PWCal}_{p, k}(\widehat{\mu}_n, \mu_\star)] \ \lesssim_{p, q} \ n^{-[\frac{1}{(2p)\vee k}\wedge(\frac{1}{p}-\frac{1}{q})]}(\log(n))^{\frac{\zeta_{p,q,k}}{p}} + n^{-\frac{1}{2 \vee p}}\sqrt{dk\log(n)} + n^{-\frac{1}{p}}dk\log(n).  
\end{equation*}
\end{theorem}
We present concentration results when $\mu_\star$ satisfies stronger conditions than Definition~\ref{def:proj-Bernstein} and~\ref{def:proj-Poincare}. 
\begin{definition}\label{def:Bernstein}
A measure $\mu \in \PScr(\br^d)$ satisfies a \emph{Bernstein tail condition} if there exists $\sigma, V > 0$ such that $\EE_{X \sim \mu}[\sup_{E \in \bs_{d, k}} \|E^\top X\|^r] \leq (1/2)\sigma^2 r! V^{r-2}$ for all $i =1, 2, \ldots, n$ and all $r \geq 2$. 
\end{definition}
\begin{theorem}\label{Theorem:concentration-bernstein}
If $\mu_\star \in \PScr(\br^d)$ satisfies a Bernstein tail condition then the following statement holds true for both $W = \underline{\PWCal}_{p, k}$ and $W = \overline{\PWCal}_{p, k}$:
\begin{equation*}
\PP(|W(\widehat{\mu}_n, \mu_\star) - \EE[W(\widehat{\mu}_n, \mu_\star)]| \geq t) \ \leq \ 2\exp\left(-\frac{t^2}{8\sigma^2 n^{1-2/p} + 4tVn^{-1/p}}\right). 
\end{equation*}
\end{theorem}
\begin{definition}\label{def:Poincare}
$\mu \in \PScr(\br^d)$ satisfies a \emph{Poincar\'{e} inequality} if there exists $M > 0$ for $X \sim \mu$ such that $\Var[f(X)] \leq M\EE[\|\nabla f(X)\|^2]$ for any $f$ satisfying $\EE[f(X)^2] < +\infty$ and $\EE[\|\nabla f(X)\|^2] < +\infty$. 
\end{definition}
\begin{theorem}\label{Theorem:concentration-poincare}
If $\mu_\star \in \PScr(\br^d)$ satisfies Poincar\'{e} inequality then the following statement holds true for both $W = \underline{\PWCal}_{p, k}$ and $W = \overline{\PWCal}_{p, k}$:
\begin{equation*}
\PP(|W(\widehat{\mu}_n, \mu_\star) - \EE[W(\widehat{\mu}_n, \mu_\star)]| \geq t) \ \leq \ 2\exp(-K^{-1}\min\{n^{\frac{1}{p}}t, \ n^{\frac{2}{2\vee p}} t^2\}),   
\end{equation*}
where $K > 0$ only depends on $M$ (cf. Definition~\ref{def:Poincare}). 
\end{theorem}

\paragraph{Discussions.} We demonstrate that the Bernstein-type tail conditions in Definition~\ref{def:proj-Bernstein} and~\ref{def:Bernstein} are not strong enough to give an effective bound for all $p \geq 1$. The similar results for the Wasserstein distance have been recently derived by~\citet{Lei-2020-Convergence} and recognized as the standard limitation for the Bernstein-type tail conditions. This is also the motivation which drives us to consider a Poincar\'{e} inequality. 

For Theorem~\ref{Theorem:PRW-main-Bernstein} and~\ref{Theorem:PRW-main-Poincare}, the first term matches that in Theorem~\ref{Theorem:IPRW-main} while the extra two terms come from bounding the gap $\EE[\sup_{E \in \bs_{d, k}} (\WCal_p(E_{\#}^\star\widehat{\mu}_n, E_{\#}^\star\mu_\star) - \EE[\WCal_p(E_{\#}^\star\widehat{\mu}_n, E_{\#}^\star\mu_\star)])]$. Compared with~\citet[Theorem~8]{Weed-2019-Estimation}, where $\mu_\star$ satisfies the Talagrand transport inequality, our conditions are weaker but our rate matches their $n^{-1/k} + n^{-1/2}\sqrt{dk\log(n)}$ rate in the large-$n$ limit when $p=1$. For Theorem~\ref{Theorem:concentration-bernstein} and~\ref{Theorem:concentration-poincare}, the latter bound is better than the former bound when $p>1$. Moreover, the tail condition in Definition~\ref{def:Bernstein} is stronger than that in Definition~\ref{def:proj-Bernstein} yet weaker than the standard Bernstein tail condition where $X \sim \mu$ inside the expectation without a $\sup$; see~\citet{Wainwright-2019-High}. The Poincar\'{e} inequality is weaker than the log-Sobolev inequality and is satisfied by various exponential measures and the measures induced by Markov processes~\citep{Ledoux-1999-Concentration}. Intutively, These two conditions handle subexponential tail while Talagrand inequality in~\citet{Weed-2019-Estimation} addresses subgaussian tail; see~\citet{Ledoux-1999-Concentration} and~\citet{Talagrand-1996-Transportation} for the details.

\subsection{Properties of MPRW and MEPRW estimators}\label{subsec:estimators}
We derive the asymptotic properties of the MPRW and MEPRW estimators under model misspecification, which is common in practice. Our setting is more general than that considered in~\citep{Nadjahi-2019-Asymptotic} and our results support the applications in real-world scenario better. Specifically, while~\citet{Nadjahi-2019-Asymptotic} focused on the well-specified setting, the statistical models can be misspecified in many real-world applications. We also use the Wasserstein distance in Assumptions~\ref{A1} and~\ref{A4} since these assumptions have been shown valid for many real-world application problems~\citep{Bernton-2019-Parameter}.     
\begin{assumption}\label{A1}
There exists a probability measure $\mu_\star \in \PScr(\br^d)$ such that the data-generating process satisfies that $\lim_{n \rightarrow +\infty} \WCal_p(\widehat{\mu}_n, \mu_\star) = 0$ almost surely.   
\end{assumption}
\begin{assumption}\label{A2}
The map $\theta \mapsto \mu_\theta$ is continuous: $\|\theta_n - \theta\|_\Theta \rightarrow 0$ implies $\mu_{\theta_n} \Rightarrow \mu_\theta$. 
\end{assumption}
\begin{assumption}\label{A3}
There exists a constant $\tau > 0$ such that the set $\Theta_\star(\tau) \subseteq \Theta$ is bounded where $\Theta_\star(\tau) = \{\theta \in \Theta: \overline{\PWCal}_{p, k}(\mu_\star, \mu_\theta) \leq \inf_{\theta \in \Theta} \overline{\PWCal}_{p, k}(\mu_\star, \mu_\theta) + \tau\}$. 
\end{assumption}
\begin{theorem}\label{Theorem:consistency-MPRW}
Under Assumption~\ref{A1}-\ref{A3}, there exists a sample space $\Omega$ with $\PP(\Omega)=1$ such that, for all $\omega \in \Omega$, 
\begin{equation*}
\lim_{n \rightarrow +\infty} \inf_{\theta \in \Theta} \overline{\PWCal}_{p, k}(\widehat{\mu}_n(\omega), \mu_\theta) = \inf_{\theta \in \Theta} \overline{\PWCal}_{p, k}(\mu_\star, \mu_\theta), 
\end{equation*}
and 
\begin{equation*}
\limsup_{n \rightarrow +\infty} \argmin_{\theta \in \Theta} \overline{\PWCal}_{p, k}(\widehat{\mu}_n(\omega), \mu_\theta) \subseteq \argmin_{\theta \in \Theta} \overline{\PWCal}_{p, k}(\mu_\star, \mu_\theta). 
\end{equation*}
In addition, $\argmin_{\theta \in \Theta} \overline{\PWCal}_{p, k}(\widehat{\mu}_n(\omega), \mu_\theta) \neq \emptyset$ for all $n \geq n(\omega)$ with some $n(\omega) > 0$. 
\end{theorem}
\begin{assumption}\label{A4}
If $\|\theta_n - \theta\|_\Theta \rightarrow 0$, then $\EE[\WCal_p(\widehat{\mu}_{\theta_n, n}, \mu_{\theta_n}) | X_{1:n}] \rightarrow 0$. 
\end{assumption}
In the next result, we present an analogous version of Theorem~\ref{Theorem:consistency-MPRW} for the MEPRW estimator as $\min\{n, m\} \rightarrow +\infty$. For the simplicity, we set $m \mydefn m(n)$ such that $m(n) \rightarrow +\infty$ as $n \rightarrow +\infty$. 
\begin{theorem}\label{Theorem:consistency-MEPRW}
Under Assumption~\ref{A1}-\ref{A4}, there exists a sample space $\Omega$ with $\PP(\Omega)=1$ such that, for all $\omega \in \Omega$, 
\begin{equation*}
\lim_{n \rightarrow +\infty}\inf_{\theta \in \Theta} \EE[\overline{\PWCal}_{p, k}(\widehat{\mu}_n(\omega), \widehat{\mu}_{\theta, m(n)})|X_{1:n}] = \inf_{\theta \in \Theta} \overline{\PWCal}_{p, k}(\mu_\star, \mu_\theta), 
\end{equation*}
and 
\begin{equation*}
\limsup_{n \rightarrow +\infty} \argmin_{\theta \in \Theta} \EE[\overline{\PWCal}_{p, k}(\widehat{\mu}_n(\omega), \widehat{\mu}_{\theta, m(n)}) \mid X_{1:n}] \subseteq \argmin_{\theta \in \Theta} \overline{\PWCal}_{p, k}(\mu_\star, \mu_\theta). 
\end{equation*}
In addition, $\argmin_{\theta \in \Theta} \EE[\overline{\PWCal}_{p, k}(\widehat{\mu}_n(\omega), \widehat{\mu}_{\theta, m(n)})|X_{1:n}] \neq \emptyset$ for $n \geq n(\omega)$ with some $n(\omega) > 0$.  
\end{theorem}
\begin{assumption}\label{A5}
There exists a constant $\tau > 0$ such that the set $\Theta_n(\tau) \subseteq \Theta$ is bounded where $\Theta_n(\tau) = \{\theta \in \Theta: \overline{\PWCal}_{p, k}(\widehat{\mu}_n, \mu_\theta) \leq \inf_{\theta \in \Theta} \overline{\PWCal}_{p, k}(\widehat{\mu}_n, \mu_\theta) + \tau\}$. 
\end{assumption}
\begin{theorem}\label{Theorem:convergence-MEPRW-MPRW}
Under Assumption~\ref{A2},~\ref{A4} and~\ref{A5}, the following statement holds true, 
\begin{equation*}
\lim_{m \rightarrow +\infty} \inf_{\theta \in \Theta} \EE[\overline{\PWCal}_{p, k}(\widehat{\mu}_n, \widehat{\mu}_{\theta, m}) \mid X_{1:n}] = \inf_{\theta \in \Theta} \overline{\PWCal}_{p, k}(\widehat{\mu}_n, \mu_\theta), 
\end{equation*}
and 
\begin{equation*}
\limsup_{m \rightarrow +\infty} \argmin_{\theta \in \Theta} \EE[\overline{\PWCal}_{p, k}(\widehat{\mu}_n, \widehat{\mu}_{\theta, m}) \mid X_{1:n}] \subseteq \argmin_{\theta \in \Theta} \overline{\PWCal}_{p, k}(\widehat{\mu}_n, \mu_\theta). 
\end{equation*}
In addition, $\argmin_{\theta \in \Theta} \EE[\overline{\PWCal}_{p, k}(\widehat{\mu}_n, \widehat{\mu}_{\theta, m})|X_{1:n}] \neq \emptyset$ for $m \geq m_n$ with some $m_n> 0$. 
\end{theorem}
To this end, the MPRW and MEPRW estimators both asymptotically converge to $\theta_\star \in \Theta$, which is a minimizer of $\theta \rightarrow \overline{\PWCal}_{p, k}(\mu_\star, \mu_\theta)$, assuming its existence. Moreover, $\theta_\star$ is not the limit of maximum likelihood estimator and satisfies $\mu_{\theta_\star} = \mu_\star$ in a well-specified setting. Our consistency results support the success of generative modelling using the max-SW distance.  

\subsection{Rate of convergence and asymptotic distribution}\label{subsec:distribution}
We investigate the asymptotic distribution of the MPRW estimator under model misspecification and establish the rate of convergence when $k = p = 1$. For any $u \in \bs^{d-1}$ and $t \in \br$, we define 
\begin{eqnarray*}
F_\theta(u, t) & = & \int_{\br^d} \one_{(-\infty, t]}(\langle u, x\rangle) \; d\mu_\theta(x), \\ 
\widehat{F}_n(u, t) & = & (1/n)|\{i \in [n]: \langle u, X_i\rangle \leq t\}|. 
\end{eqnarray*}
The functions $F_\theta(u, \cdot)$ and $\widehat{F}_n(u, \cdot)$ are the cumulative distribution functions of $u^\star_{\#}\mu_\theta$ and $u^\star_{\#}\widehat{\mu}_n$ where $u \in \bs^{d-1}$ is a unit vector. Let $L(\bs^{d-1} \times \br)$ be the class of functions on $\bs^{d-1} \times \br$ such that $f(\cdot, t)$ is continuous and $f(u, \cdot)$ is absolutely integrable, with the norm $\|f\|_L = \sup_{u \in \bs^{d-1}} \int_\br |f(u, t)| \; dt$.  
\begin{assumption}\label{A7}
There exists a measurable function $D_\star: \bs^{d-1} \times \br \rightarrow \br^{d_\theta}$ such that $\|F_\theta(u, t) - F_{\theta_\star}(u, t) - \langle\theta - \theta_\star, D_\star(u, t)\rangle\|_L = o(\|\theta - \theta_\star\|_\Theta)$. 
\end{assumption}
\begin{assumption}\label{A8}
There exists a random element $G_\star: \bs^{d-1} \times \br \mapsto \br$ such that the stochastic process $\sqrt{n}(\widehat{F}_n - F_\star)$ converges weakly in $L(\bs^{d-1} \times \br)$ to $G_\star$\footnote{As pointed by~\citet{Nadjahi-2019-Asymptotic}, one can prove that Assumption~\ref{A8} holds in general by extending~\citep[Proposition~3.5]{Dede-2009-Empirical} and~\citep[Theorem~2.1(a)]{Del-1999-Central} with some mild conditions on the tails of $u_{\#}^\star\mu_\star$. Using the same argument, this extension can also be done for $\|\cdot\|_L$ in our paper.}. 
\end{assumption}
\begin{assumption}\label{A9}
There exists a neighborhood $\NCal$ of $\theta_\star \in \Theta$ and a positive constant $c_\star$ such that $\overline{\PWCal}_{1, 1}(\mu_\theta, \mu_\star) \ \geq \ \overline{\PWCal}_{1, 1}(\mu_{\theta_\star}, \mu_\star) + c_\star\|\theta - \theta_\star\|_\Theta$ for all $\theta \in \NCal$.
\end{assumption}
\begin{remark}\label{Remark:assumption_CLT}
Assumption~\ref{A7} is strictly weaker than a norm-differentiation condition where $D_\star$ has to be nonsingular. Assumption~\ref{A8} permits model misspecification where there is no $\theta_\star \in \Theta$ such that $F_{\theta_\star} = F_\star$ and thus is more general than~\citet[A8]{Nadjahi-2019-Asymptotic}. Assumption~\ref{A9} accounts for local strong identifiability for the model $\mu_\theta$ around $\theta_\star$ and is necessary for the fast rate of $n^{-1/2}$ under model misspecification. (\citet{Bernton-2019-Parameter} assumes the analogous condition for the Wasserstein distance. However, their analysis depends on a much stronger version with $\NCal = \Theta$.) Thanks to Assumption~\ref{A9}, we do not require the condition that the parameters are weakly separable in the PRW sense.
\end{remark} 
\begin{remark}\label{Remark:well_specified}
In well-specified setting where there exists $\theta_\star \in \Theta$ such that $F_\star = F_{\theta_\star}$, it is straightforward to derive the norm-differentiation condition from Assumption~\ref{A7} and~\ref{A9}. This is not true, however, under model misspecification. Moreover, there are minor technical issues in the proof of~\citet[Theorem~B.8]{Bernton-2019-Parameter}; see Appendix~\ref{appendix:minor}. Fixing them would be straightforward but require additional assumptions. Fortunately, we can overcome this gap using some new techniques. Thus, with some refinement, our results can be interpreted as an improvement of~\citet{Bernton-2019-Parameter} with fewer assumptions.   
\end{remark}
\begin{figure*}[!t]
\centering
\includegraphics[width=0.32\textwidth]{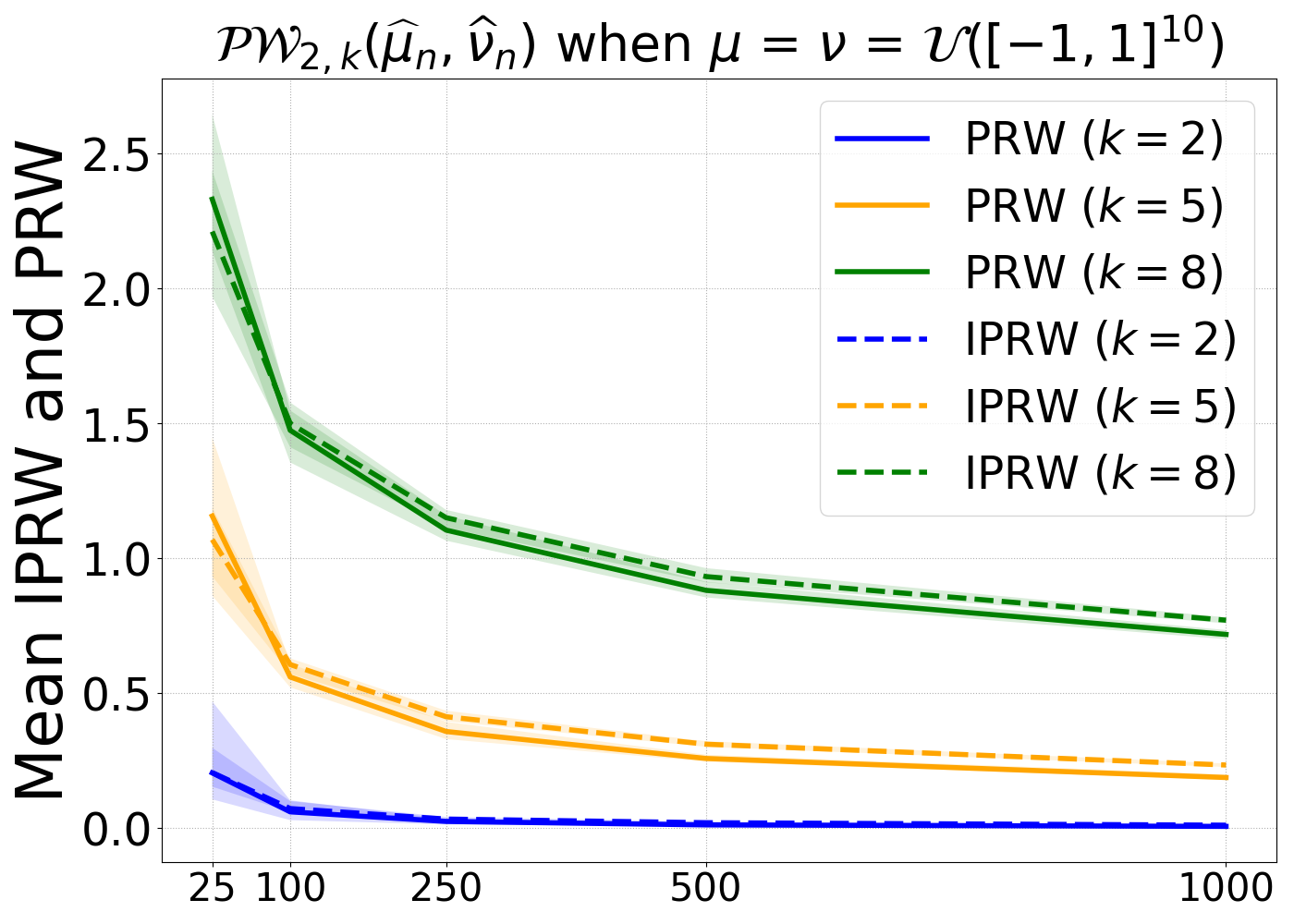}
\includegraphics[width=0.32\textwidth]{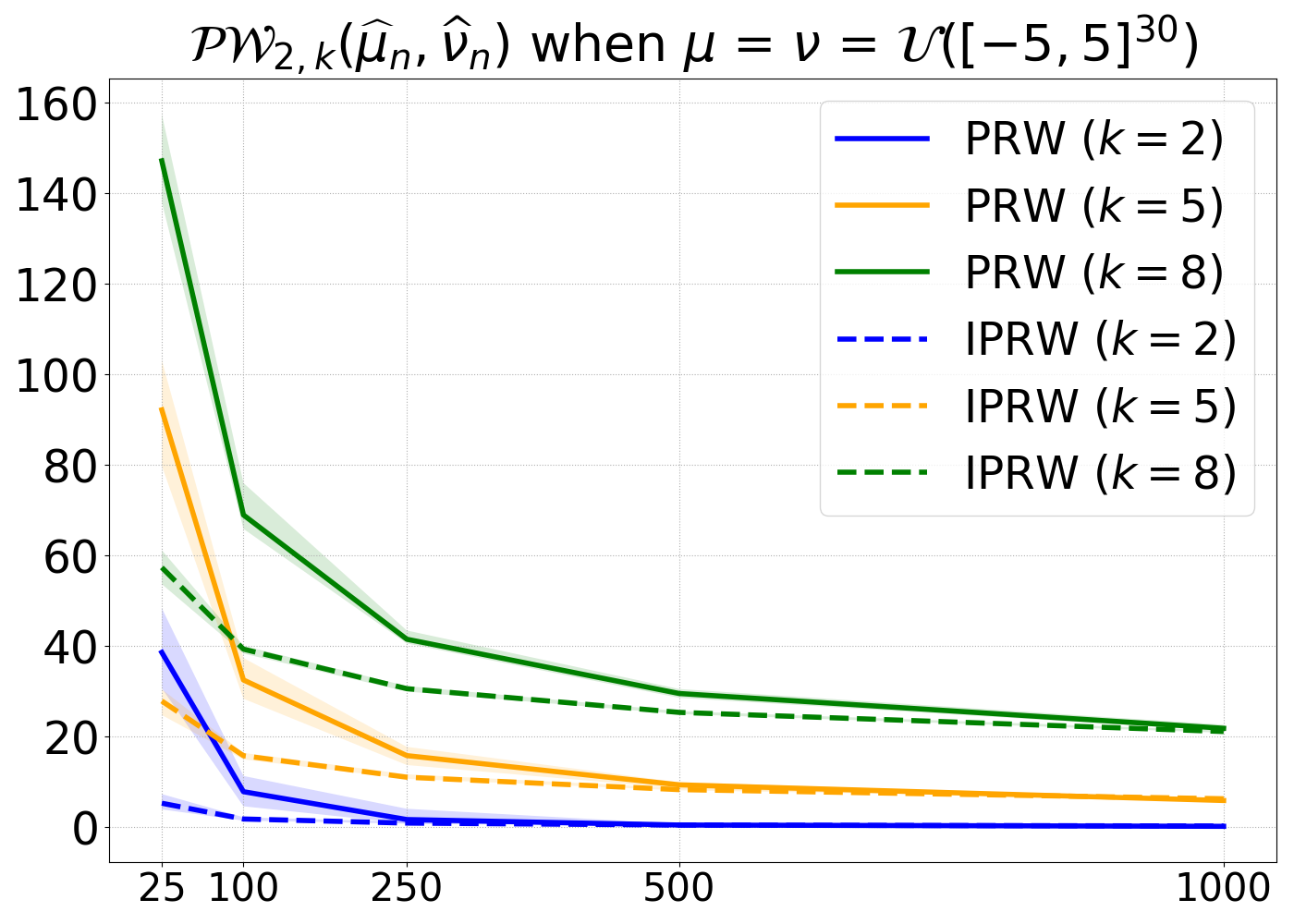}
\includegraphics[width=0.32\textwidth]{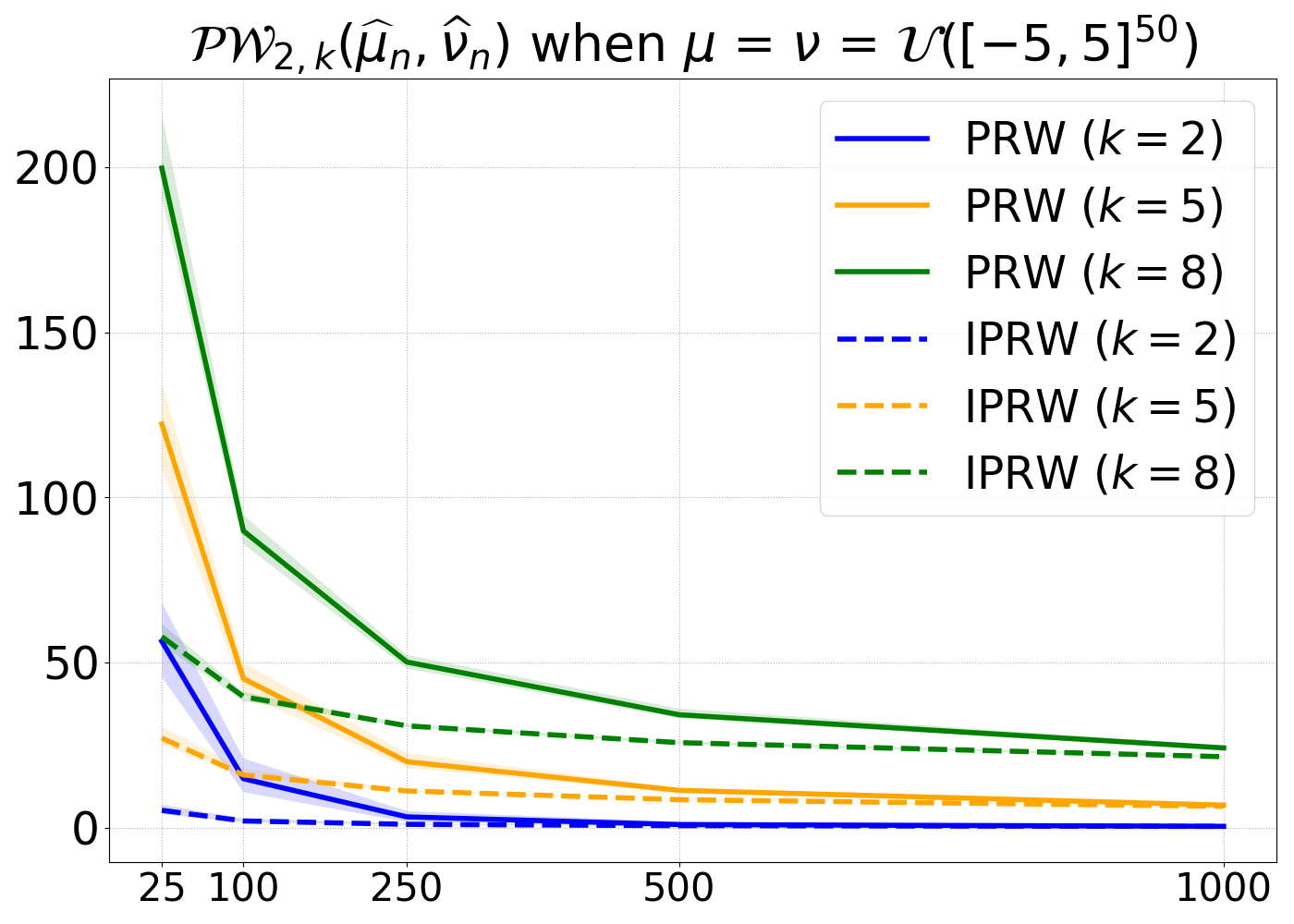}
\includegraphics[width=0.32\textwidth]{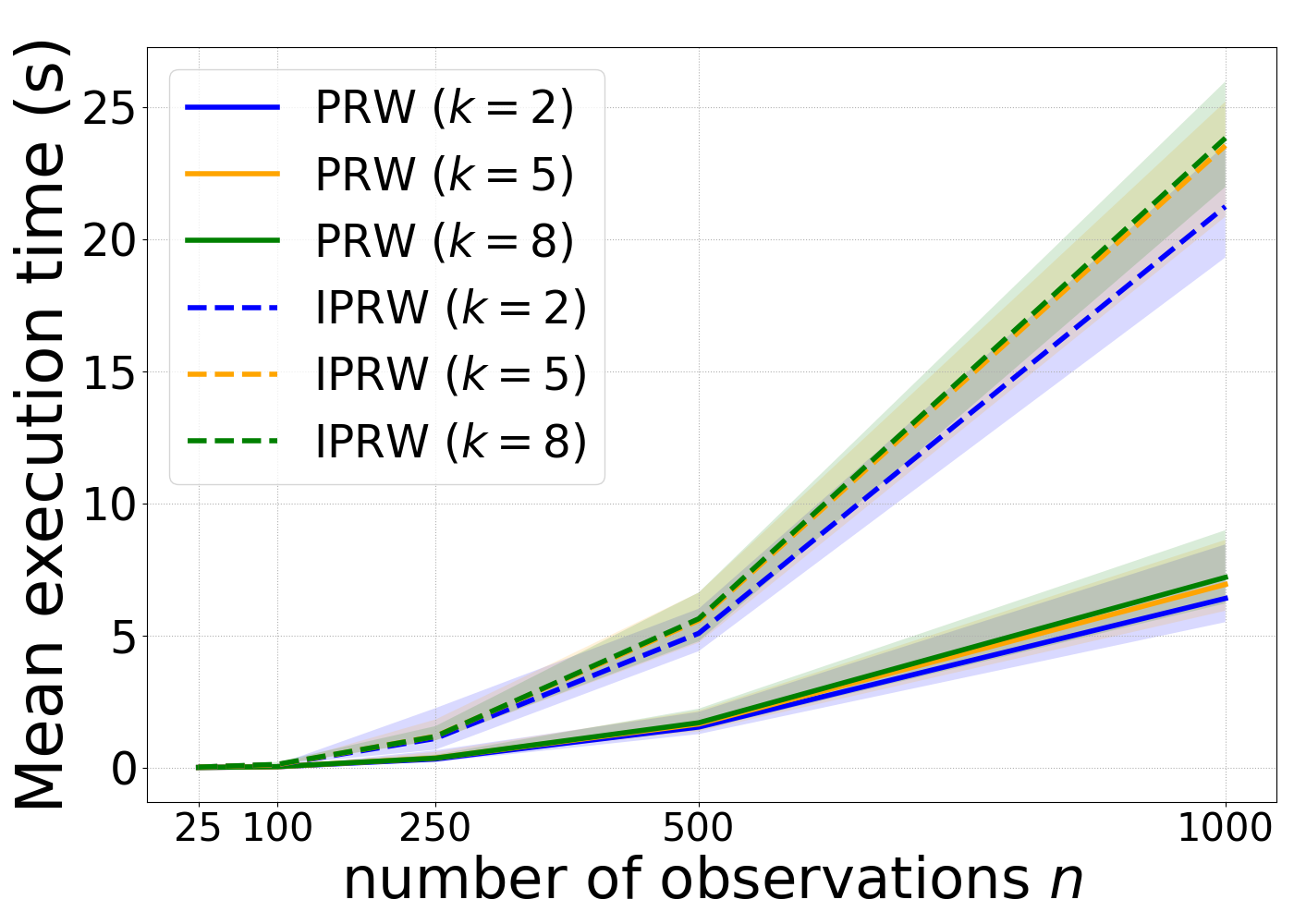}
\includegraphics[width=0.32\textwidth]{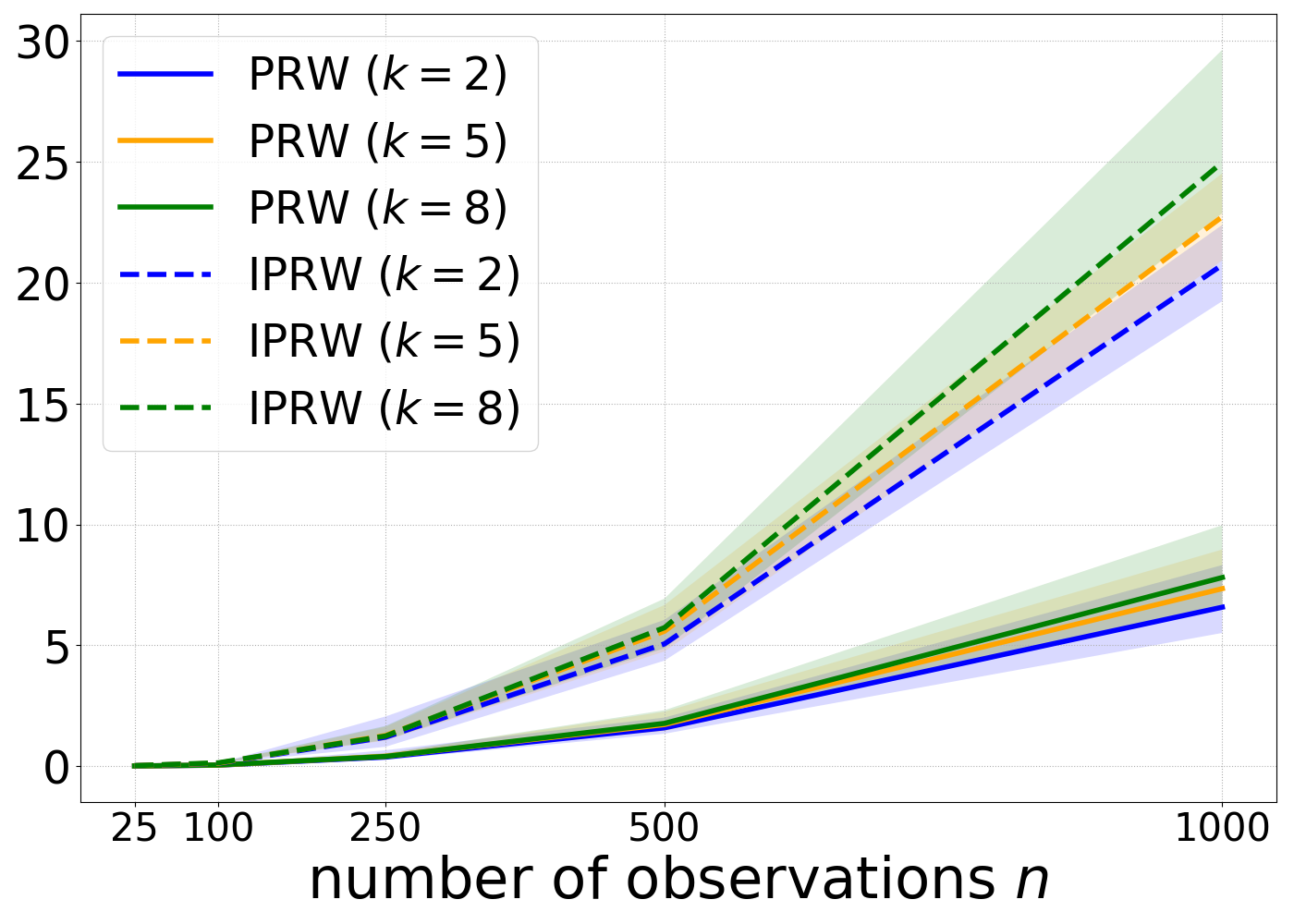}
\includegraphics[width=0.32\textwidth]{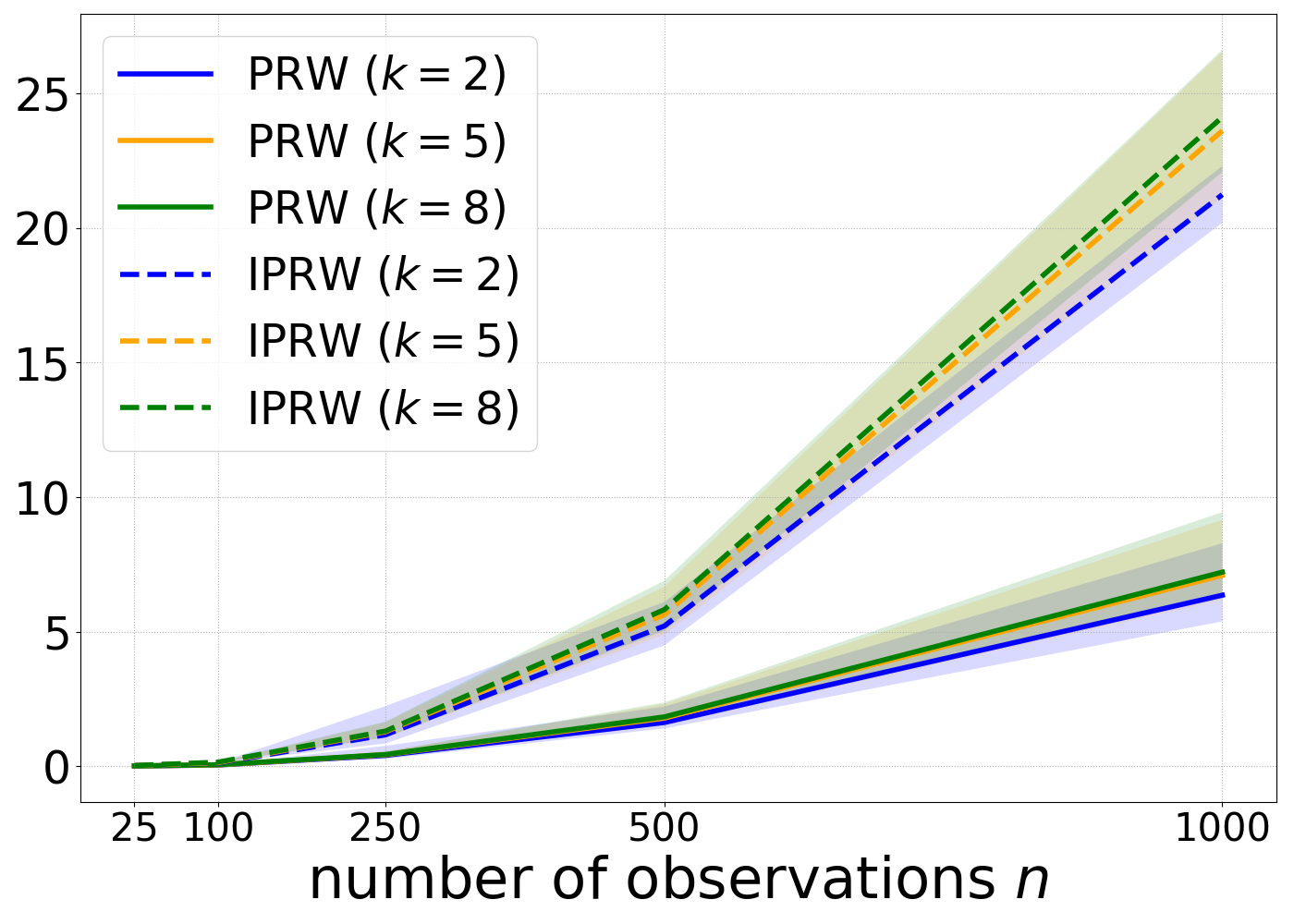}
\caption{\small{Mean values (Top) and mean computational time (Bottom) of the IPRW and PRW distances of order 2 between empirical measures $\widehat{\mu}_n$ and $\widehat{\nu}_n$ as the number of points $n$ varies. Results are averaged over 100 runs.}}\label{fig:exp1}\vspace*{-1em}
\end{figure*}
To study the asymptotic distributions in the misspecified setting, we employ definitions from~\citet[Section~7]{Pollard-1980-Minimum}. (Note, however, that our proof technique is different from~\citet{Pollard-1980-Minimum}, which depends on the nonsingularity of $D_\star$ and requires $\mu_\star = \mu_{\theta_\star}$ for some $\theta_\star$ in the interior of $\Theta$.)
\begin{definition}[Hausdorff metric]
Let $\SCal$ be the class of convex and compact sets in $L(\bs^{d-1} \times \br)$ equipped with $\|\cdot\|_L$. The Hausdorff metric on $\SCal$ is defined by $d_H(S_1, S_2) = \inf\{\delta > 0: S_1 \subseteq S_2^\delta, S_2 \subseteq S_1^\delta\}$, where $S^\delta = \cup_{x \in S} \{z \in L(\bs^{d-1} \times \br): \|z - x\|_L \leq \delta\}$. 
\end{definition}
\begin{definition}[Approximate MPRW estimators] The set of approximate MPRW estimators is defined by $M_n = \{\theta \in \Theta: \overline{\PWCal}_{1, 1}(\widehat{\mu}_n, \mu_\theta) \leq \inf_{\theta' \in \Theta} \overline{\PWCal}_{1, 1}(\widehat{\mu}_n, \mu_{\theta'}) + \eta_n / \sqrt{n}\}$, where $\eta_n > 0$ such that $\PP(\eta_n \rightarrow 0) = 1$ and $M_n$ is nonempty. 
\end{definition}
\begin{theorem}\label{Theorem:general}
Suppose Assumption~\ref{A1}-\ref{A3} and~\ref{A7}-\ref{A9} hold for some $\theta_\star$ in the interior of $\Theta$ and let $G_n = \sqrt{n}(\widehat{F}_n - F_{\theta_\star})$ and $G_n^\star = G_\star + \sqrt{n}(F_\star - F_{\theta_\star})$. We also define $K(x, \beta) = \{\theta \in \NCal_1: \|x - \sqrt{n}\langle\theta - \theta_\star, D_{\theta_\star}\rangle\|_L \leq \inf_{\theta' \in \NCal_1} \|x - \sqrt{n}\langle \theta'-\theta_\star, D_{\theta_\star}\rangle\|_L + \beta\}$ where
\begin{equation*}
\NCal_1 = \left\{\theta \in \NCal: \frac{\|F_\theta - F_{\theta_\star} - \langle\theta - \theta_\star, D_\star\rangle\|_L}{\|\theta - \theta_\star\|_\Theta} \leq \frac{c_\star}{2}\right\}.  
\end{equation*}
Then there exists a sequence satisfying $\lim_{n \rightarrow +\infty} \beta_n = 0$ such that\footnote{$\PP_\star$ denotes the (inner) probability; see~\citet{Pollard-1980-Minimum} for details.} $\PP_\star(M_n \subseteq K(G_n, \beta_n)) \rightarrow 1$ as $n \rightarrow +\infty$. For any $\epsilon > 0$, we have $\PP(d_H(K(G_n^\star, 0), K(G_n, \beta_n)) < \epsilon) \rightarrow 1$ as $n \rightarrow +\infty$. 
\end{theorem}
Theorem~\ref{Theorem:general} provides the theoretical guarantee for statistical inference with the max-SW distance under model misspecification. Indeed, since $K(G_n^\star, 0) = \argmin_{\theta \in \NCal_1} \|G_\star + \sqrt{n}(F_\star - F_{\theta_\star} - \langle \theta - \theta_\star, D_{\theta_\star}\rangle)\|_L$, the results indicate that the distributional limit of the approximate MPRW estimator set is close to the limit of the sets $\argmin_{\theta \in \NCal_1} \|G_\star + \sqrt{n}(F_\star - F_{\theta_\star} - \langle \theta - \theta_\star, D_{\theta_\star}\rangle)\|_L$ in the Hausdorff metric. Note that $d>1$ is allowed but we need $k=1$. This is necessary for our techniques since the current analysis heavily depends on the explicit form of PRW using cumulative distribution functions. Deriving CLT when $k>1$ is important but out of the scope of this paper.
\begin{remark}
In the well-specified setting, Assumption~\ref{A9} can be replaced by Assumption~\ref{A6}-\ref{A10}. Under certain conditions, we derive the CLT (cf. Theorem~\ref{Theorem:well-specificied}) which is analogous to~\citet[Theorem~6]{Nadjahi-2019-Asymptotic} for the minimum sliced Wasserstein estimators. We refer to Theorem~\ref{Theorem:well-specificied} in Appendix~\ref{appendix:A} for a simplified version in well-specified setting.
\end{remark}
\paragraph{Discussions.} We make some additional remarks on the relationship between our work and the existing works by~\citet{Bernton-2019-Parameter} and~\citet{Nadjahi-2019-Asymptotic}. Since PRW is a type of Wasserstein, the consistency proof roadmap is essentially similar to that in~\citet{Bernton-2019-Parameter} and~\citet{Nadjahi-2019-Asymptotic}. However, we remark that (i) the sample complexity bounds of PRW are new; (ii) the techniques for CLT in misspecified settings did not appear in~\citet{Nadjahi-2019-Asymptotic} and complete the analysis in~\citet{Bernton-2019-Parameter}. Remark~\ref{Remark:assumption_CLT} states that our Assumption~\ref{A9} is weaker than that is used in~\citet{Bernton-2019-Parameter}. In particular, $\NCal$ is the neighborhood defined in Assumption~\ref{A9} and accounts for a \textbf{local} strong identifiability. In contrast,~\citet{Bernton-2019-Parameter} requires a \textbf{global} strong identifiability ($\NCal=\Theta$). Remark~\ref{Remark:well_specified} states that our setting is more general than the well-specified setting which is discussed by~\citet{Nadjahi-2019-Asymptotic} from a technical point of view.

\begin{figure*}[!t]
\centering
\begin{subfigure}{0.32\textwidth}
\includegraphics[width=1.0\textwidth]{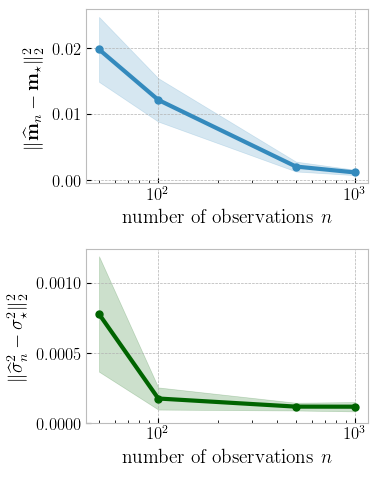}
\caption{\small{MPRW vs. $n$}}
\end{subfigure}
\begin{subfigure}{0.32\textwidth}
\includegraphics[width=1.0\textwidth]{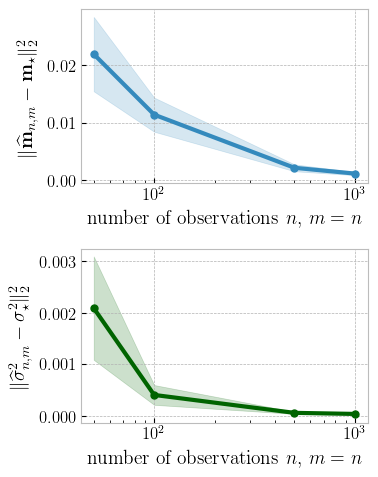}
\caption{\small{MEPRW vs. $n=m$}}
\end{subfigure}
\begin{subfigure}{0.32\textwidth}
\includegraphics[width=1.0\textwidth]{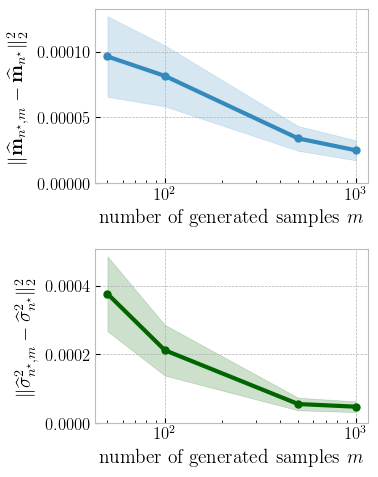}
\caption{\small{MEPRW with $n=2000$ vs. $m$}}
\end{subfigure}
\caption{\small{Minimal PRW and expected PRW estimations using Gaussian models and $n$ samples from the mixture of 8 Gaussian distributions. Results are averaged over 100 runs and shaded areas represent standard deviation.}}\label{fig:exp2_consistent}\vspace*{-1em}
\end{figure*}

\section{Experiments}
We empirically validate our theoretical findings through several experiments on synthetic and real data. Given the space limit, we present the experimental setup in Appendix~\ref{appendix:setup} and explain an optimization algorithm for computing the PRW distance and estimators in Appendix~\ref{appendix:computational}. We defer the additional results on other dataset to Appendix~\ref{appendix:results}.  

Let $\mu = \nu = \UCal([-v, v]^d)$ be an uniform distribution over a hypercube and we study the convergence and computation of $\underline{\PWCal}_{2, k}(\widehat{\mu}_n, \widehat{\nu}_n)$ and $\overline{\PWCal}_{2, k}(\widehat{\mu}_n, \widehat{\nu}_n)$ for $n \in \{20, 100, 250, 500, 1000\}$. Figure~\ref{fig:exp1} presents average distances and computational times for $(d, v) \in \{(10, 1), (30, 5), (50, 5)\}$, where the shaded areas show the max-min values over 100 runs. First, the IPRW distance is smaller than the PRW distance for small $n$ when $d$ and $v$ are large. This confirms that the IPRW distance is independent of $d$ (cf. Theorem~\ref{Theorem:IPRW-main}). Second, the PRW distance nearly matches the IPRW distance when $n$ is large. This confirms Theorem~\ref{Theorem:PRW-main-Poincare} since the uniform distribution with its bounded domain satisfies the Poincar\'{e} inequality. Finally, the computation of the PRW distance is faster than that of the IPRW distance.

\begin{wrapfigure}{r}{0.5\textwidth}
\includegraphics[width=0.45\textwidth]{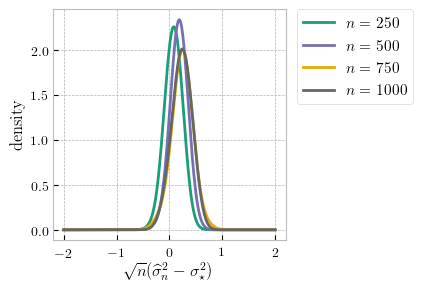}
\vspace*{-1em}\caption{Probability density of estimation of centered and rescaled $\widehat{\sigma}_n$ on the Gaussian model.}\label{fig:exp2_clt}\vspace*{-1em}
\end{wrapfigure}
Consider the parametric inference using Gaussian models $\MCal = \{\NCal(\textbf{m}, \sigma^2\textbf{I}):\textbf{m}{\in}\br^2, \sigma^2{>}0\}$ and a collection of i.i.d. observations generated from a mixture of 8 Gaussian distributions in $\br^2$. This simple setting is useful since the closed-form expression of Gaussian density makes the computation of the MPRW estimator of order 1 tractable. Following the setup in~\citet[Section~4]{Nadjahi-2019-Asymptotic}, we illustrate the consistency of the MPRW and MEPRW estimators of order 1 and the convergence of MEPRW estimator of order 1 to MPRW estimator of order 1. Results are shown in Figure~\ref{fig:exp2_consistent}; they are consistent with Theorem~\ref{Theorem:consistency-MPRW},~\ref{Theorem:consistency-MEPRW} and~\ref{Theorem:convergence-MEPRW-MPRW}, where $\textbf{m}_\star = \widehat{\textbf{m}}_{10^5}$. Despite the model misspecification, our estimators still converge as the number of observations increases and the MEPRW estimator converges to the MPRW estimator as we generate more samples. We also verify our central limit theorem by estimating the density of $\widehat{\sigma}_n^2$ with a kernel density estimator\footnote{The approach we apply here is the same as used by~\citet{Nadjahi-2019-Asymptotic}.} over $100$ runs. Figure~\ref{fig:exp2_clt} shows the distribution centered and rescaled by $\sqrt{n}$ for each $n$, where $\sigma_\star^2 = \widehat{\sigma}_{10^5}^2$, and it confirms the convergence rate we derived in Theorem~\ref{Theorem:general}; see Appendix~\ref{appendix:results} for the case with 12 or 25 distributions.

\begin{wrapfigure}{r}{0.5\textwidth}
\vspace{-1.5em}\includegraphics[width=0.4\textwidth]{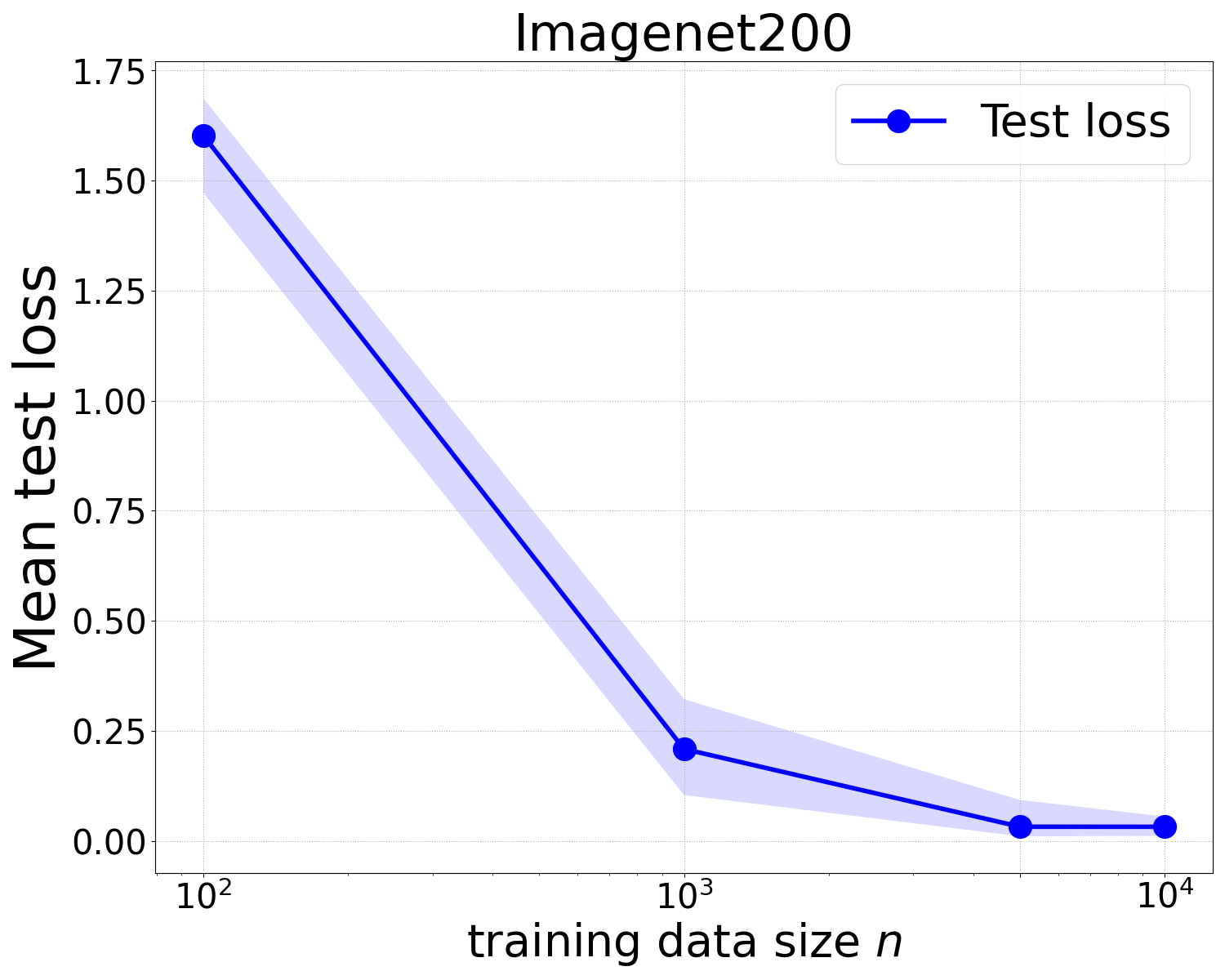}
\vspace{-1em}\caption{Mean test loss for different value of $(n, m)$ on \textsc{ImageNet200}.}\label{fig:exp3_imagenet}\vspace{-1em}
\end{wrapfigure}

We conduct experiments on image generation using the PRW generator of order 2, as an alternative to the SW generator~\citep{Deshpande-2018-Generative}. Here we focus on the case of $k=1$, where the PRW generator is exactly max-SW generator. We train the neural networks (NNs) with $(n, m) \in \{(100, 20), (1000, 40), (5000, 60), (10000, 100)\}$ where $n$ is the number of training samples and $m$ is the number of generated samples. We compare their testing losses to that of a NN trained using $n=10^5$ (i.e. whole training dataset) and $m=200$. All testing losses are evaluated using the trained models on the the testing dataset ($n=10^4$) with $m=250$ generated samples. Figure~\ref{fig:exp3_imagenet} presents the mean testing loss on \textsc{ImageNet200} over 10 runs, where the shaded areas show the max-min values over the runs. 

\paragraph{Discussions.} First, PRW has better discriminative power than max-SW or SW since it considers high-order summaries and extract more geometric information from two high-dimensional distributions, in order to distinguish them better; see~\citet{Paty-2019-Subspace} for the details. Moreover, we have presented in Figure~\ref{fig:exp1} (top row) and Figure~\ref{fig:exp1_appendix} (top row) that the PRW/IPRW value increase as $k$ increases. Thus, PRW/IPRW based on larger $k$-dimensional projections have better discriminative power.

Second, the IPRW computation generally requires many random projections and is thus more time-consuming than PRW for a desired accuracy when $k=1$; see~\citet[Page 4]{Kolouri-2019-Generalized}. Fortunately, it may require much fewer for certain application problems when the intrinsic dimension of data distribution is small, and is easily amenable to parallel computation. Thus, IPRW can serve as a practical alternative to PRW. Moreover, the reported PRW and IPRW values in Figure~\ref{fig:exp1} and Figure~\ref{fig:exp1_appendix} (appendix) are computed by using 30 iterations for PRW and 100 projections for IPRW. Therefore, the statistical/simulation error contributes to the flip of order between IPRW and PRW when their true values are close. 

Finally, our experimental results show that the max-sliced Wasserstein estimator works well in practice and converges to some point as the number of samples grow. This supports our consistency results since the max-sliced Wasserstein distance is PRW with $k=1$. Note that there are many existing works on the empirical comparison between max-SW and SW using generative modeling and we refer the interested readers to~\citet{Kolouri-2019-Generalized} and the reference therein.

\section{Conclusion}\label{sec:conclusions}
We study in this paper the statistical aspect of the projection robust Wasserstein (PRW) distance. Our work provides an enhanced understanding of two PRW distances and the associated minimal distance estimators under model misspecification, complementing the existing literature~\citep{Weed-2019-Estimation, Bernton-2019-Parameter, Nadjahi-2019-Asymptotic, Nadjahi-2020-Statistical}. Experiments on synthetic and real datasets highlight some aspects of our theoretical results. Future work includes theory for entropic PRW and the applications of PRW with $k \geq 2$ to deep generative models.   

\section{Acknowledgments}
We would like to thank four anonymous referees for constructive suggestions that improve the quality of this paper. Elynn Y. Chen is supported by National Science Foundation under the grant number DMS-1803241. This work was supported in part by the Mathematical Data Science program of the Office of Naval Research under grant number N00014-18-1-2764.

\bibliographystyle{plainnat}
\bibliography{ref}

\newpage\onecolumn
\appendix
\section{Further Results on the MPRW and MEPRW Estimators}\label{appendix:A}
In this section, we discuss the measurability of the MPRW and MEPRW estimators. 
For a generic function $f$ on the domain $\XCal$, we define $\delta$-$\argmin_{x \in \XCal} f = \{x \in \XCal: f(x) \leq \inf_{x \in \XCal} f + \delta\}$. 
Our results are summarized in the following two theorems. 
\begin{theorem}\label{Theorem:measurability-MPRW}
Under Assumption~\ref{A1}, for any $n \geq 1$ and $\delta > 0$, there exists a Borel measurable function $\widehat{\theta}_n: \Omega \rightarrow \Theta$ such that 
\begin{equation*}
\widehat{\theta}_n(\omega) \in \left\{\begin{array}{ll}
\argmin_{\theta \in \Theta} \overline{\PWCal}_{p, k}(\widehat{\mu}_n(\omega), \mu_\theta) & \textnormal{if this set is nonempty}, \\
\delta\textnormal{-}\argmin_{\theta \in \Theta} \overline{\PWCal}_{p, k}(\widehat{\mu}_n(\omega), \mu_\theta) & \textnormal{otherwise}. 
\end{array}\right.
\end{equation*}
\end{theorem}
\begin{theorem}\label{Theorem:measurability-MEPRW}
Under Assumption~\ref{A1}, for any $n \geq 1$, $m \geq 1$ and $\delta > 0$, there exists a Borel measurable function $\widehat{\theta}_{n, m}: \Omega \rightarrow \Theta$ such that 
\begin{equation*}
\widehat{\theta}_{n, m}(\omega) \in \left\{\begin{array}{ll}
\argmin_{\theta \in \Theta} \EE[\overline{\PWCal}_{p, k}(\widehat{\mu}_n(\omega), \widehat{\mu}_{\theta, m}) \mid X_{1:n}] & \textnormal{if this set is nonempty}, \\
\delta\textnormal{-}\argmin_{\theta \in \Theta} \EE[\overline{\PWCal}_{p, k}(\widehat{\mu}_n(\omega), \widehat{\mu}_{\theta, m}) \mid X_{1:n}] & \textnormal{otherwise}. 
\end{array}\right.
\end{equation*}
\end{theorem}
We also present the asymptotic distribution of the goodness-of-fit statistics as well as the MPRW estimator in the well-specified setting and establish the rate of convergence. For this we require the well separability of the model in Assumption~\ref{A6} and the non-singularity of $D_\star$ in Assumption~\ref{A10} to take place of the local strong identifiability in Assumption~\ref{A9}. 
\begin{assumption}\label{A6}
For any $\epsilon > 0$, there exists $\delta > 0$ so that $\inf_{\theta \in \Theta: \|\theta - \theta_\star\|_\Theta \geq \epsilon} \PWCal_{1, 1}(\mu_{\theta_\star}, \mu_\theta) > \delta$. 
\end{assumption}
\begin{assumption}\label{A10}
There exists a non-singular $D_\star$ such that Assumption~\ref{A7} holds true.   
\end{assumption}
\begin{theorem}\label{Theorem:well-specificied}
Suppose that $\mu_\star = \mu_{\theta_\star}$ for some $\theta_\star$ in the interior of $\Theta$. Under Assumption~\ref{A1}-\ref{A3},~\ref{A7}-\ref{A8} and~\ref{A6}-\ref{A10}, the goodness-of-fit statistics satisfies
\begin{equation*}
\sqrt{n}\inf_{\theta \in \Theta} \overline{\PWCal}_{1, 1}(\widehat{\mu}_n, \mu_\theta) \Rightarrow \inf_{\theta \in \Theta}\max_{u \in \bs^{d-1}} \int_\br |G_\star(u, t) - \langle\theta, D_\star(u, t)\rangle| \; dt, \quad \textnormal{as } n \rightarrow +\infty.   
\end{equation*}
Suppose also that the random map $\theta \rightarrow \max_{u \in \bs^{d-1}} \int_\br |G_\star(u, t) - \langle\theta, D_\star(u, t)\rangle| \; dt$ has a unique infimum almost surely. Then the MPRW estimator of order 1 satisfies 
\begin{equation*}
\sqrt{n}(\widehat{\theta}_n - \theta_\star) \Rightarrow \argmin_{\theta \in \Theta}\max_{u \in \bs^{d-1}} \int_\br |G_\star(u, t) - \langle\theta, D_\star(u, t)\rangle| \; dt, \quad \textnormal{as } n \rightarrow +\infty.   
\end{equation*}
Both the weak convergence results are valid for the metric induced by the norm $\|\cdot\|_L$. 
\end{theorem}

\section{Postponed Proofs in Subsection~\ref{subsec:topology}}
This section lays out the detailed proofs for Lemma~\ref{lemma:PRW-equivalence}, Theorem~\ref{Theorem:PRW-convergence} and \ref{Theorem:PRW-lsc}. 

\subsection{Preliminary technical results}
For completeness, we collect several preliminary technical results\footnote{For the Prokhorov's theorem, we only present the results on the Euclidean space. For more results on general separable metric space, we refer the interested readers to~\citet{Billingsley-2013-Convergence}.} which will be used in the proofs.
\begin{theorem}[Prokhorov's theorem]\label{appendix:Thm-B1}
Let $\PScr(\br^d)$ denote the collection of all probability measures defined on $\br^d$ with the Borel $\sigma$-algebra and $\{\mu_i\}_{i \in \bn}$ is a tight sequence in $\PScr(\br^d)$. Then every subsequence of $\{\mu_i\}_{i \in \bn}$ has a subsequence that converges weakly in $\PScr(\br^d)$. Moreover, if every weakly convergent subsequence has the same limit, the whole sequence converges weakly to this limit.
\end{theorem}
\begin{theorem}[Theorem~4.1 in \citet{Villani-2008-Optimal}]\label{appendix:Thm-B2}
Let $(\XCal, \mu)$ and $(\YCal, \nu)$ be two Polish probability spaces; let $a: \XCal \rightarrow \br \cup \{-\infty\}$ and $b: \YCal \rightarrow \br \cup \{-\infty\}$ be upper semi-continuous such that $a$ and $b$ are absolutely integrable with respect to the measures $\mu$ and $\nu$ respectively. Let $c: \XCal \times \YCal \rightarrow \br \cup \{+\infty\}$ be lower semi-continuous, such that $c(x, y) \geq a(x) + b(y)$ for all $x, y$. Then there exists an optimal coupling $\pi \in \Pi(\mu, \nu)$ which minimizes the total cost $\EE[c(X, Y)]$.
\end{theorem}
\begin{lemma}[Lemma~4.4 in \citet{Villani-2008-Optimal}]\label{appendix:Lem-B3}
Let $\XCal$ and $\YCal$ be two Polish spaces. Let $P \subseteq \PScr(\XCal)$ and $Q \subseteq \PScr(\YCal)$ be tight subsets of $\PScr(\XCal)$ and $\PScr(\YCal)$ respectively. Then the set of all transportation plans whose marginals lie in $P$ and $Q$ respectively, is itself tight in $\PScr(\XCal \times \YCal)$.
\end{lemma}
\begin{theorem}[Theorem~6.9 in \citet{Villani-2008-Optimal}]\label{appendix:Thm-B4}
Let $(\XCal, d)$ be a Polish space and $p \in [1, +\infty)$. The Wasserstein distance $\WCal_p$ metrizes the weak convergence in $\PScr_p(\XCal)$. That is, if $\{\mu_i\}_{i \in \bn_n}$ is a sequence of measures in $\PScr_p(\XCal)$ and $\mu \in \PScr_p(\XCal)$, then $\mu_i \Rightarrow \mu$ if and only if $\WCal_p(\mu_i, \mu) \rightarrow 0$. 
\end{theorem}
\begin{definition}[Lower semi-continuity]
We say that $f: \XCal \rightarrow \br$ is lower semi-continuous if for any $x_0 \in \XCal$ and any $y < f(x_0)$, there exists a neighborhood $U$ of $x_0$ such that $f(x) > y$ for all $x$ in $U$. In the case of a metric space, this is equivalent to $\liminf_{x \rightarrow x_0} f(x) \geq f(x_0)$ for any $x_0 \in \XCal$.  
\end{definition}

\subsection{Proof of Lemma~\ref{lemma:PRW-equivalence}}
We first show that, for any $\mu \in \PScr_p(\br^d)$ and $\nu \in \PScr_p(\br^d)$, the following inequality holds true, 
\begin{equation}\label{inequality-weak-first}
\underline{\PWCal}_{p, k}(\mu, \nu) \leq \overline{\PWCal}_{p, k}(\mu, \nu) \leq \WCal_p(\mu, \nu). 
\end{equation}
Indeed, by the definition of $\underline{\PWCal}_{p, k}$ and $\overline{\PWCal}_{p, k}$, the first inequality is trivial. For the second inequality, we derive from the definition of $\overline{\PWCal}_{p, k}$ that 
\begin{eqnarray*}
\overline{\PWCal}_{p, k}^p(\mu, \nu) = \sup_{E \in \bs_{d, k}} \WCal_p^p(E_{\#}^\star\mu, E_{\#}^\star\nu) = \sup_{E \in \bs_{d, k}} \inf_{\pi \in \Pi(\mu, \nu)} \int_{\br^d \times \br^d} \|E^\top(x-y)\|^p \; d\pi(x, y). 
\end{eqnarray*}
Since $E \in \bs_{d, k}$, we have $\|E^\top(x-y)\| \leq \|x-y\|$. Thus, we have $\overline{\PWCal}_{p, k}^p(\mu, \nu) \leq \WCal_p^p(\mu, \nu)$. Putting these pieces together yields Eq.~\eqref{inequality-weak-first}. For any sequence $\{\mu_i\}_{i \in \bn} \subseteq \PScr_p(\br^d)$ and $\mu \in \PScr_p(\br^d)$, we conclude from Eq.~\eqref{inequality-weak-first} that $\WCal_p(\mu_i, \mu) \rightarrow 0$ implies $\overline{\PWCal}_{p, k}(\mu_i, \mu) \rightarrow 0$ and $\underline{\PWCal}_{p, k}(\mu_i, \mu) \rightarrow 0$. 

The remaining step is to show that $\underline{\PWCal}_{p, k}(\mu_i, \mu) \rightarrow 0$ implies $\WCal_p(\mu_i, \mu) \rightarrow 0$. Indeed, we first prove that $\underline{\PWCal}_{p, k}(\mu_i, \mu) \rightarrow 0$ implies $\mu_i \Rightarrow \mu$. Let $Z_i \sim \mu_i$, we have $E^\top Z_i \sim E_{\#}^\star\mu_i$. By the definition of the IPRW distance (cf. Definition~\ref{def:IPRW}) and using the fact that $\underline{\PWCal}_{p, k}(\mu_i, \mu) \rightarrow 0$, we have $(\|E^\top Z_i)\|^p)_{i \in \bn}$ is uniformly integrable for all $E \in \bs_{d, k}$. Since $\bs_{d, k}$ is compact, there exists a finite set $\{E_1, E_2, \ldots, E_I\} \subseteq \bs_{d, k}$ so that $\|x\| \leq \sum_{j=1}^I \|E_j^\top x\|$ for all $x \in \br^d$. 
Therefore, we have
\begin{equation*}
\|Z_i\|^p \leq \left(\sum_{j=1}^I \|E_j^\top Z_i\|\right)^p \leq I^p\left(\max_{1 \leq j \leq I} \|E_j^\top Z_i\|^p\right) \leq I^p\left(\sum_{j=1}^I \|E_j^\top Z_i\|^p\right).  
\end{equation*}
Therefore, we deduce that $(\|Z_i\|^p)_{i \in \bn}$ is uniformly integrable which implies the tightness of $\{\mu_i\}_{i \in \bn}$. Using the Prokhorov's theorem (cf. Theorem~\ref{appendix:Thm-B1}), we obtain that every subsequence of $\{\mu_i\}_{i \in \bn}$ has a weakly convergent subsequence. 

The next step is to show that all the weakly convergent subsequences converge to the same probability measure $\mu$. We fix an arbitrary subsequence and for simplicity abbreviate the subscripts and still denote it by $\{\mu_i\}_{i \in \bn}$. Let $\tilde{\mu}_i$ be the limit of any given weakly convergent subsequence $(\mu_{i_j})_{j \in \bn}$, we need to prove that $\tilde{\mu}_i = \mu$. In particular, we define the characteristic function for any probability measure $\nu$ as follows, 
\begin{equation*}
\Phi_\nu(z) \mydefn \int_{\br^d} e^{\textsf{i}\langle z, x\rangle} \; d\nu(x) \quad \text{for all } z \in \br^d. 
\end{equation*}
Since $\mu_{i_j} \Rightarrow \tilde{\mu}_i$, we have $\Phi_{\mu_{i_j}}(z) \rightarrow \Phi_{\tilde{\mu}_i}(z)$ for all $z \in \br^d$. Thus, we need to show that $\Phi_{\mu_{i_j}}(z) \rightarrow \Phi_{\mu}(z)$ for all $z \in \br^d$. This is trivial when $z=\zero_d$ since $\Phi_{\mu_{i_j}}(\zero_d) = \Phi_{\mu}(\zero_d) = 1$ for all $j \in \bn$. Otherwise, let $r \mydefn \|z\|$ and $v \mydefn z/\|z\|$, we have
\begin{equation*}
\lim_{j \rightarrow +\infty} \Phi_{\mu_{i_j}}(z) = \lim_{j \rightarrow +\infty} \int_{\br^d} e^{\textsf{i}\langle z, x\rangle} \; d\mu_{i_j}(x) = \lim_{j \rightarrow +\infty} \int_{\br^d} e^{\textsf{i} r\langle v, x\rangle} \; d\mu_{i_j}(x). 
\end{equation*}
Since $\|v\|=1$, we define $\bar{E} \in \bs_{d, k}$ whose first column is $v$. Let $\bar{r}$ be a $k$-dimensional vector whose first coordinate is $r$ and other coordinates are zero. Then we have $r\langle v, x\rangle = \langle \bar{r}, \bar{E}^\top x\rangle$. Putting these pieces together yields that 
\begin{equation*}
\lim_{j \rightarrow +\infty} \Phi_{\mu_{i_j}}(z) = \lim_{j \rightarrow +\infty} \int_{\br^k} e^{\textsf{i}\langle\bar{r}, y\rangle} \; d \bar{E}_{\#}^\star\mu_{i_j}(y). 
\end{equation*}
For such fixed $\bar{E}$, we claim that $\WCal_p(\bar{E}_{\#}^\star\mu_{i_j}, \bar{E}_{\#}^\star\mu) \rightarrow 0$ holds true. More specifically, $\underline{\PWCal}_{p, k}(\mu_{i_j}, \mu) \rightarrow 0$ implies that $\int \mathcal{W}_p^p(E_\#^\star\mu_{i_j}, E_\#^\star\mu) d\sigma(E) \rightarrow 0$. Since $\mathcal{W}_p^p(E_\#^\star\mu_{i_j}, E_\#^\star\mu)$ is non-negative, it is easy to derive that $\mathcal{W}_p(\bar{E}_\#^\star\mu_{i_j}, \bar{E}_\#^\star\mu) \nrightarrow 0$ for almost every $E$. Nonetheless, by the continuity of $\mathcal{W}_p^p(E_\#^\star\mu_{i_j}, E_\#^\star\mu)$ with respect to $E$, we can obtain that $\mathcal{W}_p(\bar{E}_\#^\star\mu_{i_j}, \bar{E}_\#^\star\mu) \nrightarrow 0$ for all fixed $E$. Indeed, by the proof by contradiction, we assume that $\mathcal{W}_p(E_\#^\star\mu_{i_j}, E_\#^\star\mu) \nrightarrow 0$ for some fixed $E$. Then, there exists a neighborhood $S$ of $E$ (it is fixed) such that $\int_S \mathcal{W}_p^p(E_\#^\star\mu_{i_j}, E_\#^\star\mu) d\sigma(E) \nrightarrow 0$. This contradicts $\int \mathcal{W}_p^p(E_\#^\star\mu_{i_j}, E_\#^\star\mu) d\sigma(E) \rightarrow 0$ since the inside term is non-negative. Thus, we achieve the desired claim. 

Using Theorem~\ref{appendix:Thm-B4}, we have $\bar{E}_{\#}^\star\mu_{i_j} \Rightarrow \bar{E}_{\#}^\star\mu$. Since $r\langle v, x\rangle = \langle\bar{r}, \bar{E}^\top x\rangle$, we have
\begin{equation*}
\lim_{j \rightarrow +\infty} \int_{\br^k} e^{\textsf{i}\langle\bar{r}, x\rangle} \; d\bar{E}_{\#}^\star\mu_{i_j}(x) = \int_{\br^k} e^{\textsf{i}\langle\bar{r}, x\rangle} \; d\bar{E}_{\#}^\star\mu(x) = \int_{\br^d} e^{\textsf{i}r\langle v, x\rangle} \; d\mu(x) = \int_{\br^d} e^{\textsf{i}\langle z, x\rangle} \; d\mu(x). 
\end{equation*}
Putting these pieces together yields that $\Phi_{\mu_{i_j}}(z) \rightarrow \Phi_{\mu}(z)$ for all $z \in \br^d/\{\zero_d\}$ and $\tilde{\mu}_i = \mu$ for all $i \in \bn$. Using the Prokhorov's theorem again yields that the whole sequence $\{\mu_i\}_{i \in \bn}$ has the limit $\mu$ in weak sense. Therefore, $\underline{\PWCal}_{p, k}(\mu_i, \mu) \rightarrow 0$ implies $\mu_i \Rightarrow \mu$. Since the Wasserstein distances metrize the weak convergence (cf. Theorem~\ref{appendix:Thm-B4}), we conclude that $\underline{\PWCal}_{p, k}(\mu_i, \mu) \rightarrow 0$ implies $\WCal_p(\mu_i, \mu) \rightarrow 0$. This completes the proof. 

\subsection{Proof of Theorem~\ref{Theorem:PRW-convergence}}
By Lemma~\ref{lemma:PRW-equivalence}, we have $\underline{\PWCal}_{p, k}(\mu_i, \mu) \rightarrow 0$ if and only if $\overline{\PWCal}_{p, k}(\mu_i, \mu) \rightarrow 0$ if and only if $\WCal_p(\mu_i, \mu) \rightarrow 0$. By Theorem~\ref{appendix:Thm-B4}, we have $\mu_i \Rightarrow \mu$ if and only if $\WCal_p(\mu_i, \mu) \rightarrow 0$. Putting these pieces together yields the desired result. 

\subsection{Proof of Theorem~\ref{Theorem:PRW-lsc}}
Fixing $E \in \bs_{d, k}$, the mapping $x \mapsto E^\top x$ is continuous from $\br^d$ to $\br^k$. Since $\mu_i \Rightarrow \mu$ and $\nu_i \Rightarrow \nu$, the continuous mapping theorem implies that $E_{\#}^\star\mu_i \Rightarrow E_{\#}^\star\mu$ and $E_{\#}^\star\nu_i \Rightarrow E_{\#}^\star\nu$. The next step is the key ingredient in the proof and we hope to show that 
\begin{equation}\label{inequality-weak-fourth}
\WCal_p^p(E_{\#}^\star\mu, E_{\#}^\star\nu) \leq \liminf_{i \rightarrow +\infty} \WCal_p^p(E_{\#}^\star\mu_i, E_{\#}^\star\nu_i) \quad \textnormal{for all } E \in \bs_{d, k}.  
\end{equation}
From Theorem~\ref{appendix:Thm-B2}, there exists a coupling $\pi_i \in \Pi(E_{\#}^\star\mu_i, E_{\#}^\star\nu_i)$ such that $\WCal_p^p(E_{\#}^\star\mu_i, E_{\#}^\star\nu_i) = \int_{\br^k \times \br^k} \|x-y\|^p \; d\pi_i(x, y)$. By the definition of $\liminf$, there exists a subsequence of $\{\pi_i\}_{i \in \bn}$ such that $\int_{\br^k \times \br^k} \|x-y\|^p \; d\pi_i(x, y)$ converges to $\liminf_{i \rightarrow +\infty} \WCal_p^p(E_{\#}^\star\mu_i, E_{\#}^\star\nu_i)$. For the simplicity, we still denote it by $\{\pi_i\}_{i \in \bn}$. By Lemma~\ref{appendix:Lem-B3} and Prokhorov's theorem (cf. Theorem~\ref{appendix:Thm-B1}), $\{\pi_i\}_{i \in \bn}$ is sequentially compact in weak sense. Thus, there exists a subsequence $\{\pi_{i_j}\}_{j \in \bn}$ such that $\pi_{i_j} \Rightarrow \tilde{\pi} \in \PScr(\br^k \times \br^k)$. Putting these pieces together yields that 
\begin{equation*}
\liminf_{i \rightarrow +\infty} \WCal_p^p(E_{\#}^\star\mu_i, E_{\#}^\star\nu_i) = \int_{\br^k \times \br^k} \|x-y\|^p \; d\tilde{\pi}(x, y). 
\end{equation*}
By the definition of the Wasserstein distance, it suffices to show that $\tilde{\pi} \in \Pi(E_{\#}^\star\mu, E_{\#}^\star\nu)$. Indeed, let $f: \br^k \rightarrow \br$ be a continuous and bounded function, we have
\begin{equation*}
\int_{\br^k \times \br^k} f(x) \; d\tilde{\pi}(x, y) = \lim_{j \rightarrow +\infty} \int_{\br^k \times \br^k} f(x) \; d\pi_{i_j}(x, y). 
\end{equation*}  
Since $\pi_{i_j} \in \Pi(E_{\#}^\star\mu_{i_j}, E_{\#}^\star\nu_{i_j})$ and $E_{\#}^\star\mu_i \Rightarrow E_{\#}^\star\mu$, we have
\begin{equation*}
\lim_{j \rightarrow +\infty} \int_{\br^k \times \br^k} f(x) \; d\pi_{i_j}(x, y) = \lim_{j \rightarrow +\infty} \int_{\br^k} f(x) \; dE_{\#}^\star\mu_{i_j}(x) = \int_{\br^k} f(x) \; dE_{\#}^\star\mu(x). 
\end{equation*}
Since $E_{\#}^\star\nu_i \Rightarrow E_{\#}^\star\nu$, the same argument implies that $\int_{\br^k \times \br^k} f(y) \; d\tilde{\pi}(x, y) = \int_{\br^k} f(y) \; dE_{\#}^\star\nu(y)$. Putting these pieces together yields Eq.~\eqref{inequality-weak-fourth}.

For the IPRW distance, we derive from Eq.~\eqref{inequality-weak-fourth} and the Fatou's lemma that 
\begin{equation*}
\underline{\PWCal}_{p, k}^p(\mu, \nu) = \int_{\bs_{d, k}} \WCal_p^p(E_{\#}^\star\mu, E_{\#}^\star\nu) \; d\sigma(E) \leq \liminf_{i \rightarrow +\infty}  \int_{\bs_{d, k}} \WCal_p^p(E_{\#}^\star\mu_i, E_{\#}^\star\nu_i) \; d\sigma(E) = \liminf_{i \rightarrow +\infty} \underline{\PWCal}_{p, k}^p(\mu_i, \nu_i).   
\end{equation*}
Since $\underline{\PWCal}_{p, k}(\mu, \nu)$ and $\underline{\PWCal}_{p, k}(\mu_i, \nu_i)$ are both nonnegative, we take the $p$-th root of both sides of the above inequality and have $\underline{\PWCal}_{p, k}(\mu, \nu) \leq \liminf_{i \rightarrow +\infty} \underline{\PWCal}_{p, k}(\mu_i, \nu_i)$. 

For the PRW distance, we derive from Eq.~\eqref{inequality-weak-fourth} and the fact that the supremum of a sequence of lower semi-continuous mappings is lower semi-continuous that 
\begin{equation*}
\overline{\PWCal}_{p, k}^p(\mu, \nu) = \sup_{E \in \bs_{d, k}} \WCal_p^p(E_{\#}^\star\mu, E_{\#}^\star\nu) \leq \liminf_{i \rightarrow +\infty} \overline{\PWCal}_{p, k}^p(\mu_i, \nu_i). 
\end{equation*}
where the first equality holds true since the Wasserstein distance is nonnegative. Since $\overline{\PWCal}_{p, k}(\mu, \nu)$ and $\overline{\PWCal}_{p, k}(\mu_i, \nu_i)$ are both nonnegative, we have $\overline{\PWCal}_{p, k}(\mu, \nu) \leq \liminf_{i \rightarrow +\infty} \overline{\PWCal}_{p, k}(\mu_i, \nu_i)$. 

\section{Postponed Proofs in Subsection~\ref{subsec:convergence}}
In this section, we provide the detailed proofs for Theorem~\ref{Theorem:IPRW-main}-\ref{Theorem:concentration-poincare}. 

\subsection{Preliminary technical results}
To facilitate reading, we collect several preliminary technical results which will be used in the postponed proofs in subsection~\ref{subsec:convergence}. 
\begin{theorem}[Tonelli's theorem]\label{Theorem:appendix-Tonelli}
if $(\XCal, A, \mu)$ and $(\YCal, B, \nu)$ are $\sigma$-finite measure spaces, while $f: \XCal \times \YCal \rightarrow [0, +\infty]$ is non-negative measurable function, then
\begin{equation*}
\int_\XCal \left( \int_\YCal f(x, y) \; dy\right) \; dx = \int_\YCal \left(\int_\XCal f(x,y) \; dx\right) \; dy = \int_{\XCal \times \YCal} f(x,y) \; d(x,y).
\end{equation*}
\end{theorem}
The following proposition provides the state-of-the-art general bound for the Wasserstein distance between the true measure and its empirical version in $\br^d$. Note that we do not assume any additional structures of the true measure. Similar results can be found in many classical works, e.g.,~\citet[Theorem~1]{Fournier-2015-Rate},~\citet[Theorem~1]{Weed-2019-Sharp} and~\citet[Theorem~3.1]{Lei-2020-Convergence}. Since $p \geq 1$, we present the following results which directly follows the proof of~\citet[Theorem~3.1]{Lei-2020-Convergence}.  
\begin{proposition}\label{Proposition:appendix-Wasserstein}
Let $\mu_\star \in \PScr_q(\br^d)$ and $M_q \mydefn M_q(\mu_\star) < +\infty$. Then we have
\begin{equation}
\EE[\WCal_p(\widehat{\mu}_n, \mu_\star)] \leq (\EE[\WCal_p^p(\widehat{\mu}_n, \mu_\star)])^{1/p} \lesssim_{p, q} n^{-[\frac{1}{(2p)\vee d}\wedge(\frac{1}{p}-\frac{1}{q})]}(\log(n))^{\frac{\zeta'_{p,q,d}}{p}}, \quad \textnormal{for all } n \geq 1.  
\end{equation}
where $\lesssim_{p, q}$ refers to ``less than" with a constant depending only on $(p, q)$ and 
\begin{equation*}
\zeta'_{p, q, d} = \left\{\begin{array}{ll}
2 & \textnormal{if } d = q = 2p, \\
1 & \textnormal{if } ``d \neq 2p \text{ and } q = \frac{dp}{d-p}" \textnormal{ or } ``q > d = 2p", \\
0 & \textnormal{otherwise.}
\end{array}\right. 
\end{equation*}
\end{proposition}
The following proposition provides a bound for the covering number of $\bs_{d, k}$ in the operator norm of a matrix, denoted by $\|\cdot\|_\op$. This is a straightforward consequence of the classical results on the covering number of the unit sphere in $\br^d$ in Euclidean norm. For the proof details, we refer the interested readers to~\citet[Lemma~4]{Weed-2019-Estimation}. For the background materials on the covering number, we refer the interested readers to~\citet[Chapter~5]{Wainwright-2019-High}. For the ease of presentation, we provide a formal definition of covering number of $\bs_{d, k}$ in $\|\cdot\|_\op$ as follows. 

For any $\epsilon \in (0, 1)$, the $\epsilon$-covering number of $\bs_{d, k}$ in $\|\cdot\|_\op$ is defined by 
\begin{equation*}
N(\bs_{d, k}, \epsilon, \|\cdot\|_\op) = \inf\left\{N \in \bn: \exists x_1, x_2, \ldots, x_N \in \bs_{d, k}, \ \st \ \bs_{d, k} \subseteq \bigcup_{i=1}^N \BB(x_i, \epsilon)\right\}, 
\end{equation*}
where $\BB(x, r) = \{y \in \bs_{d, k}: \|y - x\|_\op \leq r\}$ is the ball of radius $r > 0$ centered at $x \in \bs_{d, k}$ in the operator norm of a matrix. 
\begin{proposition}\label{Proposition:appendix-covering}
There exists a universal constant $c > 0$ such that for all $\epsilon \in (0, 1)$, the $\epsilon$-covering number of $\bs_{d, k}$ in $\|\cdot\|_\op$ satisfies that $N(\bs_{d, k}, \epsilon, \|\cdot\|_\op) \leq (c\sqrt{k}\epsilon^{-1})^{dk}$. 
\end{proposition}
The following theorem~\citep{Lei-2020-Convergence} summarizes the concentration results assuming the Bernstein tail condition under product measure. Indeed, let $\{X_i\}_{i \in [n]}$ be independent samples from probability measure $\mu_i$ on spaces $\XCal_i$ and $X'_i$ be independent copies of $X_i$ for all $i \in [n]$. Denote $X = (X_1, \ldots, X_n)$ and $X'_{(i)} = (X_1, \ldots, X'_i, \ldots, X_n)$ which is identical to $X$ except for $X'_i$. Let $f: \prod_{i=1}^n \XCal_i \rightarrow \br$ be a function such that $\EE[|f(X)|] < +\infty$, and define $D_i = f(X) - f(X'_{(i)})$.
\begin{theorem}\label{Theorem:appendix-BM}
Suppose that there exists some $\sigma_i, M > 0$ so that $\EE[|D_i|^k \mid X_{-i}] \leq (1/2)\sigma_i^2 k! M^{k-2}$ for all $k \geq 2$. Then the following statement holds, 
\begin{equation*}
\PP(f(X) - \EE(f(X)) > t) \leq \exp\left(-\frac{t^2}{2(\sum_{i=1}^n \sigma_i^2) + 2tM}\right). 
\end{equation*}
\end{theorem}
The following theorem summarizes the concentration results assuming the Poincar\'{e} inequality under product measure. We denote by $\|\nabla_i f\|$ the length of the gradient with respect to the $i^\textnormal{th}$ coordinate. 
\begin{theorem}[Corollary~4.6 in~\citet{Ledoux-1999-Concentration}]\label{Theorem:appendix-Poincare}
Denote by $\mu^n$ the product of $\mu$ on $\otimes_{i=1}^n \br^d$ and $\mu \in \PScr(\br^d)$ satisfies the Poincar\'{e} inequality (cf. Definition~\ref{def:Poincare}). For every function $f$ on $\otimes_{i=1}^n \br^d$ satisfying $\EE(|f(X)|) < +\infty$, and $\sum_{i=1}^n \|\nabla_i f(X)\|^2 \leq \alpha^2$ and $\max_{1 \leq i \leq n} \|\nabla_i f(X)\| \leq \beta$ almost surely. Then the following statement holds true for $X \sim \mu^n$ that, 
\begin{equation*}
\PP(f(X) - \EE(f(X)) > t) \leq \exp\left(-\frac{1}{K}\min\left\{\frac{t}{\beta}, \ \frac{t^2}{\alpha^2}\right\}\right), 
\end{equation*}
where $K > 0$ only depends on the constant $M$ in the Poincar\'{e} inequality. 
\end{theorem}

\subsection{Proof of Theorem~\ref{Theorem:IPRW-main}}
Note that $\mu_\star \in \PScr_q(\br^d)$ and $M_q \mydefn M_q(\mu_\star) < +\infty$. Fixing $E \in \bs_{d, k}$, we have $E_{\#}^\star\mu_\star \in \PScr_q(\br^k)$ and $M_q(E_{\#}^\star\mu_\star) \leq M_q < +\infty$. Then Proposition~\ref{Proposition:appendix-Wasserstein} implies that 
\begin{equation*}
\left(\EE[\WCal_p^p(E_{\#}^\star\widehat{\mu}_n, E_{\#}^\star\mu_\star)]\right)^{1/p} \lesssim_{p, q} n^{-[\frac{1}{(2p)\vee k}\wedge(\frac{1}{p}-\frac{1}{q})]}(\log(n))^{\frac{\zeta'_{p,q,k}}{p}} \quad \textnormal{for all } n \geq 1.  
\end{equation*}
Since $\WCal_p(E_{\#}^\star\widehat{\mu}_n, E_{\#}^\star\mu_\star) \geq 0$ for any $E \in \bs_{d, k}$ and $\mu_\star \in \PScr_q(\br^d)$, Theorem~\ref{Theorem:appendix-Tonelli} implies that 
\begin{equation*}
\EE[\underline{\PWCal}_{p, k}^p(\widehat{\mu}_n, \mu_\star)] = \EE\left[\int_{\bs_{d, k}} \WCal_p^p(E_{\#}^\star\widehat{\mu}_n, E_{\#}^\star\mu_\star) \; d\sigma(E)\right] = \int_{\bs_{d, k}} \EE[\WCal_p^p(E_{\#}^\star\widehat{\mu}_n, E_{\#}^\star\mu_\star)] \; d\sigma(E). 
\end{equation*}
Note that $\zeta_{p,q,k} = \zeta'_{p,q,k}$ where $\zeta_{p,q,k}$ is defined in Theorem~\ref{Theorem:IPRW-main}. Moreover, $p \geq 1$. By the Jensen's inequality, we have 
\begin{equation}
\EE[\underline{\PWCal}_{p, k}(\widehat{\mu}_n, \mu_\star)] \leq (\EE[\underline{\PWCal}_{p, k}^p(\widehat{\mu}_n, \mu_\star)])^{1/p}. 
\end{equation}
Putting these pieces together yields the desired result. 

\subsection{Proof of Theorem~\ref{Theorem:PRW-main-Bernstein}}
By the definition of $\overline{\PWCal}_{p, k}(\widehat{\mu}_n, \mu_\star)$, we have
\begin{equation}\label{inequality-PRW-main-first}
\EE[\overline{\PWCal}_{p, k}(\widehat{\mu}_n, \mu_\star)] \leq \sup_{E \in \bs_{d, k}} \EE[\WCal_p(E_{\#}^\star\widehat{\mu}_n, E_{\#}^\star\mu_\star)] + \EE\left[\sup_{E \in \bs_{d, k}} \left(\WCal_p(E_{\#}^\star\widehat{\mu}_n, E_{\#}^\star\mu_\star) - \EE[\WCal_p(E_{\#}^\star\widehat{\mu}_n, E_{\#}^\star\mu_\star)]\right)\right].
\end{equation}
Using the same arguments for proving Theorem~\ref{Theorem:IPRW-main}, we have
\begin{equation}\label{inequality-PRW-main-second}
\sup_{E \in \bs_{d, k}} \EE[\WCal_p(E_{\#}^\star\widehat{\mu}_n, E_{\#}^\star\mu_\star)] \lesssim_{p, q} n^{-[\frac{1}{(2p)\vee k}\wedge(\frac{1}{p}-\frac{1}{q})]}(\log(n))^{\frac{\zeta_{p,q,k}}{p}} \quad \textnormal{for all } n \geq 1.
\end{equation}
The remaining step is to bound the gap $\EE[\sup_{E \in \bs_{d, k}} (\WCal_p(E_{\#}^\star\widehat{\mu}_n, E_{\#}^\star\mu_\star) - \EE[\WCal_p(E_{\#}^\star\widehat{\mu}_n, E_{\#}^\star\mu_\star)])]$. We first claim that $\WCal_p(E_{\#}^\star\widehat{\mu}_n, E_{\#}^\star\mu_\star) - \EE[\WCal_p(E_{\#}^\star\widehat{\mu}_n, E_{\#}^\star\mu_\star)]$ is sub-exponential with parameters $(2\sigma n^{1/2-1/p}, 2Vn^{-1/p})$ for all $E \in \bs_{d, k}$ if the true measure $\mu_\star$ satisfies the projection Bernstein-type tail condition (cf. Definition~\ref{def:proj-Bernstein}). Indeed, let $f(X)=\WCal_p(E_{\#}^\star\widehat{\mu}_n, E_{\#}^\star\mu_\star)$, we have
\begin{equation*}
D_i = f(X) - f(X'_{(i)}) \leq \WCal_p(E_{\#}^\star\widehat{\mu}_n, E_{\#}^\star\widehat{\mu}'_n) \leq n^{-1/p}\left(\|E_{\#}^\star(X_i) - E_{\#}^\star(X'_i)\|\right). 
\end{equation*}
By the triangle inequality and using the projection Bernstein-type tail condition, we have
\begin{equation*}
\EE[|D_i|^k \mid X_{-i}] \leq 2^k n^{-k/p}(\EE_{X \sim E_{\#}^\star\mu}[|X|^k]) \leq 2^{k-1} n^{-k/p} \sigma^2 k! V^{k-2} = \frac{(2n^{-1/p}\sigma)^2k!(2n^{-1/p}V)^{k-2}}{2}. 
\end{equation*}
This implies that the condition in Theorem~\ref{Theorem:appendix-BM} holds true with $\sigma_i = 2n^{-1/p}\sigma$ and $M = 2n^{-1/p}V$. Equipped with Theorem~\ref{Theorem:appendix-BM} yields that 
\begin{equation*}
\PP\left(\WCal_p(E_{\#}^\star\widehat{\mu}_n, E_{\#}^\star\mu_\star) - \EE[\WCal_p(E_{\#}^\star\widehat{\mu}_n, E_{\#}^\star\mu_\star)] \geq t\right) \leq \exp\left(-\frac{t^2}{8\sigma^2 n^{1-2/p} + 4tVn^{-1/p}}\right). 
\end{equation*}
For the simplicity, let $Z_E = \WCal_p(E_{\#}^\star\widehat{\mu}_n, E_{\#}^\star\mu_\star) - \EE[\WCal_p(E_{\#}^\star\widehat{\mu}_n, E_{\#}^\star\mu_\star)]$. Then we have $\EE[Z_E] = 0$ and $\PP(Z_E \geq t) \leq \exp(-t^2/(8\sigma^2 n^{1-2/p} + 4tVn^{-1/p}))$. This together with the definition of $Z_E$ and~\citet[Theorem~2.2]{Wainwright-2019-High} yields the desired claim. 

We then interpret $\{Z_E\}_{E \in \bs_{d, k}}$ as an empirical process indexed by $E \in \bs_{d, k}$ and claim that there exists a random variable $L$ satisfying $\EE[L] \leq 4M_q(\mu_\star)$ so that $|Z_U - Z_V| \leq L\|U - V\|_\op$ for all $U, V \in \bs_{d, k}$. More specifically, it follows from the definition that 
\begin{equation*}
Z_U - Z_V \ = \ \left(\mathcal{W}_p(U_\#^\star\widehat{\mu}_n, U_\#^\star\mu) - \mathcal{W}_p(V_\#^\star\widehat{\mu}_n, V_\#^\star\mu)\right) - \EE\left[\mathcal{W}_p(U_\#^\star\widehat{\mu}_n, U_\#^\star\mu)-\mathcal{W}_p(V_\#^\star\widehat{\mu}_n, V_\#^\star\mu)\right]. 
\end{equation*}
Since the Wasserstein distance is nonnegative and satisfies the triangle inequality, we have
\begin{eqnarray*}
\mathcal{W}_p(U_\#^\star\widehat{\mu}_n, U_\#^\star\mu) - \mathcal{W}_p(V_\#^\star\widehat{\mu}_n, V_\#^\star\mu) & = & \mathcal{W}_p(U_\#^\star\widehat{\mu}_n, U_\#^\star\mu)-\mathcal{W}_p(U_\#^\star\widehat{\mu}_n, V_\#^\star\mu)+\mathcal{W}_p(U_\#^\star\widehat{\mu}_n, V_\#^\star\mu)-\mathcal{W}_p(V_\#^\star\widehat{\mu}_n, V_\#^\star\mu) \\ 
& & \hspace*{-5em} \leq \ \mathcal{W}_p(U_\#^\star\mu, V_\#^\star\mu)+\mathcal{W}_p(U_\#^\star\widehat{\mu}_n, V_\#^\star\widehat{\mu}_n)
\end{eqnarray*}
Putting these pieces together yields that 
\begin{equation*}
Z_U - Z_V \ \leq \ \WCal_p(U_{\#}^\star\widehat{\mu}_n, V_{\#}^\star\widehat{\mu}_n) + \WCal_p(U_{\#}^\star\mu_\star, V_{\#}^\star\mu_\star) + \EE\left[\WCal_p(U_{\#}^\star\widehat{\mu}_n, V_{\#}^\star\widehat{\mu}_n) + \WCal_p(U_{\#}^\star\mu_\star, V_{\#}^\star\mu_\star)\right]. 
\end{equation*}
Since the Wasserstein distance is symmetrical, we have
\begin{equation*}
Z_V - Z_U \ \leq \ \WCal_p(U_{\#}^\star\widehat{\mu}_n, V_{\#}^\star\widehat{\mu}_n) + \WCal_p(U_{\#}^\star\mu_\star, V_{\#}^\star\mu_\star) + \EE\left[\WCal_p(U_{\#}^\star\widehat{\mu}_n, V_{\#}^\star\widehat{\mu}_n) + \WCal_p(U_{\#}^\star\mu_\star, V_{\#}^\star\mu_\star)\right]. 
\end{equation*}
Therefore, we conclude that 
\begin{equation*}
|Z_U - Z_V| \ \leq \ \WCal_p(U_{\#}^\star\widehat{\mu}_n, V_{\#}^\star\widehat{\mu}_n) + \WCal_p(U_{\#}^\star\mu_\star, V_{\#}^\star\mu_\star) + \EE\left[\WCal_p(U_{\#}^\star\widehat{\mu}_n, V_{\#}^\star\widehat{\mu}_n) + \WCal_p(U_{\#}^\star\mu_\star, V_{\#}^\star\mu_\star)\right]. 
\end{equation*}
Let $X \sim \mu$, we have
\begin{eqnarray*}
|Z_U - Z_V| & \leq & 2\left(\EE(\|(U-V)X\|^p)\right)^{1/p} + \left(\frac{1}{n}\sum_{i=1}^n \|(U-V)X_i\|^p \right)^{1/p} + \EE\left[\left(\frac{1}{n}\sum_{i=1}^n \|(U-V)X_i\|^p \right)^{1/p}\right] \\
& \leq & \|U-V\|_\op\left(2(\EE(\|X\|^p))^{1/p} + \left(\frac{1}{n}\sum_{i=1}^n \|X_i\|^p \right)^{1/p} + \EE\left[\left(\frac{1}{n}\sum_{i=1}^n \|X_i\|^p \right)^{1/p}\right]\right) \\ 
& \mydefn & L\|U-V\|_\op. 
\end{eqnarray*}
Note that $X_{1:n} = (X_1, \ldots, X_n)$ are independent and identically distributed samples according to $\mu_\star$. By the Jensen's inequality and using the fact that $q > p \geq 1$, we have
\begin{equation*}
\EE[L] \leq 4(\EE(\|X\|^p))^{1/p} \leq 4(\EE(\|X\|^q))^{1/q} = 4M_q(\mu_\star). 
\end{equation*}
Thus, by a standard $\epsilon$-net argument, we obtain that 
\begin{equation*}
\EE[\sup_{E \in \bs_{d, k}} Z_E] \leq \inf_{\epsilon > 0}\left\{\epsilon\EE[L] + 4\sigma n^{1/2-1/p}\sqrt{\log(N(\bs_{d, k}, \epsilon, \|\cdot\|_\op))} + 2Vn^{-1/p}\log(N(\bs_{d, k}, \epsilon, \|\cdot\|_\op))\right\}
\end{equation*}
Proposition~\ref{Proposition:appendix-covering} shows that there exists a universal constant $c > 0$ such that 
\begin{equation*}
\log(N(\bs_{d, k}, \epsilon, \|\cdot\|_\op)) \leq dk\log\left(\frac{c\sqrt{k}}{\epsilon}\right). 
\end{equation*}
Putting these pieces together and choosing $\epsilon = \sqrt{k}n^{-1/p}$ (it is chosen to achieve the tight bound) yields that 
\begin{eqnarray*}
\EE\left[\sup_{E \in \bs_{d, k}} Z_E \right] & \lesssim_{p, q} & \inf_{\epsilon > 0}\left\{\epsilon + n^{1/2-1/p}\sqrt{dk\log\left(\frac{\sqrt{k}}{\epsilon}\right)} + n^{-1/p}dk\log\left(\frac{\sqrt{k}}{\epsilon}\right)\right\} \\
& \lesssim_{p, q} & n^{1/2-1/p}\sqrt{dk\log(n)} + n^{-1/p}dk\log(n). 
\end{eqnarray*}
Therefore, we conclude that 
\begin{equation*}
\EE\left[\sup_{E \in \bs_{d, k}} \left(\WCal_p(E_{\#}^\star\widehat{\mu}_n, E_{\#}^\star\mu_\star) - \EE[\WCal_p(E_{\#}^\star\widehat{\mu}_n, E_{\#}^\star\mu_\star)]\right)\right] \lesssim_{p, q} n^{1/2-1/p}\sqrt{dk\log(n)} + n^{-1/p}dk\log(n).   
\end{equation*}
This together with Eq.~\eqref{inequality-PRW-main-first} and Eq.~\eqref{inequality-PRW-main-second} yields the desired inequality. 

\subsection{Proof of Theorem~\ref{Theorem:PRW-main-Poincare}}
Using the same arguments in Theorem~\ref{Theorem:PRW-main-Bernstein}, we obtain Eq.~\eqref{inequality-PRW-main-first} and Eq.~\eqref{inequality-PRW-main-second}. So it suffices to bound the gap $\EE[\sup_{E \in \bs_{d, k}} (\WCal_p(E_{\#}^\star\widehat{\mu}_n, E_{\#}^\star\mu_\star) - \EE[\WCal_p(E_{\#}^\star\widehat{\mu}_n, E_{\#}^\star\mu_\star)])]$ under different condition.  

We first claim that $\WCal_p(E_{\#}^\star\widehat{\mu}_n, E_{\#}^\star\mu_\star) - \EE[\WCal_p(E_{\#}^\star\widehat{\mu}_n, E_{\#}^\star\mu_\star)]$ is sub-exponential with parameters $(\sqrt{K/2} n^{-1/(2\vee p)}, (K/2)n^{-1/p})$ for all $E \in \bs_{d, k}$ if the true measure $\mu_\star$ satisfies the projection Poincar\'{e} inequality (cf. Definition~\ref{def:proj-Poincare}). Indeed, we consider $X = (X_1, \ldots, X_n)$ and $X' = (X'_1, \ldots, X'_n)$ where $X_i, X'_i$ are independent samples from $E_{\#}^\star\mu_\star$. Let $f(X)=\WCal_p(E_{\#}^\star\widehat{\mu}_n, E_{\#}^\star\mu_\star)$, we have $\EE(|f(X)|) < +\infty$. By the triangle inequality, we have  
\begin{equation*}
|f(X) - f(X')| \leq n^{-1/p}\left(\sum_{i=1}^n \|X_i - X'_i\|^p\right)^{1/p} \leq n^{-\frac{1}{2\vee p}}\|X - X'\|.  
\end{equation*}
This implies that the following statement holds almost surely,  
\begin{equation*}
\sum_{i=1}^n \|\nabla_i f(X)\|^2 \leq n^{-\frac{2}{2\vee p}} \quad \textnormal{and} \quad \max_{1 \leq i \leq n} \|\nabla_i f(X)\| \leq n^{-\frac{1}{p}}, \quad \textnormal{almost surely}.  
\end{equation*}
In addition, the probability measure $E_{\#}^\star\mu_\star \in \PScr(\br^k)$ is assumed to satisfy the Poincar\'{e} inequality. Equipped with Theorem~\ref{Theorem:appendix-Poincare} yields that 
\begin{equation*}
\PP\left(\WCal_p(E_{\#}^\star\widehat{\mu}_n, E_{\#}^\star\mu_\star) - \EE[\WCal_p(E_{\#}^\star\widehat{\mu}_n, E_{\#}^\star\mu_\star)] \geq t\right) \leq \exp\left(-\frac{1}{K}\min\left\{\frac{t}{n^{-1/p}}, \ \frac{t^2}{n^{-2/(2\vee p)}}\right\}\right),   
\end{equation*}
For the simplicity, let $Z_E = \WCal_p(E_{\#}^\star\widehat{\mu}_n, E_{\#}^\star\mu_\star) - \EE[\WCal_p(E_{\#}^\star\widehat{\mu}_n, E_{\#}^\star\mu_\star)]$. Then we have $\EE[Z_E] = 0$ and $\PP(Z_E \geq t) \leq \exp(-K^{-1}\min\{n^{1/p}t, n^{2/(2\vee p)}t^2\})$. This together with the definition of $Z_E$ and~\citet[Theorem~2.2]{Wainwright-2019-High} yields the desired claim. 

Using the same argument in Theorem~\ref{Theorem:PRW-main-Bernstein}, we can interpret $\{Z_E\}_{E \in \bs_{d, k}}$ as an empirical process indexed by $E \in \bs_{d, k}$ and show that there exists a random variable $L$ satisfying $\EE[L] \leq 4M_q(\mu_\star)$ so that $|Z_U - Z_V| \leq L\|U - V\|_\op$ for all $U, V \in \bs_{d, k}$. By a standard $\epsilon$-net argument, we obtain that 
\begin{equation*}
\EE[\sup_{E \in \bs_{d, k}} Z_E] \leq \inf_{\epsilon > 0}\left\{\epsilon\EE[L] + \sqrt{2K} n^{-1/(2 \vee p)}\sqrt{\log(N(\bs_{d, k}, \epsilon, \|\cdot\|_\op))} + (K/2)n^{-1/p}\log(N(\bs_{d, k}, \epsilon, \|\cdot\|_\op))\right\}. 
\end{equation*}
Combining Proposition~\ref{Proposition:appendix-covering} and choosing $\epsilon = \sqrt{k}n^{-1/p}$ (it is chosen to achieve the tight bound) yields that 
\begin{eqnarray*}
\EE\left[\sup_{E \in \bs_{d, k}} Z_E \right] & \lesssim_{p, q} & \inf_{\epsilon > 0}\left\{\epsilon + n^{-1/(2 \vee p)}\sqrt{dk\log\left(\frac{\sqrt{k}}{\epsilon}\right)} + n^{-1/p}dk\log\left(\frac{\sqrt{k}}{\epsilon}\right)\right\} \\
& \lesssim_{p, q} & n^{-1/(2 \vee p)}\sqrt{dk\log(n)} + n^{-1/p}dk\log(n). 
\end{eqnarray*}
Therefore, we conclude that 
\begin{equation*}
\EE\left[\sup_{E \in \bs_{d, k}} \left(\WCal_p(E_{\#}^\star\widehat{\mu}_n, E_{\#}^\star\mu_\star) - \EE[\WCal_p(E_{\#}^\star\widehat{\mu}_n, E_{\#}^\star\mu_\star)]\right)\right] \lesssim_{p, q} n^{-1/(2 \vee p)}\sqrt{dk\log(n)} + n^{-1/p}dk\log(n).   
\end{equation*}
This together with Eq.~\eqref{inequality-PRW-main-first} and Eq.~\eqref{inequality-PRW-main-second} yields the desired inequality. 

\subsection{Proof of Theorem~\ref{Theorem:concentration-bernstein}}
Since the arguments in this proof hold true for both IPRW and PRW distances, we denote $W = \underline{\PWCal}_{p, k}$ or $W = \overline{\PWCal}_{p, k}$ for short. Let $f(X) = W(\widehat{\mu}_n, \mu_\star)$, we have
\begin{equation*}
D_i = f(X) - f(X'_{(i)}) \leq W(\widehat{\mu}_n, \widehat{\mu}'_n) \leq n^{-1/p}\left(\sup_{E \in \bs_{d, k}} \|E_{\#}^\star(X_i) - E_{\#}^\star(X'_i)\|\right). 
\end{equation*}
By the triangle inequality, we have
\begin{equation*}
\EE\left[|D_i|^k \mid X_{-i}\right] \leq 2^k n^{-k/p}\left(\EE\left[\sup_{E \in \bs_{d, k}, X \sim E_{\#}^\star\mu} |X|^k\right]\right). 
\end{equation*}
Since the true measure $\mu_\star$ satisfies the Bernstein-type tail condition (cf. Definition~\ref{def:Bernstein}), we have
\begin{equation*}
\EE\left[|D_i|^k \mid X_{-i}\right] \leq 2^{k-1} n^{-k/p} \sigma^2 k! V^{k-2} = \frac{(2n^{-1/p}\sigma)^2k!(2n^{-1/p}V)^{k-2}}{2} 
\end{equation*}
This implies that the condition in Theorem~\ref{Theorem:appendix-BM} holds true with $\sigma_i = 2n^{-1/p}\sigma$ and $M = 2n^{-1/p}V$. Equipped with Theorem~\ref{Theorem:appendix-BM} yields the desired inequality. 

\subsection{Proof of Theorem~\ref{Theorem:concentration-poincare}}
Since the arguments in this proof hold true for both IPRW and PRW distances, we denote $W = \underline{\PWCal}_{p, k}$ or $W = \overline{\PWCal}_{p, k}$ for short. We consider $X = (X_1, X_2, \ldots, X_n)$ and $X' = (X'_1, X'_2, \ldots, X'_n)$ where $X_i, X'_i$ are independent samples from $\mu_\star$. Let $f(X) = W(\widehat{\mu}_n, \mu_\star)$, we have $\EE(|f(X)|) < +\infty$. By the triangle inequality, we have 
\begin{equation*}
|f(X) - f(X')| \leq n^{-1/p}\left(\sum_{i=1}^n \|X_i-X'_i\|^p\right)^{1/p} \leq n^{-\frac{1}{2\vee p}}\|X-X'\|.  
\end{equation*}
This implies that the following statement holds almost surely, 
\begin{equation*}
\sum_{i=1}^n \|\nabla_i f(X)\|^2 \leq n^{-\frac{2}{2\vee p}} \quad \textnormal{and} \quad \max_{1 \leq i \leq n} \|\nabla_i f(X)\| \leq n^{-\frac{1}{p}}.  
\end{equation*}
In addition, the true measure $\mu_\star$ satisfies the Poincar\'{e} inequality (cf. Definition~\ref{def:Poincare}). Equipped with Theorem~\ref{Theorem:appendix-Poincare} yields the desired inequality. 

\section{Postponed Proofs in Subsection~\ref{subsec:estimators}}
In this section, we provide the detailed proofs for Theorem~\ref{Theorem:consistency-MPRW}-\ref{Theorem:convergence-MEPRW-MPRW} and Theorem~\ref{Theorem:measurability-MPRW}-\ref{Theorem:measurability-MEPRW}. Our results are derived analogously to the proof in~\citet{Bernton-2019-Parameter} for the estimators based on Wasserstein distance and the proof in~\citet{Nadjahi-2019-Asymptotic} for the estimators based on sliced-Wasserstein distance.

\subsection{Preliminary technical results}
To facilitate the reading, we collect several preliminary technical results which will be used in the postponed proofs in subsection~\ref{subsec:estimators}.
\begin{theorem}[Theorem~2.43 in~\citet{Aliprantis-2006-Infinite}]\label{Theorem:appendix-compact}
A real-valued lower semi-continuous function on a compact space attains a minimum value, and the nonempty set of minimizers is compact. Similarly, an upper semicontinuous function on a compact set attains a maximum value, and the nonempty set of maximizers is compact.
\end{theorem}
\begin{definition}[epiconvergence]
Let $\XCal$ be a metric space and $\{f_i\}_{i \in \bn}$ be a sequence of real-valued function from $\XCal$ to $\br$. We say that the sequence $\{f_i\}_{i \in \bn}$ epiconverges to a function $f: \XCal \rightarrow\ br$ if for each $x\in \XCal$, the following statement holds true, 
\begin{align*}
\liminf_{i \rightarrow +\infty} f_i(x_i) \geq f(x) & \textnormal{ for every sequence } \{x_i\}_{i \in \bn} \textnormal{ such that } x_i \rightarrow x, \\
\limsup_{i \rightarrow +\infty} f_i(x_i) \leq f(x) & \textnormal{ for some sequence } \{x_i\}_{i \in \bn} \textnormal{ such that } x_i\rightarrow x.
\end{align*}
\end{definition}
\begin{proposition}[Proposition~7.29 in~\citet{Rockafellar-2009-Variational}]\label{Proposition:appendix-epi}
Let $\XCal$ be a metric space and $\{f_i\}_{i \in \bn}$ be a sequence of real-valued function from $\XCal$ to $\br$ with a lower semi-continuous function $f: \XCal \rightarrow \br$. Then the sequence $\{f_i\}_{i \in \bn}$ epiconverges to $f$ if and only if 
\begin{align*}
\liminf_{i \rightarrow +\infty} (\inf_{x \in K} f_i(x)) \geq \inf_{x \in K} f(x) & \textnormal{ for every compact set } K \subseteq \XCal, \\
\limsup_{i \rightarrow +\infty} (\sup_{x \in O} f_i(x)) \leq \sup_{x \in O} f(x) & \textnormal{ for every open set } O \subseteq \XCal.
\end{align*}
\end{proposition}
Recall that $\delta$-$\argmin_{x \in \XCal} f = \{x \in \XCal: f(x) \leq \inf_{x \in \XCal} f + \delta\}$ for a generic function $f: \XCal \rightarrow \br$. The following theorem gives asymptotic properties for the infimum and $\delta$-argmin of epiconvergent functions and thus a standard approach to prove the existence and consistency of the estimators. 
\begin{theorem}[Theorem~7.31 in~\citet{Rockafellar-2009-Variational}]\label{Theorem:appendix-final}
Let $\XCal$ be a metric space and $\{f_i\}_{i \in \bn}$ be a sequence of function which epiconverges to a lower semi-continuous function $f$ with $\inf_{x \in \XCal} f \in (-\infty, +\infty)$. Then we have the following statements, 
\begin{enumerate}
\item $\inf_{x \in \XCal} f_i \rightarrow \inf_{x \in \XCal} f$ if and only if for every $\delta > 0$ there exists a compact set $B \subseteq \XCal$ and $N \in \bn$ such that $\inf_{x \in B} f_i \leq \inf_{x \in \XCal} f_i + \delta$ for all $i \geq N$. 
\item $\limsup_{i \rightarrow +\infty} (\delta\textnormal{-}\argmin_{x \in \XCal} f_i) \subseteq \delta\textnormal{-}\argmin_{x \in \XCal} f$ for any $\delta \geq 0$ and $\limsup_{i \rightarrow +\infty} (\delta_i\textnormal{-}\argmin_{x \in \XCal} f_i) \subseteq \argmin_{x \in \XCal} f$ whenever $\delta_i \downarrow 0$.  
\item Assume that $\inf_{x \in \XCal} f_i \rightarrow \inf_{x \in \XCal} f$, there exists a sequence $\delta_i \downarrow 0$ such that $\delta_i \textnormal{-}\argmin_{x \in \XCal} f_i \rightarrow \argmin_{x \in \XCal} f$. Conversely, if $\argmin_{x \in \XCal} f \neq \emptyset$ and if such a sequence exists, then $\inf_{x \in \XCal} f_i \rightarrow \inf_{x \in \XCal} f$. 
\end{enumerate}
\end{theorem}
The following theorem summarizes the well-known Skorokhod's representation theorem.
\begin{theorem}[Skorokhod's representation theorem]\label{Theorem:skorokhod}
Let $\{\mu_n\}_{n \in \bn}$ be a sequence of probability measures on a metric space $\SCal$ such that $\mu_n$ converges weakly to some probability measure $\mu_\infty$ on $\SCal$ as $n \rightarrow \infty$. Suppose also that the support of $\mu_\infty$ is separable. Then there exist random variables $X_n$ defined on a common probability space $(\Omega, \FCal, \PP)$ such that the law of $X_n$ is $\mu_n$ for all $n$ (including $n = \infty$) and such that $X_n$ converges to $X_\infty$ almost surely.
\end{theorem} 
The following theorem presents the classical results which lead to a standard approach for proving the measurability of the estimators. Note that the projection $\proj(D) = \{x \in \XCal: \exists y \in \YCal, \st (x, y) \in D\}$ for each $D \subseteq \XCal \times \YCal$ and the section $D_x = \{y \in \YCal: (x, y) \in D\}$ for each $x \in \proj(D)$. 
\begin{theorem}[Corollary~1 in~\citet{Brown-1973-Measurable}]\label{Theorem:brown-1973}
Let $\XCal, \YCal$ be complete separable metric spaces and $f$ be a real-valued Borel measurable function defined on a Borel subset $D$ of $\XCal \times \YCal$. Suppose that for each $x \in \proj(D)$, the section $D_x$ is $\sigma$-compact and $f(x, \cdot)$ is lower semi-continuous with respect to the relative topology on $D_x$. Then 
\begin{enumerate}
\item The sets $G = \proj(D)$ and $I = \{x \in G: \exists y \in D_x \ \st \ y = \argmin_{z \in \YCal} f(x, z)\}$ are Borel. 
\item For each $\epsilon > 0$, there exists a Borel measure function $\varphi_\epsilon$ satisfying, for $x \in G$ that, 
\begin{equation*}
f(x, \varphi_\epsilon(x)) \left\{\begin{array}{lcl}
= \inf_{y \in G} f(x, y), & & x \in I, \\
\leq \epsilon + \inf_{y \in G} f(x, y), & \textnormal{if} & x \notin I \textnormal{ and } \inf_{y \in G} f(x, y) \neq -\infty, \\
\leq -\epsilon^{-1}, & & x \notin I \textnormal{ and } \inf_{y \in G} f(x, y) = -\infty.  
\end{array}\right. 
\end{equation*}
\end{enumerate}
\end{theorem}
To show that the MEPRW estimator is measurable, we establish the lower semi-continuity of the expectation of empirical PRW distance in the following lemma. 
\begin{lemma}\label{lemma:EPRW-lsc}
The expected empirical PRW distance is lower semi-continuous in the usual weak topology. If the sequences $\{\mu_i\}_{i \in \bn}, \{\nu_i\}_{i \in \bn} \subseteq \PScr(\br^d)$ satisfying that $\mu_i \Rightarrow \mu \in \PScr(\br^d)$ and $\nu_i \Rightarrow \nu \in \PScr(\br^d)$, we have $\EE[\overline{\PWCal}_{p, k}(\mu, \widehat{\nu}_m)] \leq \liminf_{i \rightarrow +\infty} \EE[\overline{\PWCal}_{p, k}(\mu_i, \widehat{\nu}_{i, m})]$, where $\widehat{\nu}_m = (1/m)\sum_{j=1}^m \delta_{Z_j}$ for i.i.d. samples $Z_{1:m}$ according to $\nu$ and $\{\widehat{\nu}_{i, m}\}_{i \in \bn}$ are defined similarly.  
\end{lemma}

\subsection{Proof of Theorem~\ref{Theorem:consistency-MPRW}}
We first prove that $\argmin_{\theta \in \Theta} \overline{\PWCal}_{p, k}(\mu_\star, \mu_\theta) \neq \emptyset$. Indeed, by Assumption~\ref{A2} and Theorem~\ref{Theorem:PRW-lsc}, the mapping $\theta \mapsto \overline{\PWCal}_{p, k}(\mu_\star, \mu_\theta)$ is lower semi-continuous. By Assumption~\ref{A3}, the set $\Theta_\star(\tau)$ is bounded for some $\tau > 0$. By the definition of $\inf$, there exists $\theta' \in \Theta$ such that $\overline{\PWCal}_{p, k}(\mu_\star, \mu_{\theta'}) = \inf_{\theta \in \Theta} \overline{\PWCal}_{p, k}(\mu_\star, \mu_\theta) + \tau/2$. This implies that $\theta' \in \Theta_\star(\tau)$ and $\Theta_\star(\tau)$ is nonempty. By the lower semi-continuity of the mapping $\theta \mapsto \overline{\PWCal}_{p, k}(\mu_\star, \mu_\theta)$, the set $\Theta_\star(\tau)$ is closed. Putting these pieces together yields that $\Theta_\star(\tau)$ is compact. Therefore, we conclude the desired result from Theorem~\ref{Theorem:appendix-compact}. 

Then we show that there exists a set $E \subseteq \Omega$ with $\PP(E)=1$ such that, for all $\omega \in E$, the sequence of mappings $\theta \mapsto \overline{\PWCal}_{p, k}(\widehat{\mu}_n(\omega), \mu_\theta)$ epiconverges to the mapping $\theta \mapsto \overline{\PWCal}_{p, k}(\mu_\star, \mu_\theta)$ as $n \rightarrow +\infty$. Indeed, we only need to prove that the conditions in Proposition~\ref{Proposition:appendix-epi} hold true.  

Fix $K \subseteq \Theta$ as a compact set. By the lower semi-continuity of the mapping $\theta \mapsto \overline{\PWCal}_{p, k}(\widehat{\mu}_n(\omega), \mu_\theta)$ (cf. Assumption~\ref{A2} and Theorem~\ref{Theorem:PRW-lsc}), Theorem~\ref{Theorem:appendix-compact} implies that 
\begin{equation*}
\inf_{\theta \in K} \overline{\PWCal}_{p, k}(\widehat{\mu}_n(\omega), \mu_\theta) = \overline{\PWCal}_{p, k}(\widehat{\mu}_n(\omega), \mu_{\theta_n})
\end{equation*}
for some sequence $\theta_n = \theta_n(\omega) \in K$. Thus, we have
\begin{equation*}
\liminf_{n \rightarrow +\infty} \inf_{\theta \in K} \overline{\PWCal}_{p, k}(\widehat{\mu}_n(\omega), \mu_\theta) = \liminf_{n \rightarrow +\infty} \overline{\PWCal}_{p, k}(\widehat{\mu}_n(\omega), \mu_{\theta_n}). 
\end{equation*}
By the definition of $\liminf$, there exists a subsequence of $\{\theta_n\}_{n \in \bn}$ such that $\overline{\PWCal}_{p, k}(\widehat{\mu}_n(\omega), \mu_{\theta_n})$ converges to $\liminf_{n \rightarrow +\infty} \overline{\PWCal}_{p, k}(\widehat{\mu}_n(\omega), \mu_{\theta_n})$ along this subsequence. By the compactness of $K$, this subsequence must have a convergent subsubsequence. We denote this subsubsequence as $\{\theta_{n_j}\}_{j \in \bn}$ and its limit as $\bar{\theta} \in K$. Then 
\begin{equation*}
\liminf_{n \rightarrow +\infty} \overline{\PWCal}_{p, k}(\widehat{\mu}_n(\omega), \mu_{\theta_n}) = \lim_{j \rightarrow +\infty} \overline{\PWCal}_{p, k}(\widehat{\mu}_{n_j}(\omega), \mu_{\theta_{n_j}}). 
\end{equation*}
Since $\omega \in E$ where $\PP(E)=1$, Assumption~\ref{A1} and~\ref{A2} imply $\widehat{\mu}_{n_j}(\omega) \Rightarrow \mu_\star$ and $\mu_{\theta_{n_j}} \Rightarrow \mu_{\bar{\theta}}$. These pieces together with the lower semi-continuity of the PRW distance (cf. Theorem~\ref{Theorem:PRW-lsc}) yields that $\lim_{j \rightarrow +\infty} \overline{\PWCal}_{p, k}(\widehat{\mu}_{n_j}(\omega), \mu_{\theta_{n_j}}) \geq \overline{\PWCal}_{p, k}(\mu_\star, \mu_{\bar{\theta}})$. Putting these pieces together yields that 
\begin{equation*}
\liminf_{n \rightarrow +\infty} \inf_{\theta \in K} \overline{\PWCal}_{p, k}(\widehat{\mu}_n(\omega), \mu_\theta) \geq \inf_{\theta \in K} \overline{\PWCal}_{p, k}(\mu_\star, \mu_\theta). 
\end{equation*}
Fix $O \subseteq \Theta$ as an arbitary open set. By the definition of $\inf$, there exists a sequence $\theta'_n = \theta'_n(\omega) \in O$ such that $\overline{\PWCal}_{p, k}(\mu_\star, \mu_{\theta'_n}) \rightarrow \inf_{\theta \in O} \overline{\PWCal}_{p, k}(\mu_\star, \mu_\theta)$. In addition, $\inf_{\theta \in O} \overline{\PWCal}_{p, k}(\widehat{\mu}_n(\omega), \mu_\theta) \leq \overline{\PWCal}_{p, k}(\widehat{\mu}_n(\omega), \mu_{\theta'_n})$. Thus, we have
\begin{eqnarray*}
\limsup_{n \rightarrow +\infty} \inf_{\theta \in O} \overline{\PWCal}_{p, k}(\widehat{\mu}_n(\omega), \mu_\theta) & \leq & \limsup_{n \rightarrow +\infty} \overline{\PWCal}_{p, k}(\widehat{\mu}_n(\omega), \mu_{\theta'_n})  \\
& & \hspace*{-8em} \leq \ \limsup_{n \rightarrow +\infty} \overline{\PWCal}_{p, k}(\widehat{\mu}_n(\omega), \mu_\star) + \limsup_{n \rightarrow +\infty} \overline{\PWCal}_{p, k}(\mu_\star, \mu_{\theta'_n}).  
\end{eqnarray*}
Since $\omega \in E$ where $\PP(E)=1$, Assumption~\ref{A1} implies $\limsup_{n \rightarrow +\infty} \overline{\PWCal}_{p, k}(\widehat{\mu}_n(\omega), \mu_\star) = 0$. By the definition of $\theta'_n$, $\limsup_{n \rightarrow +\infty} \overline{\PWCal}_{p, k}(\mu_\star, \mu_{\theta'_n})=\inf_{\theta \in O} \overline{\PWCal}_{p, k}(\mu_\star, \mu_\theta)$. Putting these pieces together yields that $\limsup_{n \rightarrow +\infty} \inf_{\theta \in O} \overline{\PWCal}_{p, k}(\widehat{\mu}_n(\omega), \mu_\theta) \leq \inf_{\theta \in O} \overline{\PWCal}_{p, k}(\mu_\star, \mu_\theta)$. 

Proposition~\ref{Proposition:appendix-epi} guarantees that there exists a set $E \subseteq \Omega$ with $\PP(E)=1$ such that, for all $\omega \in E$, the sequence of mappings $\theta \mapsto \overline{\PWCal}_{p, k}(\widehat{\mu}_n(\omega), \mu_\theta)$ epiconverges to the mapping $\theta \mapsto \overline{\PWCal}_{p, k}(\mu_\star, \mu_\theta)$ as $n \rightarrow +\infty$. Then the second statement of Theorem~\ref{Theorem:appendix-final} implies that 
\begin{equation}\label{result:consistency-MPRW-first}
\limsup_{n \rightarrow +\infty} \argmin_{\theta \in \Theta} \overline{\PWCal}_{p, k}(\widehat{\mu}_n(\omega), \mu_\theta) \subseteq \argmin_{\theta \in \Theta} \overline{\PWCal}_{p, k}(\mu_\star, \mu_\theta).  
\end{equation}
The next step is to show that, for every $\delta > 0$, there exists a compact set $B \subseteq \Theta$ and $N \in \bn$ such that $\inf_{\theta \in B} \overline{\PWCal}_{p, k}(\widehat{\mu}_n(\omega), \mu_\theta) \leq \inf_{\theta \in \Theta} \overline{\PWCal}_{p, k}(\widehat{\mu}_n(\omega), \mu_\theta) + \delta$. In what follows, we prove a stronger statement which states that the above inequality holds true with $\delta = 0$. Indeed, by the same reasoning for the open set case in the proof of epiconvergence, we have
\begin{equation*}
\limsup_{n \rightarrow +\infty} \inf_{\theta \in \Theta} \overline{\PWCal}_{p, k}(\widehat{\mu}_n(\omega), \mu_\theta) \leq \inf_{\theta \in \Theta} \overline{\PWCal}_{p, k}(\mu_\star, \mu_\theta). 
\end{equation*}
By Assumption~\ref{A3} and using previous argument, $\Theta_\star(\tau)$ is nonempty and compact for some $\tau > 0$. The above inequality implies that there exists $n_1(\omega) > 0$ such that, for all $n \geq n_1(\omega)$, the set $\{\theta \in \Theta: \overline{\PWCal}_{p, k}(\widehat{\mu}_n(\omega), \mu_\theta) \leq \inf_{\theta' \in \Theta} \overline{\PWCal}_{p, k}(\mu_\star, \mu_{\theta'}) + \tau/2\}$ is nonempty. For any $\theta$ in this set and let $n \geq n_1(\omega)$, we have
\begin{equation*}
\overline{\PWCal}_{p, k}(\mu_\star, \mu_\theta) \leq \overline{\PWCal}_{p, k}(\mu_\star, \widehat{\mu}_n(\omega)) + \inf_{\theta \in \Theta} \overline{\PWCal}_{p, k}(\mu_\star, \mu_\theta) + \frac{\tau}{2}. 
\end{equation*}
By Assumption~\ref{A1}, there exists $n_2(\omega) > 0$ such that, for all $n \geq n_2(\omega)$, we have
\begin{equation*}
\overline{\PWCal}_{p, k}(\mu_\star, \widehat{\mu}_n(\omega)) \leq \WCal_p(\mu_\star, \widehat{\mu}_n(\omega)) \leq \frac{\tau}{2}. 
\end{equation*}
Putting these pieces together yields that, for all $n \geq \max\{n_1(\omega), n_2(\omega)\}$, we have $\overline{\PWCal}_{p, k}(\mu_\star, \mu_\theta) \leq \inf_{\theta \in \Theta} \overline{\PWCal}_{p, k}(\mu_\star, \mu_\theta)  + \tau$. This implies that, for all $n \geq \max\{n_1(\omega), n_2(\omega)\}$ that, 
\begin{equation*}
\left\{\theta \in \Theta: \overline{\PWCal}_{p, k}(\widehat{\mu}_n(\omega), \mu_\theta) \leq \inf_{\theta' \in \Theta} \overline{\PWCal}_{p, k}(\mu_\star, \mu_{\theta'}) + \frac{\tau}{2}\right\} \subseteq \Theta_\star(\tau). 
\end{equation*}
Therefore, we have $\inf_{\theta \in \Theta} \overline{\PWCal}_{p, k}(\widehat{\mu}_n(\omega), \mu_\theta) = \inf_{\theta \in \Theta_\star(\tau)} \overline{\PWCal}_{p, k}(\widehat{\mu}_n(\omega), \mu_\theta)$. This together with the compactness of  $\Theta_\star(\tau)$ yields the desired result. 

The first statement of Theorem~\ref{Theorem:appendix-final} implies that
\begin{equation}\label{result:consistency-MPRW-second}
\inf_{\theta \in \Theta} \overline{\PWCal}_{p, k}(\widehat{\mu}_n(\omega), \mu_\theta) \rightarrow \inf_{\theta \in \Theta} \overline{\PWCal}_{p, k}(\mu_\star, \mu_\theta), \quad \textnormal{as } n \rightarrow +\infty. 
\end{equation}
By Assumption~\ref{A2} and Theorem~\ref{Theorem:PRW-lsc}, the mapping $\theta \mapsto \overline{\PWCal}_{p, k}(\widehat{\mu}_n(\omega), \mu_\theta)$ is lower semi-continuous. Theorem~\ref{Theorem:appendix-compact} implies $\argmin_{\theta \in \Theta} \overline{\PWCal}_{p, k}(\widehat{\mu}_n(\omega), \mu_\theta)$ are nonempty for all $n \geq \max\{n_1(\omega), n_2(\omega)\}$. Together with Eq.~\eqref{result:consistency-MPRW-first} and~\eqref{result:consistency-MPRW-second} yields the desired results. 

Finally, we remark that these results hold true for $\delta_n$-$\argmin_{\theta \in \Theta} \overline{\PWCal}_{p, k}(\widehat{\mu}_n, \mu_\theta)$ with $\delta_n \rightarrow 0$. For Eq.~\eqref{result:consistency-MPRW-first} and~\eqref{result:consistency-MPRW-second}, the analogous results can be derived by using the second and third statements of Theorem~\ref{Theorem:appendix-final}. To show that $\delta_n$-$\argmin_{\theta \in \Theta} \overline{\PWCal}_{p, k}(\widehat{\mu}_n, \mu_\theta)$ is nonempty, we notice it contains the nonempty set $\argmin_{\theta \in \Theta} \overline{\PWCal}_{p, k}(\widehat{\mu}_n, \mu_\theta)$.

\subsection{Proof of Theorem~\ref{Theorem:consistency-MEPRW}}
Following up the same approach used for analyzing Theorem~\ref{Theorem:consistency-MPRW}, it is straightforward to derive that $\argmin_{\theta \in \Theta} \overline{\PWCal}_{p, k}(\mu_\star, \mu_\theta) \neq \emptyset$. Then we show that there exists a set $E \subseteq \Omega$ with $\PP(E)=1$ such that, for all $\omega \in E$, the sequences $\theta \mapsto \EE[\overline{\PWCal}_{p, k}(\widehat{\mu}_n(\omega), \widehat{\mu}_{\theta, m(n)}) \mid X_{1:n}]$ epiconverges $\theta \mapsto \overline{\PWCal}_{p, k}(\mu_\star, \mu_\theta)$ as $n \rightarrow +\infty$. Indeed, it suffices to verify the conditions in Proposition~\ref{Proposition:appendix-epi}. 

Fix $K \subseteq \Theta$ as an arbitrary compact set. By Assumption~\ref{A2} and Lemma~\ref{lemma:EPRW-lsc}, the mapping $\theta \mapsto \EE[\overline{\PWCal}_{p, k}(\widehat{\mu}_n(\omega), \widehat{\mu}_{\theta, m(n)}) \mid X_{1:n}]$ is lower semi-continuous. Then Theorem~\ref{Theorem:appendix-compact} implies that 
\begin{equation*}
\inf_{\theta \in K} \EE[\overline{\PWCal}_{p, k}(\widehat{\mu}_n(\omega), \widehat{\mu}_{\theta, m(n)}) \mid X_{1:n}] = \EE[\overline{\PWCal}_{p, k}(\widehat{\mu}_n(\omega), \widehat{\mu}_{\theta_n, m(n)}) \mid X_{1:n}]
\end{equation*}
for some sequence $\theta_n = \theta_n(\omega) \in K$. Thus, we have
\begin{equation*}
\liminf_{n \rightarrow +\infty} \inf_{\theta \in K} \EE\left[\overline{\PWCal}_{p, k}(\widehat{\mu}_n(\omega), \widehat{\mu}_{\theta, m(n)}) \mid X_{1:n}\right] = \liminf_{n \rightarrow +\infty} \EE\left[\overline{\PWCal}_{p, k}(\widehat{\mu}_n(\omega), \widehat{\mu}_{\theta_n, m(n)}) \mid X_{1:n}\right]. 
\end{equation*}
Following up the same approach used in the proof of Theorem~\ref{Theorem:consistency-MPRW}, there exists a subsequence of $\{\theta_n\}_{n \in \bn}$, denoted by $\{\theta_{n_j}\}_{j \in \bn}$ with the limit $\bar{\theta} \in K$, such that 
\begin{eqnarray*}
\liminf_{n \rightarrow +\infty} \EE[\overline{\PWCal}_{p, k}(\widehat{\mu}_n(\omega), \widehat{\mu}_{\theta_n, m(n)}) \mid X_{1:n}] & = & \lim_{j \rightarrow +\infty} \EE[\overline{\PWCal}_{p, k}(\widehat{\mu}_{n_j}(\omega), \widehat{\mu}_{\theta_{n_j}, m(n_j)}) \mid X_{1:n_j}] \\
& & \hspace*{-16em} \geq \ \liminf_{j \rightarrow +\infty} \EE[\overline{\PWCal}_{p, k}(\widehat{\mu}_{n_j}(\omega), \mu_{\theta_{n_j}})] - \limsup_{j \rightarrow +\infty} \EE[\overline{\PWCal}_{p, k}(\mu_{\theta_{n_j}}, \widehat{\mu}_{\theta_{n_j}, m(n_j)}) \mid X_{1:n_j}]. 
\end{eqnarray*}
Since $\omega \in E$ where $\PP(E)=1$, Assumption~\ref{A1} and~\ref{A2} imply $\widehat{\mu}_{n_j}(\omega) \Rightarrow \mu_\star$ and $\mu_{\theta_{n_j}} \Rightarrow \mu_{\bar{\theta}}$. These pieces together with the lower semi-continuity of the PRW distance (cf. Theorem~\ref{Theorem:PRW-lsc}) yields that $\liminf_{j \rightarrow +\infty} \overline{\PWCal}_{p, k}(\widehat{\mu}_{n_j}(\omega), \mu_{\theta_{n_j}}) \geq \overline{\PWCal}_{p, k}(\mu_\star, \mu_{\bar{\theta}})$. By Assumption~\ref{A4} and using $\theta_{n_j} \rightarrow \bar{\theta}$, we have $\limsup_{j \rightarrow +\infty} \EE[\overline{\PWCal}_{p, k}(\mu_{\theta_{n_j}}, \widehat{\mu}_{\theta_{n_j}, m(n_j)}) \mid X_{1:n_j}] \rightarrow 0$. Putting these pieces together yields that 
\begin{equation*}
\liminf_{n \rightarrow +\infty} \inf_{\theta \in K} \EE[\overline{\PWCal}_{p, k}(\widehat{\mu}_n(\omega), \widehat{\mu}_{\theta, m(n)}) \mid X_{1:n}] \geq \inf_{\theta \in K} \overline{\PWCal}_{p, k}(\mu_\star, \mu_\theta). 
\end{equation*}
Fix $O \subseteq \Theta$ as an arbitary open set. By the definition of $\inf$, there exists a sequence $\theta'_n = \theta'_n(\omega) \in O$ such that $\overline{\PWCal}_{p, k}(\mu_\star, \mu_{\theta'_n}) \rightarrow \inf_{\theta \in O} \overline{\PWCal}_{p, k}(\mu_\star, \mu_\theta)$. In addition, we have 
\begin{equation*}
\inf_{\theta \in O} \EE[\overline{\PWCal}_{p, k}(\widehat{\mu}_n(\omega), \widehat{\mu}_{\theta, m(n)}) \mid X_{1:n}] \leq \EE[\overline{\PWCal}_{p, k}(\widehat{\mu}_n(\omega), \widehat{\mu}_{\theta'_n, m(n)}) \mid X_{1:n}]. 
\end{equation*}
Thus, we have
\begin{eqnarray*}
& & \limsup_{n \rightarrow +\infty} \inf_{\theta \in O} \EE[\overline{\PWCal}_{p, k}(\widehat{\mu}_n(\omega), \widehat{\mu}_{\theta, m(n)}) \mid X_{1:n}] \ \leq \ \limsup_{n \rightarrow +\infty} \EE[\overline{\PWCal}_{p, k}(\widehat{\mu}_n(\omega), \widehat{\mu}_{\theta'_n, m(n)}) \mid X_{1:n}] \\
& \leq & \limsup_{n \rightarrow +\infty} \overline{\PWCal}_{p, k}(\widehat{\mu}_n(\omega), \mu_\star) + \limsup_{n \rightarrow +\infty} \overline{\PWCal}_{p, k}(\mu_\star, \mu_{\theta'_n}) + \limsup_{n \rightarrow +\infty} \EE[\overline{\PWCal}_{p, k}(\mu_{\theta'_n}, \widehat{\mu}_{\theta'_n, m(n)}) \mid X_{1:n}]. 
\end{eqnarray*}
Since $\omega \in E$ where $\PP(E)=1$, Assumption~\ref{A1} implies $\limsup_{n \rightarrow +\infty} \overline{\PWCal}_{p, k}(\widehat{\mu}_n(\omega), \mu_\star) = 0$. By the definition of $\theta'_n$, we have $\limsup_{n \rightarrow +\infty} \overline{\PWCal}_{p, k}(\mu_\star, \mu_{\theta'_n})=\inf_{\theta \in O} \overline{\PWCal}_{p, k}(\mu_\star, \mu_\theta)$. Using Assumption~\ref{A4} and $\lim_{j \rightarrow +\infty} \theta_{m_j}=\bar{\theta}$, we have $\limsup_{n \rightarrow +\infty} \EE[\overline{\PWCal}_{p, k}(\mu_{\theta'_n}, \widehat{\mu}_{\theta'_n, m(n)}) \mid X_{1:n}] = 0$. Putting these pieces together yields that $\limsup_{n \rightarrow +\infty} \inf_{\theta \in O} \EE[\overline{\PWCal}_{p, k}(\widehat{\mu}_n(\omega), \widehat{\mu}_{\theta, m(n)}) \mid X_{1:n}] \leq \inf_{\theta \in O} \overline{\PWCal}_{p, k}(\mu_\star, \mu_\theta)$. 

Proposition~\ref{Proposition:appendix-epi} guarantees that there exists a set $E \subseteq \Omega$ with $\PP(E)=1$ such that, for all $\omega \in E$, the sequence of mappings $\theta \mapsto \EE[\overline{\PWCal}_{p, k}(\widehat{\mu}_n(\omega), \widehat{\mu}_{\theta, m(n)}) \mid X_{1:n}]$ epiconverges to the mapping $\theta \mapsto \overline{\PWCal}_{p, k}(\mu_\star, \mu_\theta)$ as $n \rightarrow +\infty$. Then the second statement of Theorem~\ref{Theorem:appendix-final} implies that 
\begin{equation}\label{result:consistency-MEPRW-first}
\limsup_{n \rightarrow +\infty} \argmin_{\theta \in \Theta} \EE[\overline{\PWCal}_{p, k}(\widehat{\mu}_n(\omega), \widehat{\mu}_{\theta, m(n)}) \mid X_{1:n}] \subseteq \argmin_{\theta \in \Theta} \overline{\PWCal}_{p, k}(\mu_\star, \mu_\theta).  
\end{equation}
The next step is to show that, for every $\delta > 0$, there exists a compact set $B \subseteq \Theta$ and $N \in \bn$ such that $\inf_{\theta \in B} \EE[\overline{\PWCal}_{p, k}(\widehat{\mu}_n(\omega), \widehat{\mu}_{\theta, m(n)}) \mid X_{1:n}] \leq \inf_{\theta \in \Theta} \EE[\overline{\PWCal}_{p, k}(\widehat{\mu}_n(\omega), \widehat{\mu}_{\theta, m(n)}) \mid X_{1:n}] + \delta$. In what follows, we prove a stronger statement which states that the above inequality holds true with $\delta = 0$. Indeed, by the same reasoning for the open set case in the proof of epiconvergence, we have
\begin{equation*}
\limsup_{n \rightarrow +\infty} \inf_{\theta \in \Theta} \EE[\overline{\PWCal}_{p, k}(\widehat{\mu}_n(\omega), \widehat{\mu}_{\theta, m(n)}) \mid X_{1:n}] \leq \inf_{\theta \in \Theta} \overline{\PWCal}_{p, k}(\mu_\star, \mu_\theta). 
\end{equation*}
By Assumption~\ref{A3} and using previous argument, $\Theta_\star(\tau)$ is nonempty and compact for some $\tau > 0$. The above inequality implies that there exists $n_1(\omega) > 0$ such that, for all $n \geq n_1(\omega)$, the set $\{\theta \in \Theta: \EE[\overline{\PWCal}_{p, k}(\widehat{\mu}_n(\omega), \widehat{\mu}_{\theta, m(n)}) \mid X_{1:n}] \leq \inf_{\theta' \in \Theta} \overline{\PWCal}_{p, k}(\mu_\star, \mu_{\theta'}) + \tau/3\}$ is nonempty. For any $\theta$ in this set and let $n \geq n_1(\omega)$, we have
\begin{equation*}
\overline{\PWCal}_{p, k}(\mu_\star, \mu_\theta) \leq \overline{\PWCal}_{p, k}(\mu_\star, \widehat{\mu}_n(\omega)) + \EE[\overline{\PWCal}_{p, k}(\mu_\theta, \widehat{\mu}_{\theta, m(n)}) \mid X_{1:n}] + \inf_{\theta \in \Theta} \overline{\PWCal}_{p, k}(\mu_\star, \mu_\theta) + \frac{\tau}{3}. 
\end{equation*}
By Assumption~\ref{A1}, there exists $n_2(\omega) > 0$ such that, for all $n \geq n_2(\omega)$, we have
\begin{equation*}
\overline{\PWCal}_{p, k}(\mu_\star, \widehat{\mu}_n(\omega)) \leq \WCal_p(\mu_\star, \widehat{\mu}_n(\omega)) \leq \frac{\tau}{3}. 
\end{equation*}
By Assumption~\ref{A4}, there exists $n_3(\omega) > 0$ such that, for all $n \geq n_3(\omega)$, we have
\begin{equation*}
\EE[\overline{\PWCal}_{p, k}(\widehat{\mu}_{\theta, m(n)}, \mu_\theta) \mid X_{1:n}] \leq \EE[\WCal_p(\widehat{\mu}_{\theta, m(n)}, \mu_\theta) \mid X_{1:n}] \leq \frac{\tau}{3}. 
\end{equation*}
Putting these pieces together yields that, for all $n \geq \max\{n_1(\omega), n_2(\omega), n_3(\omega)\}$ that,  
\begin{equation*}
\overline{\PWCal}_{p, k}(\mu_\star, \mu_\theta) \leq \inf_{\theta \in \Theta} \overline{\PWCal}_{p, k}(\mu_\star, \mu_\theta)  + \tau. 
\end{equation*}
This implies that, for all $n \geq \max\{n_1(\omega), n_2(\omega), n_3(\omega)\}$ that, 
\begin{equation*}
\left\{\theta \in \Theta: \EE[\overline{\PWCal}_{p, k}(\widehat{\mu}_n(\omega), \widehat{\mu}_{\theta, m(n)}) \mid X_{1:n}] \leq \inf_{\theta' \in \Theta} \overline{\PWCal}_{p, k}(\mu_\star, \mu_{\theta'}) + \frac{\tau}{3}\right\} \subseteq \Theta_\star(\tau). 
\end{equation*}
Therefore, we have $\inf_{\theta \in \Theta} \EE[\overline{\PWCal}_{p, k}(\widehat{\mu}_n(\omega), \widehat{\mu}_{\theta, m(n)}) | X_{1:n}] = \inf_{\theta \in \Theta_\star(\tau)} \EE[\overline{\PWCal}_{p, k}(\widehat{\mu}_n(\omega), \widehat{\mu}_{\theta, m(n)}) | X_{1:n}]$. This together with the compactness of  $\Theta_\star(\tau)$ yields the desired result. 

The first statement of Theorem~\ref{Theorem:appendix-final} implies that
\begin{equation}\label{result:consistency-MEPRW-second}
\inf_{\theta \in \Theta} \EE[\overline{\PWCal}_{p, k}(\widehat{\mu}_n(\omega), \widehat{\mu}_{\theta, m(n)}) \mid X_{1:n}] \rightarrow \inf_{\theta \in \Theta} \overline{\PWCal}_{p, k}(\mu_\star, \mu_\theta), \quad \textnormal{as } n \rightarrow +\infty. 
\end{equation}
By Assumption~\ref{A2} and Lemma~\ref{lemma:EPRW-lsc}, the mapping $\theta \mapsto \EE[\overline{\PWCal}_{p, k}(\widehat{\mu}_n(\omega), \widehat{\mu}_{\theta, m(n)}) \mid X_{1:n}]$ is lower semi-continuous. Theorem~\ref{Theorem:appendix-compact} implies $\argmin_{\theta \in \Theta} \EE[\overline{\PWCal}_{p, k}(\widehat{\mu}_n(\omega), \widehat{\mu}_{\theta, m(n)}) \mid X_{1:n}]$ are nonempty for all $n \geq \max\{n_1(\omega), n_2(\omega), n_3(\omega)\}$. Together with Eq.~\eqref{result:consistency-MEPRW-first} and~\eqref{result:consistency-MEPRW-second} yields the desired results. 

Finally, we remark that these results hold true for $\delta_n$-$\argmin_{\theta \in \Theta} \EE[\overline{\PWCal}_{p, k}(\widehat{\mu}_n(\omega), \widehat{\mu}_{\theta, m(n)}) \mid X_{1:n}]$ with $\delta_n \rightarrow 0$. For Eq.~\eqref{result:consistency-MEPRW-first} and~\eqref{result:consistency-MEPRW-second}, the analogous results can be derived by using the second and third statements of Theorem~\ref{Theorem:appendix-final}. To show that $\delta_n$-$\argmin_{\theta \in \Theta} \EE[\overline{\PWCal}_{p, k}(\widehat{\mu}_n(\omega), \widehat{\mu}_{\theta, m(n)}) \mid X_{1:n}]$ is nonempty, we notice it contains the nonempty set $\argmin_{\theta \in \Theta} \EE[\overline{\PWCal}_{p, k}(\widehat{\mu}_n(\omega), \widehat{\mu}_{\theta, m(n)}) \mid X_{1:n}]$. 

\subsection{Proof of Theorem~\ref{Theorem:convergence-MEPRW-MPRW}}
We first prove that $\argmin_{\theta \in \Theta} \overline{\PWCal}_{p, k}(\widehat{\mu}_n, \mu_\theta) \neq \emptyset$. Indeed, by Assumption~\ref{A2} and Theorem~\ref{Theorem:PRW-lsc}, the mapping $\theta \mapsto \overline{\PWCal}_{p, k}(\widehat{\mu}_n, \mu_\theta)$ is lower semi-continuous. By Assumption~\ref{A5}, the set $\Theta_n(\tau)$ is bounded for some $\tau_n > 0$. By the definition of $\inf$, there exists $\theta'_n \in \Theta$ such that $\overline{\PWCal}_{p, k}(\widehat{\mu}_n, \mu_{\theta'_n}) = \inf_{\theta \in \Theta} \overline{\PWCal}_{p, k}(\widehat{\mu}_n, \mu_\theta) + \tau_n/2$. This implies that $\theta'_n \in \Theta_n(\tau)$ and $\Theta_n(\tau)$ is nonempty. By the lower semi-continuity of the mapping $\theta \mapsto \overline{\PWCal}_{p, k}(\widehat{\mu}_n, \mu_\theta)$, the set $\Theta_n(\tau)$ is closed. Putting these pieces together yields that $\Theta_n(\tau)$ is compact. Therefore, we conclude the desired result from Theorem~\ref{Theorem:appendix-compact}. 

Then we show that the sequences $\theta \mapsto \EE[\overline{\PWCal}_{p, k}(\widehat{\mu}_n, \widehat{\mu}_{\theta, m}) \mid X_{1:n}]$ epiconverges to $\theta \mapsto \overline{\PWCal}_{p, k}(\widehat{\mu}_n, \mu_\theta)$ as $m \rightarrow +\infty$. Indeed, it suffices to verify the conditions in Proposition~\ref{Proposition:appendix-epi}. 

Fix $K \subseteq \Theta$ as an arbitrary compact set. By Assumption~\ref{A2} and Lemma~\ref{lemma:EPRW-lsc}, the mapping $\theta \mapsto \EE[\overline{\PWCal}_{p, k}(\widehat{\mu}_n, \widehat{\mu}_{\theta, m}) \mid X_{1:n}]$ is lower semi-continuous. Then Theorem~\ref{Theorem:appendix-compact} implies that 
\begin{equation*}
\inf_{\theta \in K} \EE[\overline{\PWCal}_{p, k}(\widehat{\mu}_n, \widehat{\mu}_{\theta, m}) \mid X_{1:n}] = \EE[\overline{\PWCal}_{p, k}(\widehat{\mu}_n, \widehat{\mu}_{\theta_m, m}) \mid X_{1:n}]
\end{equation*}
for some sequence $\theta_m \in K$. Thus, we have
\begin{equation*}
\liminf_{m \rightarrow +\infty} \inf_{\theta \in K} \EE\left[\overline{\PWCal}_{p, k}(\widehat{\mu}_n, \widehat{\mu}_{\theta, m}) \mid X_{1:n}\right] = \liminf_{m \rightarrow +\infty} \EE\left[\overline{\PWCal}_{p, k}(\widehat{\mu}_n, \widehat{\mu}_{\theta_m, m}) \mid X_{1:n}\right]. 
\end{equation*}
Following up the same approach used in the proof of Theorem~\ref{Theorem:consistency-MPRW}, there exists a subsequence of $\{\theta_m\}_{m \in \bn}$, denoted by $\{\theta_{m_j}\}_{j \in \bn}$ with the limit $\bar{\theta} \in K$, such that 
\begin{eqnarray*}
\liminf_{m \rightarrow +\infty} \EE[\overline{\PWCal}_{p, k}(\widehat{\mu}_n, \widehat{\mu}_{\theta_m, m}) \mid X_{1:n}] & = & \lim_{j \rightarrow +\infty} \EE[\overline{\PWCal}_{p, k}(\widehat{\mu}_n, \widehat{\mu}_{\theta_{m_j}, m_j}) \mid X_{1:n}] \\
& & \hspace*{-16em} \geq \ \liminf_{j \rightarrow +\infty} \EE[\overline{\PWCal}_{p, k}(\widehat{\mu}_n, \mu_{\theta_{m_j}})] - \limsup_{j \rightarrow +\infty} \EE[\overline{\PWCal}_{p, k}(\mu_{\theta_{m_j}}, \widehat{\mu}_{\theta_{m_j}, m_j}) \mid X_{1:n}]. 
\end{eqnarray*}
Assumption~\ref{A1} and~\ref{A2} imply $\widehat{\mu}_{m_j} \Rightarrow \mu_\star$ and $\mu_{\theta_{m_j}} \Rightarrow \mu_{\bar{\theta}}$. Together with the lower semi-continuity of the PRW distance yields that $\liminf_{j \rightarrow +\infty} \overline{\PWCal}_{p, k}(\widehat{\mu}_n, \mu_{\theta_{m_j}}) \geq \overline{\PWCal}_{p, k}(\widehat{\mu}_n, \mu_{\bar{\theta}})$. By Assumption~\ref{A4} and using $\theta_{m_j} \rightarrow \bar{\theta}$, we have $\limsup_{j \rightarrow +\infty} \EE[\overline{\PWCal}_{p, k}(\mu_{\theta_{m_j}}, \widehat{\mu}_{\theta_{m_j}, m_j}) \mid X_{1:n}] = 0$. Thus, we conclude that $\liminf_{m \rightarrow +\infty} \EE[\overline{\PWCal}_{p, k}(\widehat{\mu}_n, \widehat{\mu}_{\theta_m, m}) \mid X_{1:n}] \geq \inf_{\theta \in K} \overline{\PWCal}_{p, k}(\widehat{\mu}_n, \mu_\theta)$.

Fix $O \subseteq \Theta$ as an arbitary open set. By the definition of $\inf$, there exists a sequence $\theta'_m \in O$ such that $\overline{\PWCal}_{p, k}(\widehat{\mu}_n, \mu_{\theta'_m}) \rightarrow \inf_{\theta \in O} \overline{\PWCal}_{p, k}(\widehat{\mu}_n, \mu_\theta)$. In addition, we have 
\begin{equation*}
\inf_{\theta \in O} \EE[\overline{\PWCal}_{p, k}(\widehat{\mu}_n, \widehat{\mu}_{\theta, m}) \mid X_{1:n}] \leq \EE[\overline{\PWCal}_{p, k}(\widehat{\mu}_n, \widehat{\mu}_{\theta'_m, m}) \mid X_{1:n}]. 
\end{equation*}
Thus, we have
\begin{eqnarray*}
\limsup_{m \rightarrow +\infty} \inf_{\theta \in O} \EE[\overline{\PWCal}_{p, k}(\widehat{\mu}_n, \widehat{\mu}_{\theta, m}) \mid X_{1:n}] & \leq & \limsup_{m \rightarrow +\infty} \EE[\overline{\PWCal}_{p, k}(\widehat{\mu}_n, \widehat{\mu}_{\theta'_m, m}) \mid X_{1:n}] \\
& & \hspace*{-10em} \leq \ \limsup_{m \rightarrow +\infty} \overline{\PWCal}_{p, k}(\widehat{\mu}_n, \mu_{\theta'_m}) + \limsup_{n \rightarrow +\infty} \EE[\overline{\PWCal}_{p, k}(\mu_{\theta'_n}, \widehat{\mu}_{\theta'_m, m}) \mid X_{1:n}]. 
\end{eqnarray*}
By the definition of $\theta'_m$, we have $\limsup_{m \rightarrow +\infty} \overline{\PWCal}_{p, k}(\widehat{\mu}_n, \mu_{\theta'_m})=\inf_{\theta \in O} \overline{\PWCal}_{p, k}(\widehat{\mu}_n, \mu_\theta)$. Using Assumption~\ref{A4} and $\lim_{j \rightarrow +\infty} \theta_{m_j}=\bar{\theta}$, we have $\limsup_{m \rightarrow +\infty} \EE[\overline{\PWCal}_{p, k}(\mu_{\theta'_m}, \widehat{\mu}_{\theta'_m, m}) \mid X_{1:n}] = 0$. Putting these pieces together yields that $\limsup_{m \rightarrow +\infty} \inf_{\theta \in O} \EE[\overline{\PWCal}_{p, k}(\widehat{\mu}_n, \widehat{\mu}_{\theta, m}) \mid X_{1:n}] \leq \inf_{\theta \in O} \overline{\PWCal}_{p, k}(\widehat{\mu}_n, \mu_\theta)$. 

Proposition~\ref{Proposition:appendix-epi} guarantees that the sequence of mappings $\theta \mapsto \EE[\overline{\PWCal}_{p, k}(\widehat{\mu}_n, \widehat{\mu}_{\theta, m}) \mid X_{1:n}]$ epiconverges to the mapping $\theta \mapsto \overline{\PWCal}_{p, k}(\widehat{\mu}_n, \mu_\theta)$ as $m \rightarrow +\infty$. Then the second statement of Theorem~\ref{Theorem:appendix-final} implies that 
\begin{equation}\label{result:convergence-MPRW-MEPRW-first}
\limsup_{m \rightarrow +\infty} \argmin_{\theta \in \Theta} \EE[\overline{\PWCal}_{p, k}(\widehat{\mu}_n, \widehat{\mu}_{\theta, m}) \mid X_{1:n}] \subseteq \argmin_{\theta \in \Theta} \overline{\PWCal}_{p, k}(\widehat{\mu}_n, \mu_\theta).  
\end{equation}
The next step is to show that, for every $\delta > 0$, there exists a compact set $B \subseteq \Theta$ and $N \in \bn$ such that $\inf_{\theta \in B} \EE[\overline{\PWCal}_{p, k}(\widehat{\mu}_n, \widehat{\mu}_{\theta, m}) \mid X_{1:n}] \leq \inf_{\theta \in \Theta} \EE[\overline{\PWCal}_{p, k}(\widehat{\mu}_n, \widehat{\mu}_{\theta, m}) \mid X_{1:n}] + \delta$. In what follows, we prove a stronger statement which states that the above inequality holds true with $\delta = 0$. Indeed, by the same reasoning for the open set case in the proof of epiconvergence, we have
\begin{equation*}
\limsup_{n \rightarrow +\infty} \inf_{\theta \in \Theta} \EE[\overline{\PWCal}_{p, k}(\widehat{\mu}_n, \widehat{\mu}_{\theta, m}) \mid X_{1:n}] \leq \inf_{\theta \in \Theta} \overline{\PWCal}_{p, k}(\widehat{\mu}_n, \mu_\theta). 
\end{equation*}
By Assumption~\ref{A5} and using previous argument, $\Theta_n(\tau)$ is nonempty and compact for some $\tau > 0$. The above inequality implies that there exists $m_1 > 0$ such that, for all $m \geq m_1$, the set $\{\theta \in \Theta: \EE[\overline{\PWCal}_{p, k}(\widehat{\mu}_n, \widehat{\mu}_{\theta, m}) \mid X_{1:n}] \leq \inf_{\theta \in \Theta} \overline{\PWCal}_{p, k}(\widehat{\mu}_n, \mu_\theta) + \tau/2\}$ is nonempty. For any $\theta$ in this set and let $m \geq m_1$, we have
\begin{equation*}
\overline{\PWCal}_{p, k}(\widehat{\mu}_n, \mu_\theta) \leq \EE[\overline{\PWCal}_{p, k}(\widehat{\mu}_{\theta, m}, \mu_\theta) \mid X_{1:n}] + \inf_{\theta \in \Theta} \overline{\PWCal}_{p, k}(\widehat{\mu}_n, \mu_\theta) + \frac{\tau}{2}. 
\end{equation*}
By Assumption~\ref{A4}, there exists $m_2 > 0$ such that, for all $m \geq m_2$, we have
\begin{equation*}
\EE[\overline{\PWCal}_{p, k}(\widehat{\mu}_{\theta, m}, \mu_\theta) \mid X_{1:n}] \leq \EE[\WCal_p(\widehat{\mu}_{\theta, m}, \mu_\theta) \mid X_{1:n}] \leq \frac{\tau}{2}. 
\end{equation*}
Putting these pieces together yields that $\overline{\PWCal}_{p, k}(\widehat{\mu}_n, \mu_\theta) \leq \inf_{\theta' \in \Theta} \overline{\PWCal}_{p, k}(\widehat{\mu}_n, \mu_{\theta'})  + \tau$ for all $m \geq \max\{m_1, m_2\}$. This implies that, for all $m \geq \max\{m_1, m_2\}$ that, 
\begin{equation*}
\left\{\theta \in \Theta: \EE[\overline{\PWCal}_{p, k}(\widehat{\mu}_n, \widehat{\mu}_{\theta, m}) \mid X_{1:n}] \leq \inf_{\theta' \in \Theta} \overline{\PWCal}_{p, k}(\widehat{\mu}_n, \mu_{\theta'}) + \frac{\tau}{2}\right\} \subseteq \Theta_n(\tau). 
\end{equation*}
Therefore, we have $\inf_{\theta \in \Theta} \EE[\overline{\PWCal}_{p, k}(\widehat{\mu}_n, \widehat{\mu}_{\theta, m}) | X_{1:n}] = \inf_{\theta \in \Theta_n(\tau)} \EE[\overline{\PWCal}_{p, k}(\widehat{\mu}_n, \widehat{\mu}_{\theta, m}) | X_{1:n}]$. This together with the compactness of  $\Theta_n(\tau)$ yields the desired result. 

The first statement of Theorem~\ref{Theorem:appendix-final} implies that
\begin{equation}\label{result:convergence-MPRW-MEPRW-second}
\inf_{\theta \in \Theta} \EE[\overline{\PWCal}_{p, k}(\widehat{\mu}_n, \widehat{\mu}_{\theta, m}) \mid X_{1:n}] \rightarrow \inf_{\theta \in \Theta} \overline{\PWCal}_{p, k}(\widehat{\mu}_n, \mu_\theta), \quad \textnormal{as } m \rightarrow +\infty.  
\end{equation}
By Assumption~\ref{A2} and Lemma~\ref{lemma:EPRW-lsc}, the mapping $\theta \mapsto \EE[\overline{\PWCal}_{p, k}(\widehat{\mu}_n, \widehat{\mu}_{\theta, m}) \mid X_{1:n}]$ is lower semi-continuous. Theorem~\ref{Theorem:appendix-compact} implies $\argmin_{\theta \in \Theta} \EE[\overline{\PWCal}_{p, k}(\widehat{\mu}_n, \widehat{\mu}_{\theta, m}) \mid X_{1:n}]$ are nonempty for all $m \geq \max\{m_1, m_2\}$. Together with Eq.~\eqref{result:convergence-MPRW-MEPRW-first} and Eq.~\eqref{result:convergence-MPRW-MEPRW-second} yields the desired results. 

Finally, we remark that these results hold true for $\delta_n$-$\argmin_{\theta \in \Theta} \EE[\overline{\PWCal}_{p, k}(\widehat{\mu}_n, \widehat{\mu}_{\theta, m}) \mid X_{1:n}]$ with $\delta_n \rightarrow 0$. For Eq.~\eqref{result:convergence-MPRW-MEPRW-first} and~\eqref{result:convergence-MPRW-MEPRW-second}, the analogous results can be derived by using the second and third statements of Theorem~\ref{Theorem:appendix-final}. To show that $\delta_n$-$\argmin_{\theta \in \Theta} \EE[\overline{\PWCal}_{p, k}(\widehat{\mu}_n, \widehat{\mu}_{\theta, m}) \mid X_{1:n}]$ is nonempty, we notice it contains the nonempty set $\argmin_{\theta \in \Theta} \EE[\overline{\PWCal}_{p, k}(\widehat{\mu}_n, \widehat{\mu}_{\theta, m}) \mid X_{1:n}]$.  

\subsection{Proof of Lemma~\ref{lemma:EPRW-lsc}}
Since $\nu_i \Rightarrow \nu \in \PScr(\br^d)$ and $\br^d$ is separable, the Skorokhod's representation theorem (cf. Theorem~\ref{Theorem:skorokhod}) implies that there exists $m$ sequences of random variables $\{\{Z_i^k\}_{i \in \bn}, k \in [m]\}$ and $m$ random variables $\{Z^k, k \in [m]\}$ such that the distribution of $Z_i^k$ is $\nu_i$, the distribution of $Z^k$ is $\nu$ and $\{Z_i^k\}_{i \in \bn}$ converges to $Z^k$ almost surely for all $k\in [m]$. 

Suppose that $\widehat{\nu}_{i, m} = (1/m)(\sum_{k=1}^m \delta_{Z_i^k})$ and $\widehat{\nu}_m = (1/m)(\sum_{k=1}^m Z^k)$, we proceed to the key part of the proof and show that $\{\widehat{\nu}_{i, m}\}_{i \in \bn}$ weakly converges to $\widehat{\nu}_m$. Indeed, it suffices to consider the deterministic case where $\widehat{\nu}_{i, m} = (1/m)(\sum_{k=1}^m \delta_{z_i^k})$ and $\widehat{\nu}_m = (1/m)(\sum_{k=1}^m z^k)$ where $\{\{z_i^k\}_{i \in \bn}, k \in [m]\}$ and $\{z^k, k \in [m]\}$ are all deterministic such that $\lim_{i \rightarrow +\infty} \left(\max_{k \in [m]} \|z_i^k - z^k\|\right) = 0$. Since the Wasserstein distance metrizes the weak convergence (cf. Theorem~\ref{appendix:Thm-B4}), we only need to show that $\lim_{i \rightarrow +\infty} \WCal_2(\widehat{\nu}_{i, m}, \widehat{\nu}_m) = 0$. By the definition of the Wasserstein distance, $\{\widehat{\nu}_{i, m}\}_{i \in \bn}$ and $\widehat{\nu}_m$, we have $\WCal_2^2(\widehat{\nu}_{i, m}, \widehat{\nu}_m) \leq \max_{k \in [m]} \|z_i^k - z^k\|^2$. Putting these pieces together yields that $\{\widehat{\nu}_{i, m}\}_{i \in \bn}$ weakly converges to $\widehat{\nu}_m$ almost surely. 

Finally, we conclude from the lower semi-continuity of the PRW distance (cf. Theorem~\ref{Theorem:PRW-lsc}) and the Fatou's lemma that 
\begin{equation*}
\EE[\overline{\PWCal}_{p, k}(\mu, \widehat{\nu}_m)] \leq \EE\left[\liminf_{i \rightarrow +\infty} \overline{\PWCal}_{p, k}(\mu_i, \widehat{\nu}_{i, m})\right] \leq \liminf_{i \rightarrow +\infty} \EE[\overline{\PWCal}_{p, k}(\mu_i, \widehat{\nu}_{i, m})]. 
\end{equation*}
This completes the proof. 

\subsection{Proof of Theorem~\ref{Theorem:measurability-MPRW}}
Using Assumption~\ref{A2} and Theorem~\ref{Theorem:PRW-lsc}, the mapping $(\mu, \theta) \mapsto \overline{\PWCal}_{p, k}(\mu, \mu_\theta)$ is lower semi-continuous in $\PCal(\br^d) \times \Theta$. It remains to verify that the conditions in Theorem~\ref{Theorem:brown-1973} are satisfied. 

We notice that the empirical measure $\widehat{\mu}_n(\omega)$ depends on $\omega \in \Omega$ only through $X_{1:n} \in \otimes_{i=1}^n \br^d$. Thus, we can write $\widehat{\mu}_n(\omega) = \widehat{\mu}_n(x)$ as a function in $\otimes_{i=1}^n \br^d$. Let $D = (\otimes_{i=1}^n \br^d) \times \Theta$, it is a Borel subset of $(\otimes_{i=1}^n \br^d) \times \br$. Since $\br^d$ is a Polish space, $\br^d \times \ldots \times \br^d$ endowed with the product topology is a Polish space. $D_x$ is $\sigma$-compact for any $x \in \proj(D)$ since $D_x \subseteq \Theta$ and $\Theta$ is $\sigma$-compact. 

Define $f(x, \theta) = \overline{\PWCal}_{p, k}(\widehat{\mu}_n(x), \mu_\theta)$, we claim that $f$ is measurable on $D$ and $f(x, \cdot)$ is lower semi-continuous on $D_x$. Indeed, we have shown that the mapping $(\mu, \theta) \mapsto \overline{\PWCal}_{p, k}(\mu, \mu_\theta)$ is lower semi-continuous and thus measurable in $\PCal(\br^d) \times \Theta$. The mapping $x \mapsto \widehat{\mu}_n(x)$ is measurable in $\otimes_{i=1}^n \br^d$. Since the composition of measurable functions is measure, $f$ is measurable on $D$. Moreover, for any $x \in \otimes_{i=1}^n \br^d$, $f(x, \cdot)$ is lower semi-continuous on $D_x$ since the mapping $(\mu, \theta) \mapsto \overline{\PWCal}_{p, k}(\mu, \mu_\theta)$ is lower semi-continuous on $D$. Putting these pieces together yields the desired results. 

\subsection{Proof of Theorem~\ref{Theorem:measurability-MEPRW}}
Using Assumption~\ref{A2} and Lemma~\ref{lemma:EPRW-lsc}, the mapping $(\nu, \theta) \mapsto \EE[\overline{\PWCal}_{p, k}(\nu, \widehat{\mu}_{\theta, m}) \mid X_{1:n}]$ is lower semi-continuous in $\PCal(\br^d) \times \Theta$. Then the proof can be done similarly to the proof of Theorem~\ref{Theorem:measurability-MPRW} using this result and Theorem~\ref{Theorem:brown-1973}. 

\section{Postponed Proofs in Subsection~\ref{subsec:distribution}}
In this section, we provide the detailed proofs for Theorem~\ref{Theorem:general} and Theorem~\ref{Theorem:well-specificied}. Our derivation is the refinement of the analysis in~\citet{Bernton-2019-Parameter} for minimal Wasserstein estimators. 

\subsection{Preliminary technical results}
To facilitate reading, we collect several preliminary technical results which will be used in the postponed proofs in subsection~\ref{subsec:distribution}. 

Let $(\XCal, \|\cdot\|_X)$ be a normed linear space and $\theta \mapsto f_\theta$ be a map from a subset $\Theta$ of $\br^d$ into $\XCal$. The statistical information comes from a sequence $\{f_n\}_{n \in \bn}$ of random elements of $\XCal$, each of which is assumed to be measurable with respect to the $\sigma$-algebra generated by the balls in $\XCal$. In some sense $f_n$ should converge to $f_{\theta_\star}$ where $\theta_\star$ is some fixed (but unknown) point in the interior of $\Theta$. To avoid the abuse of notation, we use $K_1(x, \beta)$ here. 
\begin{theorem}[Theorem~4.2 in~\citet{Pollard-1980-Minimum}]\label{Theorem:Pollard-first}
Suppose the following assumptions hold: 
\begin{enumerate}
\item $\inf_{\theta \notin N} \|f_\theta - f_{\theta_\star}\|_X > 0$ for every neighborhood $N$ of $\theta_\star$. 
\item $\theta \mapsto f_\theta$ is norm differentiable with non-singular derivative $D_{\theta_\star}$ at $\theta_\star$. 
\item There exists a random element $G_\star \in \XCal$ for which $G_n \mydefn \sqrt{n}(f_n - f_{\theta_\star}) \Rightarrow G_\star$ in the sense for the metric induced by the norm $\|\cdot\|_X$. 
\end{enumerate}
Then the limiting distribution of the goodness-of-fit statistic is given by
\begin{equation*}
\sqrt{n}\inf_{\theta \in \Theta} \|f_n - f_\theta\|_X \Rightarrow \inf_{\theta \in \Theta} \|G_\star - \langle \theta, D_{\theta_\star}\rangle\|_X. 
\end{equation*}
\end{theorem}
Let $K_1(x, \beta) = \{\theta: \|x - \langle \theta, D_{\theta_\star}\rangle\|_X \leq \inf_{\theta' \in \Theta} \|x - \langle \theta', D_{\theta_\star}\rangle\|_X + \beta\}$ and $M_n$ is defined by 
\begin{equation*}
M_n = \left\{\theta \in \Theta: \|f_n - f_\theta\|_X \leq \inf_{\theta' \in \Theta} \|f_n - f_{\theta'}\|_X + \eta_n / \sqrt{n}\right\}, 
\end{equation*}
where $\eta_n > 0$ is any sequence such that $\PP(\eta_n \rightarrow 0) = 1$ and $M_n$ is nonempty. 
\begin{theorem}[Theorem~7.2 in~\citet{Pollard-1980-Minimum}]\label{Theorem:Pollard-second}
Under the conditions of Theorem~\ref{Theorem:Pollard-first}, there exists a sequence of real number $\beta_n \downarrow 0$ satisfying
\begin{equation*}
\PP_\star(M_n \subseteq \theta_\star + n^{-1/2}K_1(G_n, \beta_n)) \rightarrow 1, \quad \textnormal{as } n \rightarrow +\infty. 
\end{equation*}
Moreover, for any $\epsilon > 0$, we have $\PP(d_H(K_1(G_n^\star, 0), K_1(G_n, \beta_n)) < \epsilon) \rightarrow 1$ as $n \rightarrow +\infty$. 
\end{theorem}

\subsection{Proof of Theorem~\ref{Theorem:general}}
First, we show that $M_n \subseteq \NCal_1$ with (inner) probability approaching 1 as $n \rightarrow +\infty$. Indeed, with inner probability approaching 1, we have
\begin{equation*}
\argmin_{\theta \in \Theta} \overline{\PWCal}_{1, 1}(\widehat{\mu}_n, \mu_\theta) \subseteq \argmin_{\theta \in \Theta} \overline{\PWCal}_{1, 1}(\mu_\star, \mu_\theta).  
\end{equation*} 
By the definition of $\overline{\PWCal}_{1, 1}$, we conclude that any minimizer of $\|\widehat{F}_n - F_\theta\|_L$ will be included in the set of minimizers of $\|F_\star - F_\theta\|_L$ with inner probability approaching 1. By Assumption~\ref{A9}, the minimizer of $\|F_\star - F_\theta\|_L$ is unique and $\NCal_1$ is the neighborhood of this minimizer. Putting these pieces together yields that the set $\inf_{\theta \in \Theta} \overline{\PWCal}_{1,1}(\widehat{\mu_n}, \mu_\theta)$ is contained in the set $\NCal_1$ with (inner) probability approaching 1 as $n \rightarrow +\infty$. By the definition of $M_n$, we achieve the desired result. 

Then we make three key claims. First, we claim that $M_n \subseteq \Theta_n$ with (inner) probability approaching 1 as $n \rightarrow +\infty$, where $\Theta_n$ is defined by 
\begin{equation*}
\Theta_n = \left\{\theta \in \Theta: \|\theta-\theta_\star\|_\Theta \leq \frac{4\sqrt{n}\|\widehat{F}_n - F_\star\|_L + 2\eta_n}{c_\star\sqrt{n}}\right\}. 
\end{equation*}
Indeed, for any $\theta \in \NCal_1$, we derive from the triangle inequality that 
\begin{equation*}
\|\widehat{F}_n - F_\theta\|_L - \|\widehat{F}_n - F_{\theta_\star}\|_L \geq \|F_\theta - F_\star\|_L - \|F_{\theta_\star} - F_\star\|_L - 2\|\widehat{F}_n - F_\star\|_L. 
\end{equation*}
Using the definition of $\overline{\PWCal}_{1,1}$ together with Assumption~\ref{A9}, we have
\begin{equation}\label{inequality-general-first}
\|\widehat{F}_n - F_\theta\|_L - \|\widehat{F}_n - F_{\theta_\star}\|_L \geq c_\star\|\theta-\theta_\star\|_\Theta - 2\|\widehat{F}_n - F_\star\|_L.
\end{equation}
Since $M_n \subseteq \NCal_1$ with (inner) probability approaching one, Eq.~\eqref{inequality-general-first} holds true for any $\theta \in M_n$ with (inner) probability approaching one. Moreover, by the definition of $M_n$, we have $\theta \in M_n$ satisfies
\begin{equation}\label{inequality-general-second}
\|\widehat{F}_n - F_\theta\|_L \leq \inf_{\theta' \in \Theta} \overline{\PWCal}_{1, 1}(\widehat{\mu}_n, \mu_{\theta'}) + \frac{\eta_n}{\sqrt{n}} \leq \|\widehat{F}_n - F_{\theta_\star}\|_L + \frac{\eta_n}{\sqrt{n}}
\end{equation}
Combining Eq.~\eqref{inequality-general-first}, Eq.~\eqref{inequality-general-second} and the definition of $\Theta_n$, we conclude that $\theta \in \Theta_n$ if $\theta \in M_n$ with (inner) probability approaching 1. This completes the proof the first claim.  

Second, we claim that $\argmin_{\theta' \in \NCal_1} \|G_n - \langle \sqrt{n}(\theta' - \theta_\star), D_{\theta_\star}\rangle\|_L \subseteq \NCal_1 \cap \Theta_n$ with (inner) probability approaching 1 as $n \rightarrow +\infty$. Indeed, by the definition of $G_n$, we have
\begin{equation*}
\|G_n - \langle \sqrt{n}(\theta' - \theta_\star), D_{\theta_\star}\rangle\|_L = \sqrt{n}\|\widehat{F}_n - F_{\theta_\star} - \langle \theta - \theta_\star, D_{\theta_\star}\rangle\|_L. 
\end{equation*}
For the simplicity of notation, we let $R_\theta = F_\theta - F_{\theta_\star} - \langle \theta - \theta_\star, D_{\theta_\star}\rangle$. By Assumption~\ref{A7}, we have $\|R_\theta\|_L = o(\|\theta - \theta_\star\|_\Theta)$. By the definition of $\NCal_1$, we have $\|R_\theta\|_L \leq (1/2)c_\star\|\theta - \theta_\star\|_\Theta$. Therefore, for any $\theta \in \NCal_1$, we have
\begin{eqnarray*}
\|\widehat{F}_n - F_{\theta_\star} - \langle \theta - \theta_\star, D_{\theta_\star}\rangle\|_L & \geq & \|\widehat{F}_n - F_\theta\|_L - \|R_\theta\|_L   \\
& & \hspace*{-10em} \overset{\textnormal{Eq.~\eqref{inequality-general-first}}}{\geq} \ \|\widehat{F}_n - F_{\theta_\star}\|_L + (1/2)c_\star\|\theta-\theta_\star\|_\Theta - 2\|\widehat{F}_n - F_\star\|_L.  
\end{eqnarray*}
This implies that, for any $\theta \in \NCal_1 \setminus \Theta_n$, we have
\begin{equation*}
\|\widehat{F}_n - F_{\theta_\star} - \langle \theta - \theta_\star, D_{\theta_\star}\rangle\|_L \geq \|\widehat{F}_n - F_{\theta_\star}\|_L \geq \inf_{\theta' \in \NCal_1 \cap \Theta_n} \|\widehat{F}_n - F_{\theta_\star} - \langle \theta' - \theta_\star, D_{\theta_\star}\rangle\|_L. 
\end{equation*}
This completes the proof of the second claim. 

Thirdly, we claim that there is an uniform control over the difference between $\theta \mapsto \sqrt{n}\|\widehat{F}_n - F_\theta\|_L$ and the convex map $\theta \mapsto \|G_n - \sqrt{n}\langle \theta - \theta_\star, D_{\theta_\star}\rangle\|_L$ over the set $\Omega_n$ with (inner) probability approaching 1 as $n \rightarrow +\infty$. Indeed, we define
\begin{equation*}
\Gamma_n = \sup_{\theta \in \Omega_n} \ |\sqrt{n}\|\widehat{F}_n - F_\theta\|_L - \|G_n - \sqrt{n}\langle \theta - \theta_\star, D_{\theta_\star}\rangle\|_L|. 
\end{equation*}
By the definition of $G_n$, we have
\begin{eqnarray*}
\Gamma_n & = & \sup_{\theta \in \Omega_n} \ |\sqrt{n}\|\widehat{F}_n - F_{\theta_\star} - \langle\theta - \theta_\star, D_{\theta_\star}\rangle - R_\theta\|_L - \sqrt{n}\|\widehat{F}_n - F_{\theta_\star}- \langle \theta - \theta_\star, D_{\theta_\star}\rangle\|_L| \\ 
& = & o\left(\sup_{\theta \in \Omega_n} \sqrt{n}\|\theta - \theta_\star\|_\Theta\right) \ = \ o(\sqrt{n}\|\widehat{F}_n - F_\star\|_L)
\end{eqnarray*}
By Assumption~\ref{A8}, we have $\Gamma_n \rightarrow 0$ as $\|\theta - \theta_\star\|_\Theta \rightarrow 0$ with (inner) probability approaching 1 as $n \rightarrow +\infty$. This completes the proof of the third claim. 

By the definition of $G_n$ and $G_n^\star$, we have $\|G_n - G_n^\star\|_L = \|\sqrt{n}(\widehat{F}_n - F_\star) - G_\star\|_L$. By Assumption~\ref{A8}, there exists a sequence $\tau_n^1 \rightarrow 0$ such that $\PP(\|G_n - G_n^\star\|_L > \tau_n^1) \rightarrow 0$. By the definition of $\Gamma_n$ and $\eta_n$, there exists two sequences $\tau_n^2 \rightarrow 0$ and $\tau_n^3 \rightarrow 0$ such that $\PP(\Gamma_n > \tau_n^2) \rightarrow 0$ and $\PP(\eta_n > \tau_n^3) \rightarrow 0$. 

Let $\beta_n = \max\{2\tau_n^1, 2\tau_n^2 + \tau_n^3\}$, we have $\beta_n \rightarrow 0$ with (inner) probability approaching 1 as $n \rightarrow +\infty$. It remains to show that $M_n \subseteq K(G_n, \beta_n)$ with (inner) probability approaching 1 as $n \rightarrow +\infty$. Indeed, we have
\begin{equation*}
\inf_{\theta' \in \NCal_1} \|G_n - \langle \sqrt{n}(\theta' - \theta_\star), D_{\theta_\star}\rangle\|_L \geq \inf_{\theta' \in \NCal_1} \sqrt{n}\|\widehat{F}_n - F_{\theta'}\|_L - \tau_n^2. 
\end{equation*}
By the definition of $M_n$, let $\theta \in M_n$, the above inequality implies
\begin{equation*}
\inf_{\theta' \in \NCal_1} \|G_n - \langle \sqrt{n}(\theta' - \theta_\star), D_{\theta_\star}\rangle\|_L \geq \sqrt{n}\|\widehat{F}_n - F_{\theta}\|_L - \tau_n^2 - \tau_n^3. 
\end{equation*}
Since $M_n \subseteq \Theta_n$ with (inner) probability approaching 1 as $n \rightarrow +\infty$, we have
\begin{equation*}
\sqrt{n}\|\widehat{F}_n - F_{\theta}\|_L \geq \|G_n - \langle \sqrt{n}(\theta - \theta_\star), D_{\theta_\star}\rangle\|_L - \tau_n^2. 
\end{equation*}
Putting these pieces together with $\beta_n \geq 2\tau_n^2 + \tau_n^3$ yields that $\theta \in K(G_n, \beta_n)$.

Finally, let $\epsilon > 0$, we prove that $\PP(d_H(K(G_n^\star, 0), K(G_n, \beta_n)) < \epsilon) \rightarrow 1$ as $n \rightarrow +\infty$. Indeed, by the triangle inequality, $\theta \in K(G_n^\star, 0)$ implies $\theta \in K(G_n, 2\|G_n - G_n^\star\|_L)$. Therefore, we conclude that $K(G_n^\star, 0) \subseteq K(G_n, \beta_n)$ with (inner) probability approaching one as $n \rightarrow +\infty$. On the other hand, $\theta \in K(G_n, \beta_n)$ implies $\theta \in K(G_n^\star, \beta_n + 2\|G_n - G_n^\star\|_L)$. By the definition of $\beta_n$, $G_n$ and $G_n^\star$, we obtain that $\beta_n + 2\|G_n - G_n^\star\|_L \rightarrow 0$ with (inner) probability approaching one as $n \rightarrow +\infty$. By the definition of the Hausdorff metric, we conclude the desired result. 

\subsection{Proof of Theorem~\ref{Theorem:well-specificied}}
Different from Theorem~\ref{Theorem:general}, the proof of Theorem~\ref{Theorem:well-specificied} is relatively straightforward and based on Theorem~\ref{Theorem:Pollard-first} and~\ref{Theorem:Pollard-second}. It is mostly because there exists $\theta_\star$ in the interior of $\Theta$ such that $F_\star = F_{\theta_\star}$. 

More specifically, we consider $f_\theta = F_\theta$ and $f_n = \widehat{F}_n$ such that 
\begin{equation*}
F_\theta(u, t) = \int_{\br^d} \one_{(-\infty, t]}(\langle u, x\rangle) \; d\mu_\theta(x), \quad \widehat{F}_n(u, t) = (1/n)|\{i \in [n]: \langle u, X_i\rangle \leq t\}|. 
\end{equation*}
Let $\XCal = L(\bs^{d-1} \times \br)$ and $\|\cdot\|_X = \|\cdot\|_L$, we can check that $(\XCal, \|\cdot\|_X)$ is a normed linear space. By the definition of $\overline{\PWCal}_{1,1}$, we have $\overline{\PWCal}_{1,1}(\widehat{\mu}_n, \mu_\theta) = \|\widehat{F}_n - F_\theta\|_X$. By Assumption~\ref{A1}, $\widehat{F}_n$ converges to $F_\star$. Moreover, in well-specified setting, $F_\star = F_{\theta_\star}$ where $\theta_\star$ is some fixed (but unknown) point in the interior of $\Theta$. Now we are ready to check the conditions of Theorem~\ref{Theorem:Pollard-first}. 

First, Assumption~\ref{A6} and $\overline{\PWCal}_{1,1}(\widehat{\mu}_n, \mu_\theta) = \|\widehat{F}_n - F_\theta\|_X$ imply C1. Furthermore, by the definition of norm differentiable, Assumption~\ref{A7} and Assumption~\ref{A10} imply C2. Finally, Assumption~\ref{A8} and $F_\star = F_{\theta_\star}$ imply C3. Therefore, we conclude from Theorem~\ref{Theorem:Pollard-first} that    
\begin{equation*}
\sqrt{n}\inf_{\theta \in \Theta} \overline{\PWCal}_{1,1}(\widehat{\mu}_n, \mu_\theta) = \sqrt{n}\inf_{\theta \in \Theta} \|\widehat{F}_n - F_\theta\|_L \Rightarrow \inf_{t \in \Theta} \|G_\star - \langle t, D_{\theta_\star}\rangle\|_L. 
\end{equation*}
in the sense for the metric induced by the norm $\|\cdot\|_L$. This together with the definition of the norm $\|\cdot\|_L$ implies the desired result for the goodness-of-fit statistics. 

On the other hand, Theorem~\ref{Theorem:Pollard-second} can be applied with specific choice of $\eta_n$. More specifically, we notice that the estimator $\widehat{\theta}_n$ is well defined by 
\begin{equation*}
\widehat{\theta}_n \mydefn \argmin_{\theta \in \Theta} \overline{\PWCal}_{1, 1}(\widehat{\mu}_n, \mu_\theta) = \argmin_{\theta \in \Theta} \|\widehat{F}_n - F_\theta\|_L. 
\end{equation*}
Let $\eta_n = 0$, the set $M_n = \{\widehat{\theta}_n\}$ is a singleton set. This implies that $\sqrt{n}(\widehat{\theta}_n - \theta_\star) \Rightarrow K_1(G_\star, 0)$ as $n \rightarrow +\infty$ under its Hausdorff metric topology. Since the random map $\theta \rightarrow \max_{u \in \bs^{d-1}} \int_\br |G_\star(u, t) - \langle\theta, D_\star(u, t)\rangle| \; dt$ has a unique infimum almost surely, we have $K_1(G_\star, 0)$ is a singleton set defined by 
\begin{equation*}
K_1(G_\star, 0) = \argmin_{\theta \in \Theta}\max_{u \in \bs^{d-1}} \int_\br |G_\star(u, t) - \langle\theta, D_\star(u, t)\rangle| \; dt. 
\end{equation*}
In this case, the Hausdorff metric is simply induced by the norm $\|\cdot\|_L$. Putting these pieces together yields the desired result for the MPRW estimator of order 1. 

\subsection{Minor Technical Issues}\label{appendix:minor}
We use the notations of~\citet[Theorem~B.8]{Bernton-2019-Parameter} throughout this subsection. Indeed, in page 38-39 of the recent arvix version of~\citet{Bernton-2019-Parameter}, the authors prove that $m(H_n) = \inf_{u \in L_n} f(H_n, u)$, implicitly assuming that the minimizer of the map $\theta \mapsto \sqrt{n}\|F_n - F_{\theta_\star} - \langle\theta - \theta_\star, D_{\theta_\star}\rangle\|_{L_1}$ is contained in the set $\NCal_1 = \{\theta \in \NCal: \|\theta - \theta_\star\|_{\mathcal{H}} \leq c_\star/2\}$. However, this result is not obvious. Indeed, it seems difficult to derive such results from the existing fact that the minimizer of $\theta \mapsto \sqrt{n}\|F_n - F_\theta\|_{L_1}$ is contained in $\NCal$. We \textbf{only} have the uniform control over the difference between $\theta \mapsto \sqrt{n}\|F_n - F_\theta\|_{L_1}$ and $\theta \mapsto \sqrt{n}\|F_n - F_{\theta_\star} - \langle\theta - \theta_\star, D_{\theta_\star}\rangle\|_{L_1}$ over the set $S_n$ instead of the whole set. So there is few relationship between the minimizers of these two mappings. Moreover, the techniques from the proof of~\citet[Theorem~7.2]{Pollard-1980-Minimum} can not be applicable to fix this issue here since the proof depends on the assumption that $\mu_\star = \mu_{\theta_\star}$ which does not hold under model misspecification yet.

\section{Computational Aspects}\label{appendix:computational}

The computation of the PRW distance is in general computationally intractable when the projection dimension is $k \geq 2$ since this amounts to solving a nonconvex max-min optimization model. Despite several pessimistic results~\citep{Paty-2019-Subspace, Weed-2019-Estimation}, we adopt the Riemannian optimization toolbox~\citep{Absil-2009-Optimization} to develop a Riemannian supergradient algorithm and empirically show that our algorithm can approximate $\overline{\PWCal}_{2, k}(\widehat{\mu}_n, \widehat{\nu}_n)$ when the projection dimension is $k \geq 2$. Part of results can be found in the appendix of concurrent work~\citep{Lin-2020-Projection} and we provide the details for the sake of completeness.  

\paragraph{Approximation of $\overline{\PWCal}_{2, k}$.} We consider the computation of $\overline{\PWCal}_{2, k}$ between empirical measures. Indeed, let $\{x_1, x_2, \ldots, x_n\} \subseteq \br^d$ and $\{y_1, y_2, \ldots, y_n\} \subseteq \br^d$ denote sets of $n$ atoms, and let $(r_1, r_2, \ldots, r_n) \in \Delta^n$ and $(c_1, c_2, \ldots, c_n) \in \Delta^n$ denote weight vectors, we define discrete measures $\widehat{\mu}_n \mydefn \sum_{i=1}^n r_i\delta_{x_i}$ and $\widehat{\nu}_n \mydefn \sum_{j=1}^n c_j\delta_{y_j}$. The computation of $\overline{\PWCal}_{2, k}(\widehat{\mu}_n, \widehat{\nu}_n)$ is equivalent to solving a structured max-min optimization model where the maximization and minimization are performed over the Stiefel manifold $\St(d, k) : = \{U \in \br^{d \times k} \mid U^\top U = I_k\}$ and the transportation polytope $\Pi(\mu, \nu) : = \{\pi \in \br_+^{n \times n} \mid r(\pi) = r, \ c(\pi) = c\}$ respectively. Formally, we have
\begin{equation}\label{prob:main}
\max\limits_{U \in \br^{d \times k}} \min\limits_{\pi \in \br_+^{n \times n}} \sum_{i=1}^n \sum_{j=1}^n \pi_{i, j}\|U^\top x_i - U^\top y_j\|^2 \quad \st \ U^\top U = I_k, \ r(\pi) = r, \ c(\pi) = c.
\end{equation}
Eq.~\eqref{prob:main} is equivalent to the non-convex nonsmooth optimization model as follows, 
\begin{equation}\label{prob:Stiefel-nonsmooth}
\max\limits_{U \in \St(d, k)} \ \left\{f(U) \mydefn \min\limits_{\pi \in \Pi(\mu, \nu)} \sum_{i=1}^n \sum_{j=1}^n \pi_{i, j} \|U^\top x_i - U^\top y_j\|^2\right\}.  
\end{equation}
Fixing $U \in \St(d, k)$, Eq.~\eqref{prob:Stiefel-nonsmooth} becomes a classical OT problem which can be either solved by the Sinkhorn iteration~\citep{Cuturi-2013-Sinkhorn} or the variant of network simplex method in the \textsc{POT} package~\citep{Flamary-2017-Pot}. The key challenge is the maximization over the Stiefel manifold $\St(d, k) : = \{U \in \br^{d \times k} \mid U^\top U = I_k\}$. 

Eq.~\eqref{prob:Stiefel-nonsmooth} is a special instance of the Stiefel manifold optimization problem. The dimension of $\St(d, k)$ is equal to $dk - k(k+1)/2$ and the tangent space at the point $Z \in \St(d, k)$ is defined by $\Tg_Z\St \mydefn \{\xi \in \br^{d \times k}: \xi^\top Z + Z^\top\xi = 0\}$. We endow $\St(d, k)$ with Riemannian metric inherited from the Euclidean inner product $\langle X, Y\rangle$ for any $X, Y \in \Tg_Z\St$ and $Z \in \St(d, k)$. Then the projection of $G \in \br^{d \times k}$ onto $\Tg_Z\St$ is given by~\citet[Example~3.6.2]{Absil-2009-Optimization}: $P_{\Tg_Z\St}(G) = G - Z(G^\top Z + Z^\top G)/2$. We make use of the notion of a \textit{retraction}, which is the first-order approximation of an exponential mapping on the manifold and which is amenable to computation~\citep[Definition~4.1.1]{Absil-2009-Optimization}. For the Stiefel manifold, we have the following definition: 
\begin{definition}\label{def:retraction}
A retraction on $\St \equiv \St(d, k)$ is a smooth mapping $\retr: \Tg\St \rightarrow \St$ from the tangent bundle $\Tg\St$ onto $\St$ such that the restriction of $\retr$ onto $\Tg_Z\St$, denoted by $\retr_Z$, satisfies that (i) $\retr_Z(0) = Z$ for all $Z \in \St$ where $0$ denotes the zero element of $\textnormal{T}\St$, and (ii) for any $Z \in \St$, it holds that $\lim_{\xi \in \Tg_Z\St, \xi \rightarrow 0} \|\retr_Z(\xi) - (Z + \xi)\|_F/\|\xi\|_F = 0$.
\end{definition}
Our algorithm uses the retraction based on the \textbf{QR decomposition} as suggested by~\citet{Liu-2019-Quadratic}. More specifically, $\retr_Z^{\textnormal{qr}}(\xi) = \textnormal{qr}(Z + \xi)$ where $\textnormal{qr}(A)$ is the Q factor of the QR factorization of $A$. 

We start with a brief overview of the Riemannian supergradient ascent algorithm for nonsmooth Stiefel optimization, denoted by $\max_{U \in \St(d, k)} F(U)$. A generic Riemannian supergradient ascent algorithm for solving this problem is given by
\begin{equation*}
U_{t+1} \ \leftarrow \ \retr_{U_t}(\gamma_{t+1}\xi_{t+1}) \quad \textnormal{ for any } \xi_{t+1} \in \subdiff F(U_t),   
\end{equation*} 
where $\subdiff F(U_t)$ is Riemannian subdifferential of $F$ at $U_t$ and $\retr$ is any retraction on $\St(d, k)$. The step size is set as $\gamma_{t+1} = \gamma_0/\sqrt{t+1}$ as suggested by~\citep{Li-2019-Nonsmooth}. By the definition of Riemannian subdifferential, $\xi_t$ can be obtained by taking $\xi \in \partial F(U)$ and by setting $\xi_t = P_{\Tg_U\St}(\xi)$. Thus, it is necessary for us to specify the subdifferential of $f$ in Eq.~\eqref{prob:Stiefel-nonsmooth}. We define $V_\pi = \sum_{i=1}^n \sum_{j=1}^n \pi_{i, j} (x_i - y_j)(x_i - y_j)^\top \in \br^{d \times d}$ which is symmetry and derive that 
\begin{equation*}
\partial f(U) \ = \ \textnormal{Conv}\{2V_{\pi^\star} U \mid \pi^\star \in \argmin\limits_{\pi \in \Pi(\mu, \nu)} \ \langle UU^\top, V_\pi\rangle\}, \quad \textnormal{ for any } U \in \br^{d \times k},   
\end{equation*}
It remains to solve an OT problem with a given $U$ at each inner loop of the maximization and use the output $\pi(U)$ to obtain a supergradient of $f$. The network simplex method can exactly solve this LP. To this end, we summarize the pseudocode of the RSGAN algorithm in Algorithm~\ref{alg:supergrad-simplex}.
\begin{algorithm}[!t]
\caption{Riemannian SuperGradient Ascent with Network Simplex Iteration (RSGAN)}\label{alg:supergrad-simplex}
\begin{algorithmic}[1]
\STATE \textbf{Input:} measures $\{(x_i, r_i)\}_{i \in [n]}$ and $\{(y_j, c_j)\}_{j \in [n]}$, dimension $k$ and tolerance $\epsilon$.  
\STATE \textbf{Initialize:} $U_0 \in \St(d, k)$ and $\gamma_0 > 0$. 
\FOR{$t = 0, 1, 2, \ldots, T-1$}
\STATE Compute $\pi_{t+1} \leftarrow \textsc{OT}(\{(x_i, r_i)\}_{i \in [n]}, \{(y_j, c_j)\}_{j \in [n]}, U_t)$.
\STATE Compute $\xi_{t+1} \leftarrow P_{\Tg_{U_t}\St}(2V_{\pi_{t+1}}U_t)$. 
\STATE Compute $\gamma_{t+1} \leftarrow \gamma_0/\sqrt{t+1}$.
\STATE Compute $U_{t+1} \leftarrow \retr_{U_t}(\gamma_{t+1}\xi_{t+1})$.  
\ENDFOR
\end{algorithmic}
\end{algorithm}
\paragraph{Approximation of $\underline{\PWCal}_{2, k}$.} We recall the definition of the IPRW distance of order 2 as follows, 
\begin{equation*}
\underline{\PWCal}_{2, k}^2(\mu, \nu) = \int_{\bs_{d, k}} \WCal_2^2(E_{\#}^\star\mu, E_{\#}^\star\nu) d\sigma(E),   
\end{equation*}
where $\sigma$ is the uniform distribution on $\bs_{d, k}$ and $E^\star$ is the linear transformation associated with $E$ for any $x \in \br^d$ by $E^\star(x) = E^\top x$. For any measurable function $f$ and $\mu \in \PScr(\br^d)$, we denote $f_{\#}\mu$ as the push-forward of $\mu$ by $f$, so that $f_{\#}\mu(A) = \mu(f^{-1}(A))$ where $f^{-1}(A) = \{x \in \br^d: f(x) \in A\}$ for any Borel set $A$. We approximate the integral by selecting a finite set of projections $\SCal \subseteq \bs_{d, k}$ and computing the empirical average: 
\begin{equation*}
\underline{\PWCal}_{2, k}^2(\mu, \nu) \approx \frac{1}{\card(\SCal)}\sum_{E \in \SCal} \WCal_2^2(E_{\#}^\star\mu, E_{\#}^\star\nu).  
\end{equation*}
The quality of this approximation depends on the sampling of $\bs_{d, k}$. In this paper, we use random projections picked uniformly on $\bs_{d, k}$, which is analogues to the approach proposed by~\citet{Bonneel-2015-Sliced} for the case of $k=1$; see \textbf{Sampling schemes} for the details. 

\paragraph{Approximation of $\overline{\PWCal}_{p, 1}$.} We recall the definition of the PRW distance of order $p$ with the projection dimension $k = 1$ as follows, 
\begin{equation*}
\overline{\PWCal}^p_{p, 1}(\mu, \nu) \mydefn \sup_{u \in \bs_{d, 1}} \WCal_p^p(u_{\#}^\star\mu, u_{\#}^\star\nu) = \sup_{u \in \bs_{d, 1}} \int_0^1 |F_{u_{\#}^\star\mu}^{-1}(t) - F_{u_{\#}^\star\nu}^{-1}(t)|^p dt.  
\end{equation*}
where $u \in \bs_{d, 1}$ is an unit $d$-dimensional vector, $u^\star$ is the linear transformation associated with $u$ for any $x \in \br^d$ by $u^\star(x) = u^\top x$, and $F_\xi^{-1}$ is the quantile function of $\xi$. This integral can be estimated using a Monte Carlo estimate and a linear interpolation of the quantile function. Following up~\citet[Appendix~4]{Nadjahi-2019-Asymptotic}, we consider two approximations of this quantity. The first one is given by, 
\begin{equation}\label{def:approx-I}
\overline{\PWCal}^p_{p, 1}(\mu, \nu) = \sup_{u \in \bs_{d, 1}} \frac{1}{K}\sum_{k=1}^K |\tilde{F}_{u_{\#}^\star\mu}^{-1}(t_k) - \tilde{F}_{u_{\#}^\star\nu}^{-1}(t_k)|^p,   
\end{equation} 
where $\{t_k\}_{k=1}^K$ are uniform and independent samples from $[0, 1]$ and $\tilde{F}_\xi^{-1}$ is a linear interpolation of $F_\xi^{-1}$ which denotes either the exact quantile function of a discrete measure $\xi$, or an approximation by a Monte Carlo procedure. The second one is given by 
\begin{equation}\label{def:approx-II}
\overline{\PWCal}^p_{p, 1}(\mu, \nu) = \sup_{u \in \bs_{d, 1}} \frac{1}{K}\sum_{k=1}^K |s_k - \tilde{F}_{u_{\#}^\star\nu}^{-1}(\tilde{F}_{u_{\#}^\star\mu}(s_k))|^p,   
\end{equation} 
where $\{s_k\}_{k=1}^K$ are uniform and independent samples from $u_{\#}^\star\mu$ and $\tilde{F}_\xi$ (resp. $\tilde{F}_\xi^{-1}$) is a linear interpolation of $F_\xi$ (resp. $F_\xi^{-1}$) which denotes either the exact cumulative distribution function (resp. quantile function) of a discrete measure $\xi$, or an approximation by a Monte Carlo procedure.

\paragraph{Sampling schemes.} We explain the methods that we use to generate the i.i.d. samples from the uniform distribution on the set of $d \times k$ orthogonal matrices, i.e., $\bs_{d, k} = \{E \in \br^{d \times k}: E^\top E = I_k\}$ and the i.i.d. samples from multivariate elliptically contoured stable distributions. 

To sample from $\bs_{d, k}$, we first construct the $(d \times k)$-dimensional matrix $Z$ by drawing each of its components from the standard normal distribution $\NCal(0, 1)$ and then perform the QR decomposition of it: $E = \textnormal{qr}(Z)$. By the definition, $E \in \bs_{d, k}$ is an uniform sample.    

To sample from multivariate elliptically contoured stable distributions, we follows the approach presented in~\citet[Appendix~4]{Nadjahi-2019-Asymptotic}. Indeed, we recall that if $Y \in \br^d$ is $\alpha$-stable and elliptically contoured, i.e., $Y \in \ECal\alpha\SCal_c(\Sigma, \textbf{m})$, then its joint characteristic function is defined as, for any $t \in \br^d$ that,  
\begin{equation}\label{def:EASC}
\EE\left[\exp(it^\top Y)\right] \ = \ \exp\left(-(t^\top\Sigma t)^{\alpha/2} + it^\top\textbf{m}\right),  
\end{equation}
where $\Sigma$ is a positive definite matrix (akin to a correlation matrix), $\textbf{m} \in \br^d$ is a location vector (equal to the mean if it exists) and $\alpha \in (0, 2)$ controls the thickness of the tail. Elliptically contoured stable distributions are scale mixtures of multivariate Gaussian distributions~\citep[Proposition~2.5.2]{Samoradnitsky-2017-Stable} with computationally intractable densities. Fortunately, it was shown by~\citet{Nolan-2013-Multivariate} that sampling from multivariate elliptically contoured stable distributions is possible: let $A \sim \SCal_{\alpha/2}(\beta, \gamma, \delta)$ be a one-dimensional positive $(\alpha/2)$-stable random variable with $\beta=1$, $\gamma=2\cos(\pi\alpha/4)^{2/\alpha}$ and $\delta=0$, and $G \sim \NCal(0, \Sigma)$. By the definition, $Y = \sqrt{A}G + \textbf{m}$ satisfies Eq.~\eqref{def:EASC} and $Y \sim \ECal\alpha\SCal_c(\Sigma, \textbf{m})$. 

\paragraph{Optimization methods.} Computing the MPRW and MEPRW estimators are intractable in general. This is mainly because the PRW distance requires a maximization over infinitely many projections. Formally, we hope to solve the following minimax optimization model,  
\begin{equation*}
\min_{\theta \in \Theta} \overline{\PWCal}_{p, 1}^p(\mu_\theta, \mu_\star) \ = \ \min_{\theta \in \Theta} \max_{u \in \bs_{d, 1}} \int_0^1 |F_{u_{\#}^\star\mu_\theta}^{-1}(t) - F_{u_{\#}^\star\mu_\star}^{-1}(t)|^p dt,  
\end{equation*}
where $\{\mu_\theta: \theta \in \Theta\}$ is the model and $\mu_\star$ is the data-generating process. Following up the approach presented in~\citet{Nadjahi-2019-Asymptotic} together with the approximation of $\overline{\PWCal}_{p, 1}$, we consider using the ADAM optimization method to minimize the (expected) PRW distance over the set of parameters while applying multiple projected supergradient ascent to find an approximate projection $u$ which maximizes over $\bs_{d, 1}$ at each inner loop. The ADAM optimization method is associated with the default parameter setting as suggested by~\citet{Kingma-2015-ADAM}. At each inner loop, we run 5 projected supergradient ascent with the learning rate $10^{-3}$. 

\textit{Gaussian models}. For the MPRW estimator, we consider the approximate $\overline{\PWCal}_{2, 1}^2$ distance based on Eq.~\eqref{def:approx-II}. Indeed, let $\mu$ denote $\NCal(\textbf{m}, \sigma^2\textbf{I})$ and $\widehat{\nu}$ denote the empirical probability measures of $n$ samples drawn from the data-generating process, we define the function $f_1(\mathbf{m}, \sigma^2, u)$ as 
\begin{equation*}
f_1(\mathbf{m}, \sigma^2, u) \ = \ \frac{1}{\card(\SCal)}\sum_{s \in \SCal} |s - \tilde{F}_{u_{\#}^\star\widehat{\nu}}^{-1}(\tilde{F}_{u_{\#}^\star\mu}(s))|^2\NCal(s; u^\top\textbf{m}, \sigma^2\textbf{I}), 
\end{equation*}
where $\SCal \subseteq \br$ and $\NCal(s; u^\top\textbf{m}, \sigma^2\textbf{I})$ refers to the density function of Gaussian of parameters $(u^\top\textbf{m}, \sigma^2\textbf{I})$ evaluated at $s \in \SCal$. We compute the explicit gradient expression of $f_1(\mathbf{m}, \sigma^2, u)$ with respect to the mean $\mathbf{m}$, the variance $\sigma^2$ and the projection vector $u$ as follows, 
\begin{eqnarray*}
\nabla_{\textbf{m}} f_1(\mathbf{m}, \sigma^2, u) & = & \frac{1}{\sigma^2\card(\SCal)}\sum_{s \in \SCal} \left(|s - \tilde{F}_{u_{\#}^\star\widehat{\nu}}^{-1}(\tilde{F}_{u_{\#}^\star\mu}(s))|^2\NCal(s; u^\top\textbf{m}, \sigma^2\textbf{I})(s - u^\top\textbf{m})u\right), \\
\nabla_{\sigma^2} f_1(\mathbf{m}, \sigma^2, u) & = & \frac{1}{2\sigma^4\card(\SCal)}\sum_{s \in \SCal} \left(|s - \tilde{F}_{u_{\#}^\star\widehat{\nu}}^{-1}(\tilde{F}_{u_{\#}^\star\mu}(s))|^2\NCal(s; u^\top\textbf{m}, \sigma^2\textbf{I})((s - u^\top\textbf{m})^2 - \sigma^2)\right), \\
\nabla_u f_1(\mathbf{m}, \sigma^2, u) & = & \frac{1}{\sigma^2\card(\SCal)}\sum_{s \in \SCal} \left(|s - \tilde{F}_{u_{\#}^\star\widehat{\nu}}^{-1}(\tilde{F}_{u_{\#}^\star\mu}(s))|^2\NCal(s; u^\top\textbf{m}, \sigma^2\textbf{I})(s - u^\top\textbf{m})\textbf{m}\right). 
\end{eqnarray*}
For the MEPRW estimator, we consider the approximate $\overline{\PWCal}_{2, 1}^2$ distance based on Eq.~\eqref{def:approx-I}. Indeed, let $\widehat{\mu}$ and $\widehat{\nu}$ denote the empirical probability measures of $m$ samples drawn from $\NCal(\textbf{m}, \sigma^2\textbf{I})$ and $n$ samples drawn from the data-generating process, we define the function $f_2(\mathbf{m}, \sigma^2, u)$ as 
\begin{equation*}
f_2(\mathbf{m}, \sigma^2, u) \ = \ \frac{1}{K}\sum_{k=1}^K |\tilde{F}_{u_{\#}^\star\widehat{\mu}}^{-1}(t_k) - \tilde{F}_{u_{\#}^\star\widehat{\nu}}^{-1}(t_k)|^2,  
\end{equation*}
where $\{t_k\}_{k=1}^K$ are uniform and independent samples from $[0, 1]$. We compute the explicit gradient expression of $f_2(\mathbf{m}, \sigma^2, u)$ with respect to the mean $\mathbf{m}$, the variance $\sigma^2$ and the projection vector $u$ as follows, 
\begin{eqnarray*}
\nabla_{\textbf{m}} f_2(\mathbf{m}, \sigma^2, u) & = & -\frac{2}{K}\sum_{k=1}^K |\tilde{F}_{u_{\#}^\star\widehat{\mu}}^{-1}(t_k) - \tilde{F}_{u_{\#}^\star\widehat{\nu}}^{-1}(t_k)|u, \\
\nabla_{\sigma^2} f_2(\mathbf{m}, \sigma^2, u) & = & -\frac{2}{K}\sum_{k=1}^K |\tilde{F}_{u_{\#}^\star\widehat{\mu}}^{-1}(t_k) - \tilde{F}_{u_{\#}^\star\widehat{\nu}}^{-1}(t_k)|\textbf{m}, \\
\nabla_u f_2(\mathbf{m}, \sigma^2, u) & = & -\frac{1}{\sigma^2 K}\sum_{k=1}^K \left(|\tilde{F}_{u_{\#}^\star\widehat{\mu}}^{-1}(t_k) - \tilde{F}_{u_{\#}^\star\widehat{\nu}}^{-1}(t_k)|(u^\top\textbf{m} - \tilde{F}_{u_{\#}^\star\widehat{\mu}}^{-1}(t_k))\right). 
\end{eqnarray*}
\textit{Elliptically contoured stable models.} When comparing the MEPRW estimator with the MPRW estimator using elliptically contoured stable models, we also approximate these estimators using the ADAM optimization method with the default parameter setting. 

We consider the approximate $\overline{\PWCal}_{2, 1}^2$ distance based on Eq.~\eqref{def:approx-I}. Indeed, let $\widehat{\mu}$ and $\widehat{\nu}$ denote the empirical probability measures of $m$ samples drawn from $\ECal\alpha\SCal_c(\textbf{I}, \textbf{m})$ and $n$ samples drawn from the data-generating process, we define the function $f_3(\mathbf{m}, u)$ as 
\begin{equation*}
f_3(\mathbf{m}, u) \ = \ \frac{1}{K}\sum_{k=1}^K |\tilde{F}_{u_{\#}^\star\widehat{\mu}}^{-1}(t_k) - \tilde{F}_{u_{\#}^\star\widehat{\nu}}^{-1}(t_k)|^2. 
\end{equation*}
where $\{t_k\}_{k=1}^K$ are uniform and independent samples from $[0, 1]$. We compute the explicit gradient expression of $f(\mathbf{m}, u)$ with respect to the location parameter $\mathbf{m}$ and the projection vector $u$ as follows, 
\begin{eqnarray*}
\nabla_{\textbf{m}} f_3(\mathbf{m}, u) & = & -\frac{2}{K}\sum_{k=1}^K |\tilde{F}_{u_{\#}^\star\widehat{\mu}}^{-1}(t_k) - \tilde{F}_{u_{\#}^\star\widehat{\nu}}^{-1}(t_k)|u, \\ 
\nabla_u f_3(\mathbf{m}, u) & = & -\frac{2}{K}\sum_{k=1}^K |\tilde{F}_{u_{\#}^\star\widehat{\mu}}^{-1}(t_k) - \tilde{F}_{u_{\#}^\star\widehat{\nu}}^{-1}(t_k)|\mathbf{m}. 
\end{eqnarray*}
\textit{Generative modeling.} We use the ADAM optimizer provided Pytorch GPU.
\begin{figure}[!t]
\centering
\includegraphics[width=0.32\textwidth]{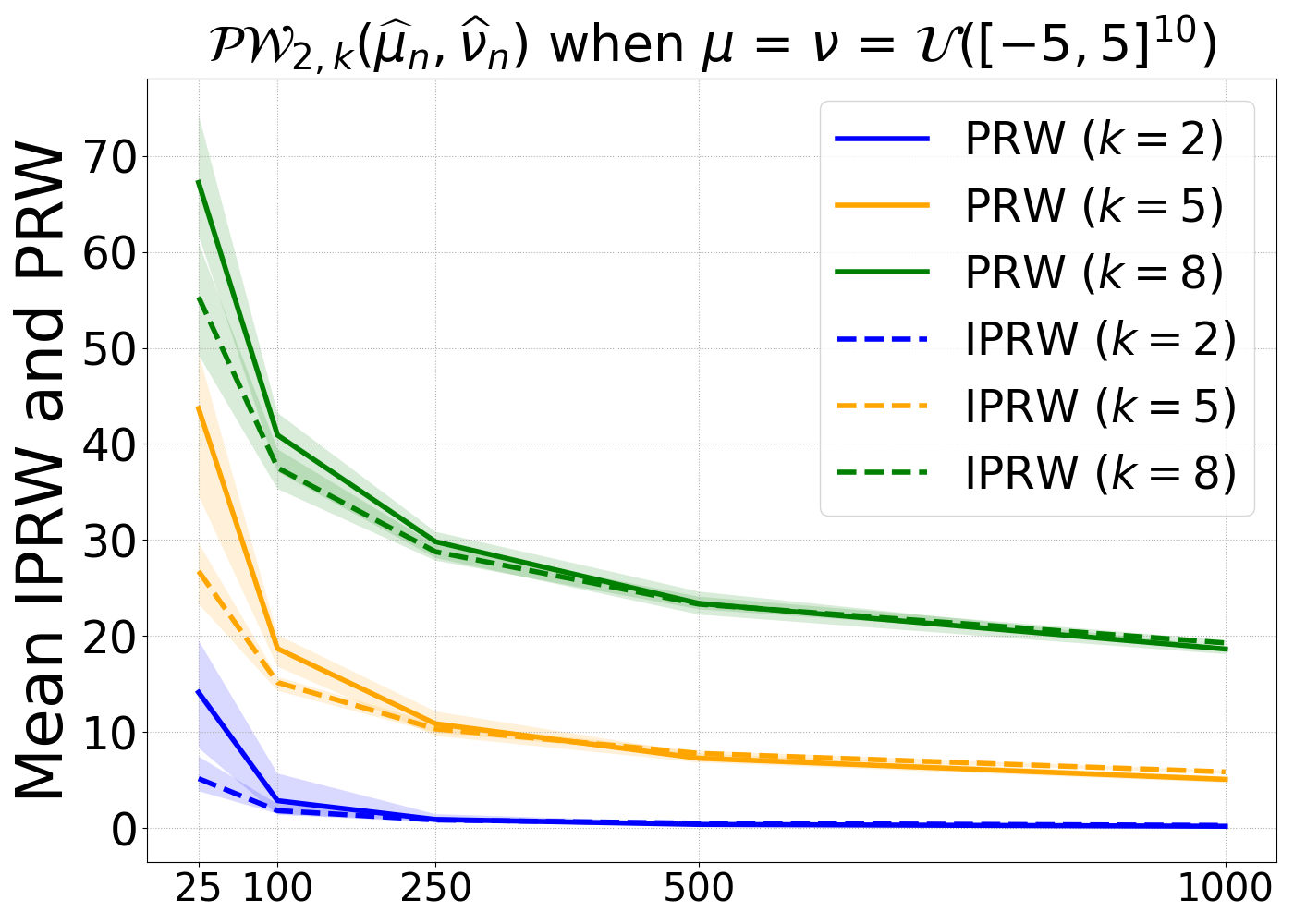}
\includegraphics[width=0.32\textwidth]{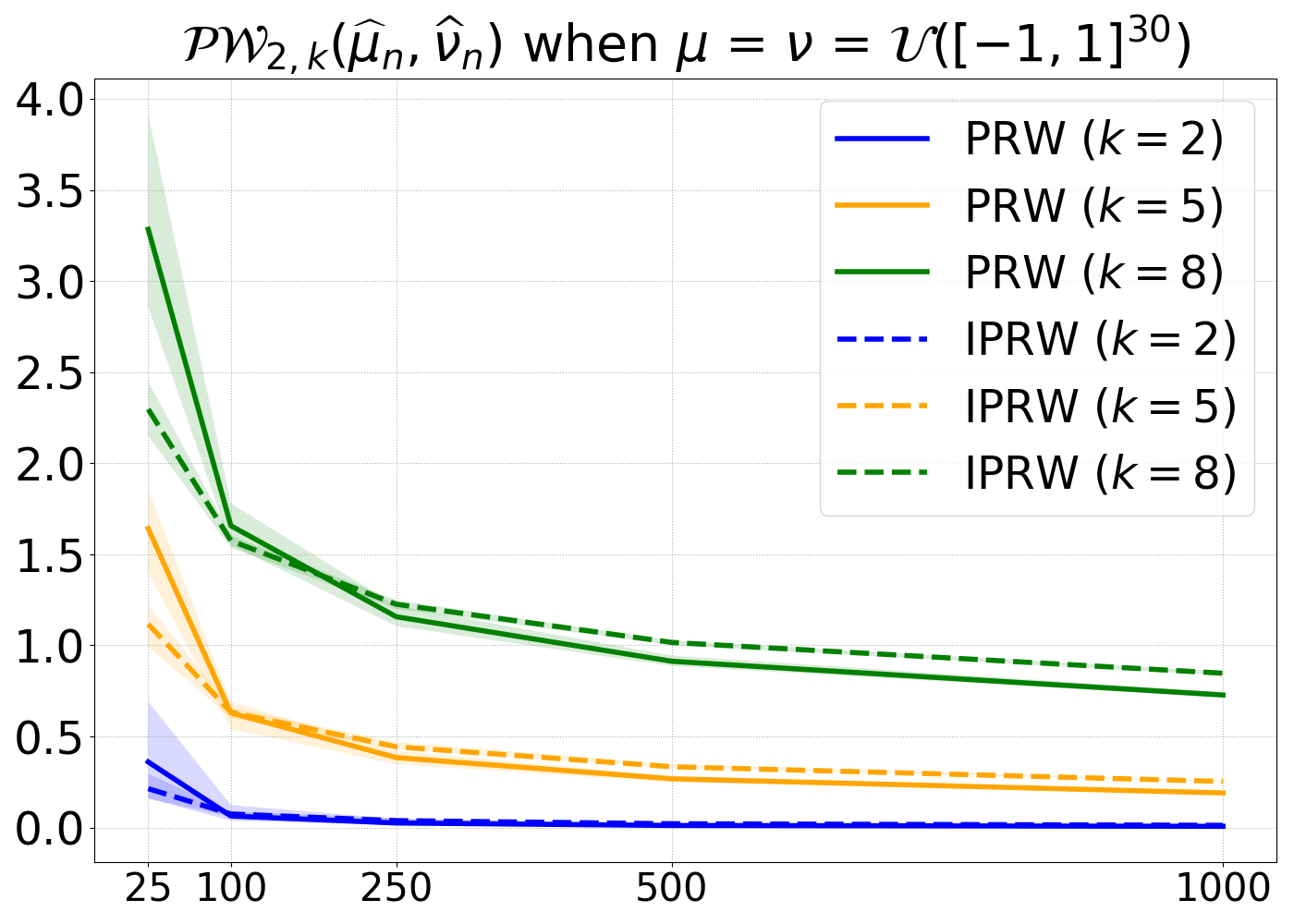}
\includegraphics[width=0.32\textwidth]{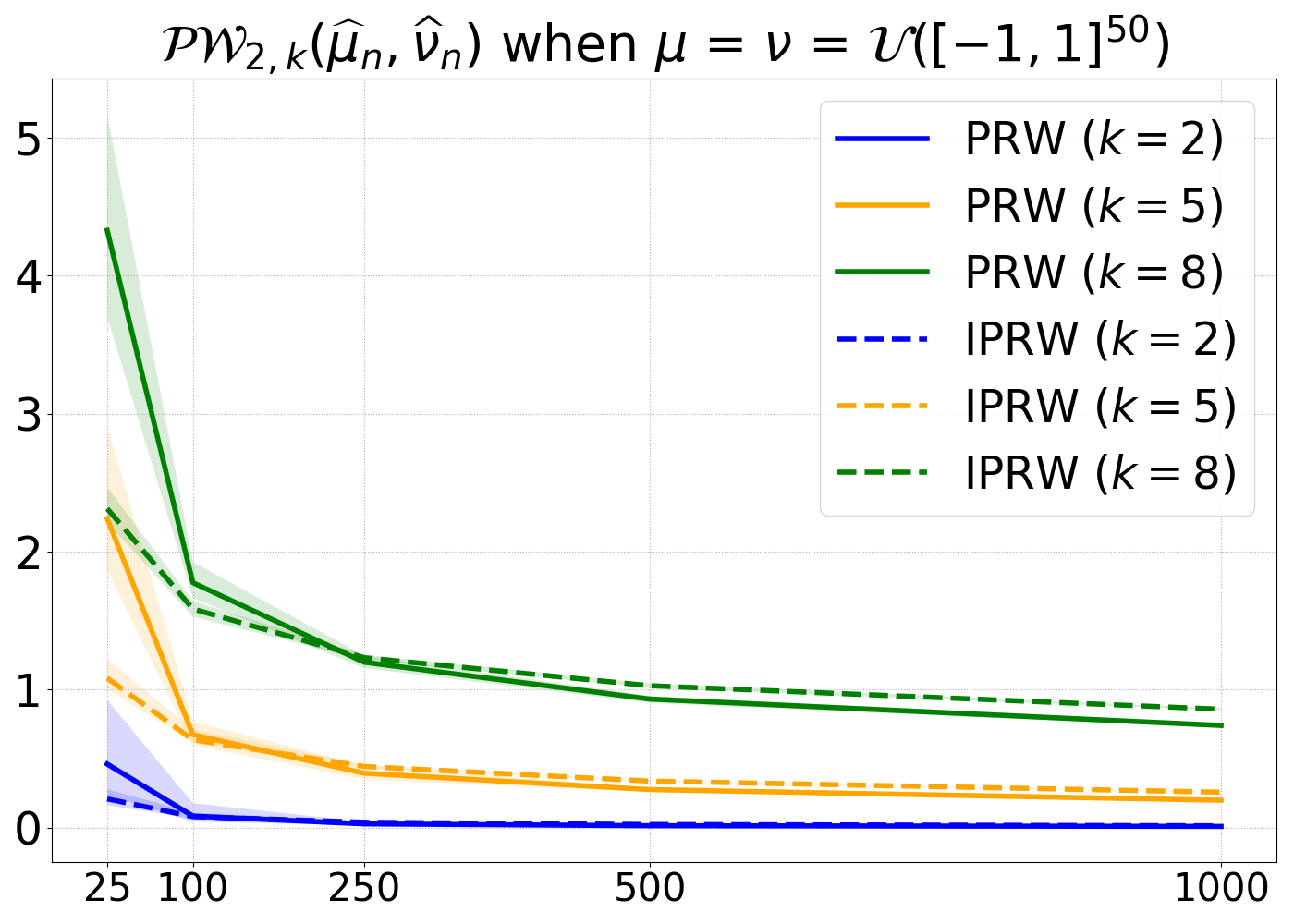}
\includegraphics[width=0.32\textwidth]{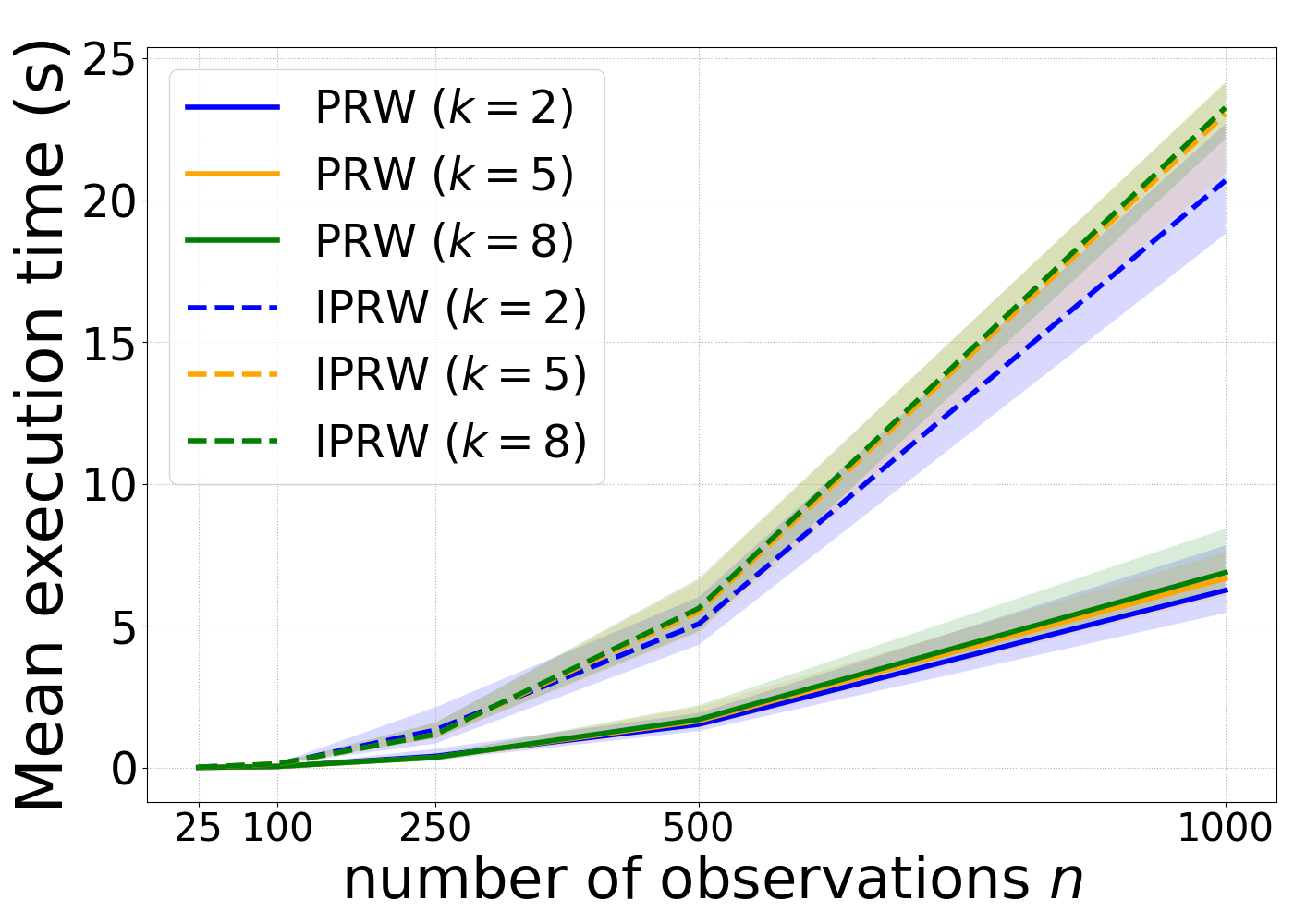}
\includegraphics[width=0.32\textwidth]{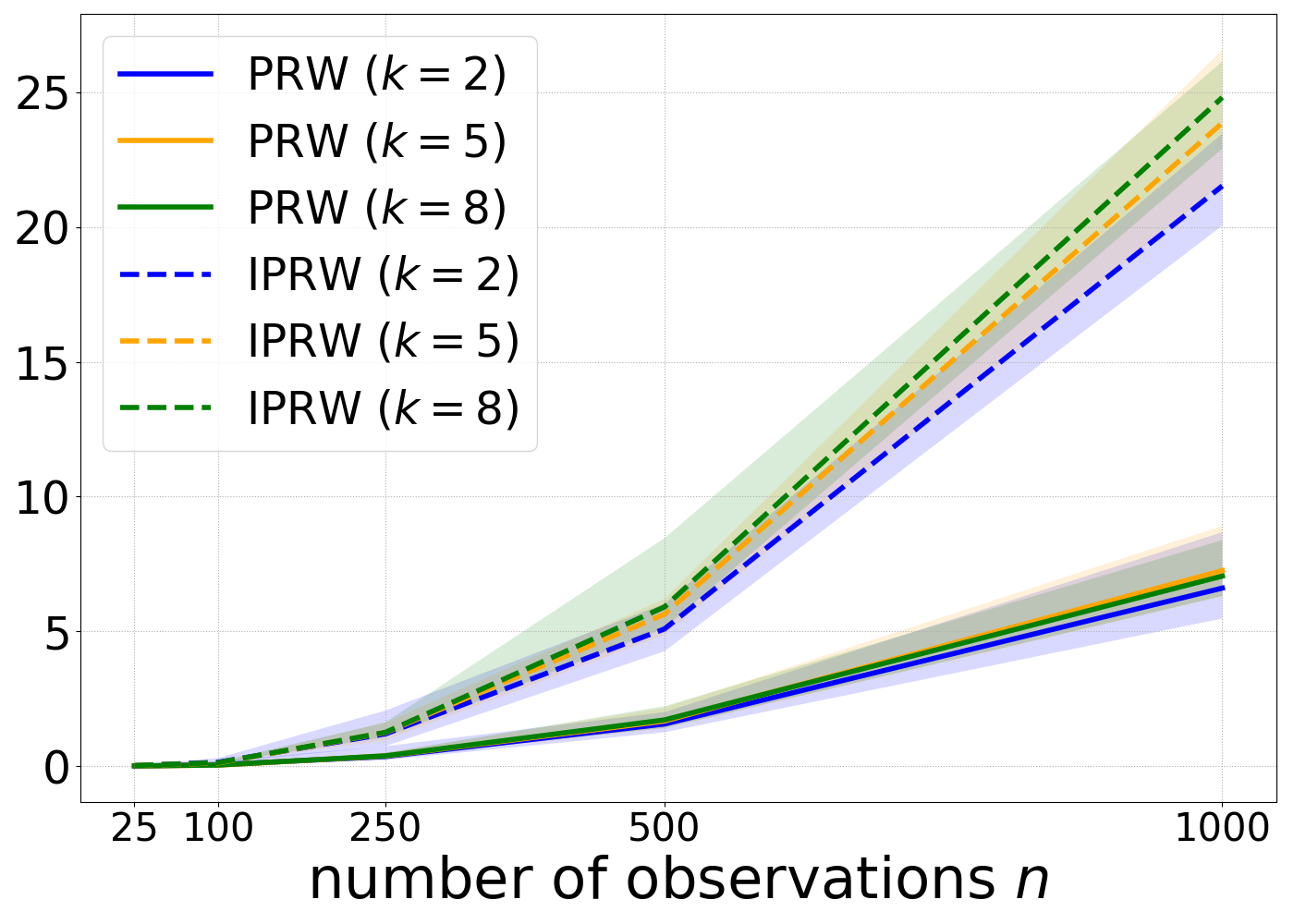}
\includegraphics[width=0.32\textwidth]{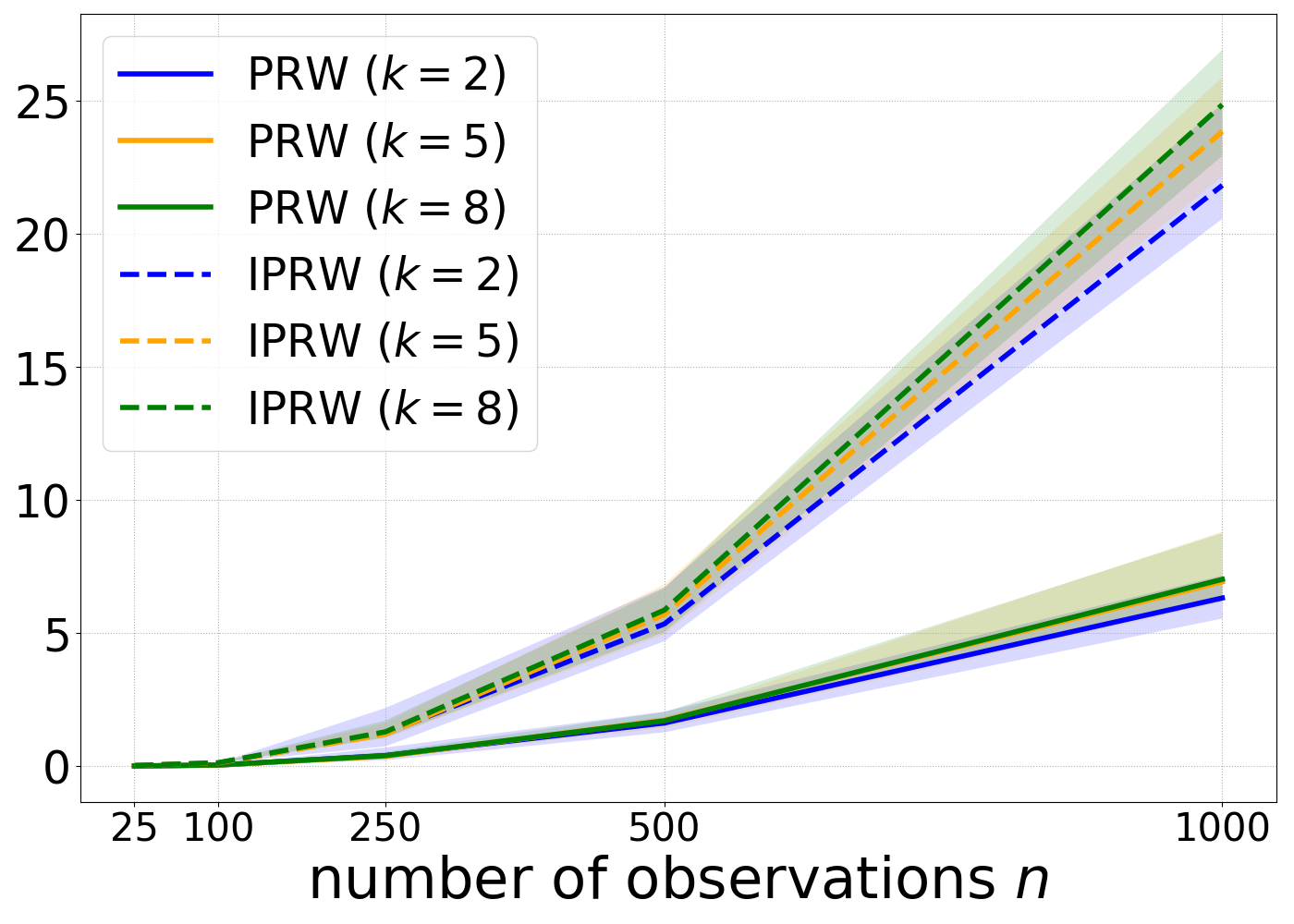}
\caption{\small{Mean values (Top) and mean computational time (Bottom) of the IPRW and PRW distances of order 2 between empirical measures $\widehat{\mu}_n$ and $\widehat{\nu}_n$ as the number of points $n$ varies. Results are averaged over 100 runs.}}\label{fig:exp1_appendix}\vspace*{-1em}
\end{figure}
\section{Experimental Setup}\label{appendix:setup}
\paragraph{Computing infrastructure.} For the experiments on the uniform distribution over hypercube, we implement in Python 3.7 with Numpy 1.18 on a workstation with an Intel Core i5-9400F (6 cores and 6 threads) and 32GB memory, equipped with Ubuntu 18.04. For the experiments on MPRW and MEPRW estimators, we implement in Python 2.7 with Numpy 1.16 and IPython 5.8 on the same machine. These experiments were not conducted with GPU. For the experiments on neural networks, we implement on the same machine with 2 GPUs (GeForce GTX 1070 and GeForce GTX 2070). 

\paragraph{Convergence and concentration.} We conduct the experiment on the uniform distribution over different hypercubes which are also used in the experiment~\citep{Paty-2019-Subspace}. In particular, we consider $\mu = \nu = \UCal([-v, v]^d)$ which is an uniform distribution over an hypercube and where $d$ and $v$ stand for the dimension and scale of the distribution respectively. $\widehat{\mu}_n$ and $\widehat{\nu}_n$ are empirical distributions corresponding to $\mu$ and $\nu$ with $n$ samples. We evaluate the PRW and IPRW distance in terms of mean values and mean computational times over $100$ runs for $(d, v) \in \{(10, 1), (10, 3), (30, 1), (30, 5), (50, 1), (50, 5)\}$. For the PRW distance, we run Algorithm~\ref{alg:supergrad-simplex} with \textsc{emd} solver in the \textsc{POT} package~\citep{Flamary-2017-Pot} and terminate the algorithm either when the maximum number of iterations $T=30$ is reached or when $\|U_{t+1} - U_t\|_F \leq 10^{-6}$. For the IPRW distance, we draw 100 uniform and independent projections from $\bs_{d, k}$ and compute each Wasserstein distances using \textsc{emd} solver in the \textsc{POT} package again. 
\begin{figure}[!t]
\centering
\begin{subfigure}{0.48\textwidth}
\includegraphics[width=1.0\textwidth]{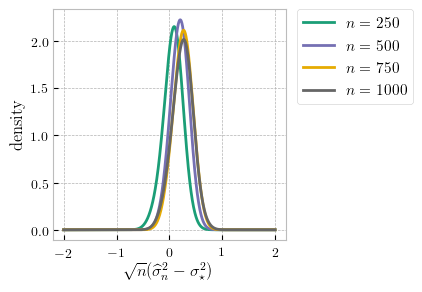}
\caption{\small{Mixture of 12 Gaussian distributions}}
\end{subfigure}
\begin{subfigure}{0.48\textwidth}
\includegraphics[width=1.0\textwidth]{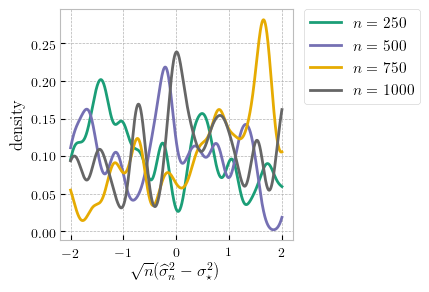}
\caption{\small{Mixture of 25 Gaussian distributions}}
\end{subfigure}
\caption{\small{Probability density of estimation of centered and rescaled $\widehat{\sigma}_n$ on the Gaussian model for different $n$.}}\label{fig:exp2_clt_appendix}\vspace*{-1em}
\end{figure}
\paragraph{Model misspecification.} We conduct the experiments on three type of data: the mixture of 8, 12 and 25 Gaussian distributions with Gaussian models $\MCal_1 = \{\NCal(\textbf{m}, \sigma^2\textbf{I}): \textbf{m} \in \br^2, \sigma^2 > 0\}$ and elliptically contoured stable models $\MCal_2 = \{\ECal\alpha\SCal_c(\textbf{I}, \textbf{m}): \textbf{m} \in \br^2\}$. For data-generating process, we fix $k$ centers $\{(a_i, b_i)\}_{1 \leq i \leq k}$. For each sample, we first randomly select $\textbf{m}$ from the centers at uniform and then draw the sample from $\NCal(2\textbf{m}, 0.01)$. For the mixture of 8 and 12 Gaussian distributions, the fixed set of centers are evenly distributed around a unit circle. For the mixture of 25 Gaussian distributions, the fixed set of centers are 25 grid points in $[-2, 2]^2$.

We use the ADAM optimization method with the default parameter setting to compute the MPRW and MEPRW estimators. At each inner loop, we run 5 projected supergradient ascent with the learning rate $10^{-3}$. For the Gaussian models, we estimate the densities of $\widehat{\sigma}_n^2$ with a kernel density estimator by computing 100 times MPRW estimator of order 1. The maximum number of ADAM iterations is set as 20000. To illustrate the consistency of MPRW and MEPRW estimators, we compute 100 times MPRW and MEPRW estimators of order 2, where the maximum number of ADAM iterations are set as 20000 and 10000 respectively. We also verify the convergence of MEPRW to MPRW by computing 100 times these estimators on a fixed set of $n=2000$ observations for different $m$ generated samples from the model. The maximum number of ADAM iterations for MPRW and MEPRW estimators are set as 20000 and 10000. For the elliptically contoured stable models, we verify the consistency property of MEPRW and the convergence of MEPRW to MPRW. For the former one, we compute 100 times MEPRW estimator of order 2 and set the maximum number of ADAM iterations as 10000. For the latter one, we compute 100 times MPRW and MEPRW estimators of order 2 on a fixed set of $n=100$ observations for different $m$ generated samples from the model. The maximum number of ADAM iterations are set as 20000 and 10000. All of these settings are consistently used on the mixture of 8, 12 and 25 Gaussian distributions.   

\paragraph{Generative modeling.} The procedure of the max-SW generator is summarized as follows: we first sample a random variable $Z$ from a fixed distribution on the base space $\ZCal$, and then transforms $Z$ through a neural network parametrized by $\theta$. This provides a parametric function $T_\theta: \ZCal \rightarrow \br^d$ which allows us to generate images from a distribution $\mu_\theta$. Our goal is to optimize the neural network parameters $\theta$ by minimizing the max-SW distance~\citep{Deshpande-2019-Max} between $\mu_\theta$ and data-generating distribution. We use a neural network with the fully-connected configuration from~\citet[Appendix D]{Deshpande-2018-Generative} and train our model with \textsc{CIFAR10}\footnote{Available in https://www.cs.toronto.edu/~kriz/cifar.html} and \textsc{ImageNet200}\footnote{Available in https://tiny-imagenet.herokuapp.com/}. The former one consists of 60000 and 10000 images of size $3 \times 32 \times 32$ for training and testing while the latter one consists of 100000 and 10000 images for training and testing. We use the minimal expected max-SW estimator of order 2 approximated with $50$ projected gradient ascent steps and $10^{-4}$ learning rate. We train for 1000 iterations with the ADAM optimizer~\citep{Kingma-2015-ADAM} and $10^{-4}$ learning rate. 
\begin{figure}[!t]
\centering
\begin{subfigure}{0.32\textwidth}
\includegraphics[width=1.0\textwidth]{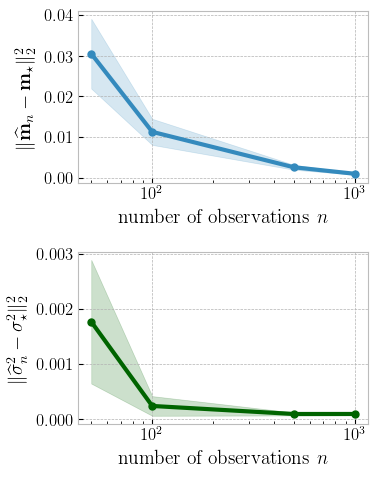}
\caption{\small{MPRW vs. $n$}}
\end{subfigure}
\begin{subfigure}{0.32\textwidth}
\includegraphics[width=1.0\textwidth]{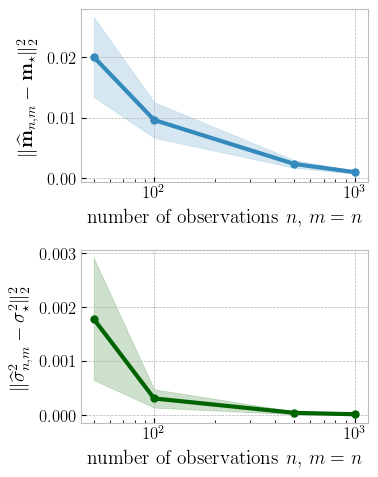}
\caption{\small{MEPRW vs. $n=m$}}
\end{subfigure}
\begin{subfigure}{0.32\textwidth}
\includegraphics[width=1.0\textwidth]{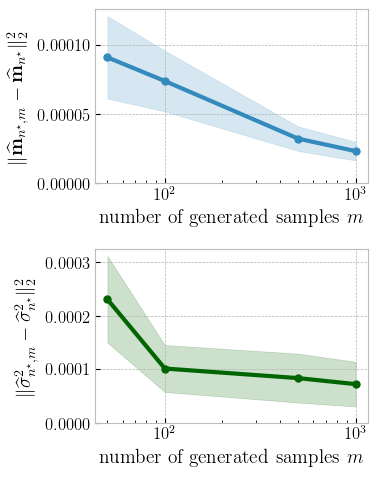}
\caption{\small{MEPRW with $n=2000$ vs. $m$}}
\end{subfigure}
\caption{\small{Minimal PRW and expected PRW estimations using Gaussian models and $n$ samples from the mixture of 12 Gaussian distributions. Results are averaged over 100 runs and shaded areas represent standard deviation.}}\label{fig:exp2_consistent_12_appendix}
\end{figure}
\begin{figure}[!t]
\centering
\begin{subfigure}{0.32\textwidth}
\includegraphics[width=1.0\textwidth]{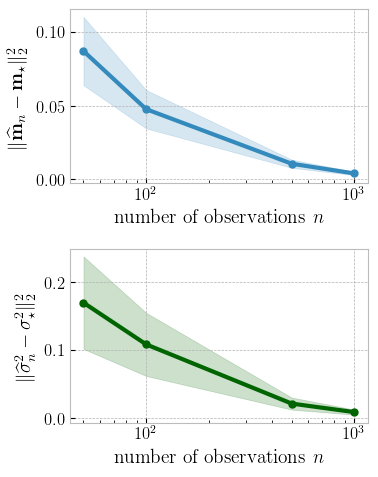}
\caption{\small{MPRW vs. $n$}}
\end{subfigure}
\begin{subfigure}{0.32\textwidth}
\includegraphics[width=1.0\textwidth]{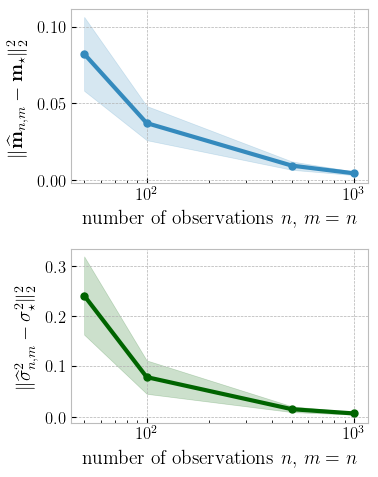}
\caption{\small{MEPRW vs. $n=m$}}
\end{subfigure}
\begin{subfigure}{0.32\textwidth}
\includegraphics[width=1.0\textwidth]{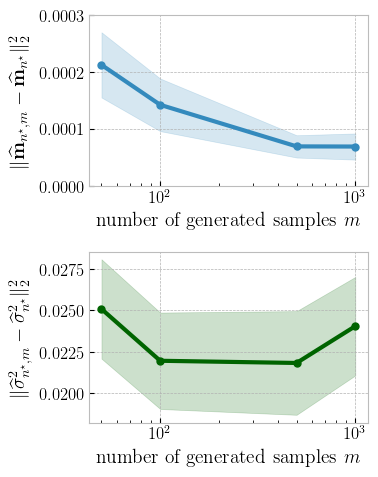}
\caption{\small{MEPRW with $n=2000$ vs. $m$}}
\end{subfigure}
\caption{\small{Minimal PRW and expected PRW estimations using Gaussian models and $n$ samples from the mixture of 25 Gaussian distributions. Results are averaged over 100 runs and shaded areas represent standard deviation.}}\label{fig:exp2_consistent_25_appendix}
\end{figure}

\section{Additional Experimental Results}\label{appendix:results}
\paragraph{Convergence and concentration.} Figure~\ref{fig:exp1_appendix} presents average distances and computational times for $(d, v) \in \{(10, 5), (30, 1), (50, 1)\}$, where the shaded areas show the max-min values over 100 runs. We also observe that the IPRW distance is smaller than the PRW distance for small $n$, especially so when $d$ and $v$ are large. The two distances are close when $n$ is large, supporting the theoretical results given by Theorem~\ref{Theorem:IPRW-main} and Theorem~\ref{Theorem:PRW-main-Poincare} in practice. The computation of the PRW distance is relatively faster than that of the IPRW distance in these computations. 
\begin{figure}[!ht]
\centering
\begin{subfigure}{0.4\textwidth}
\includegraphics[width=1.0\textwidth]{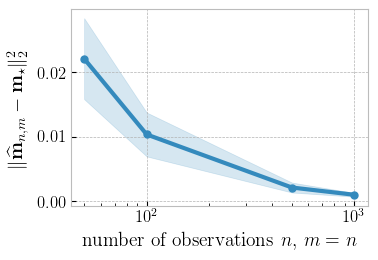}
\includegraphics[width=1.0\textwidth]{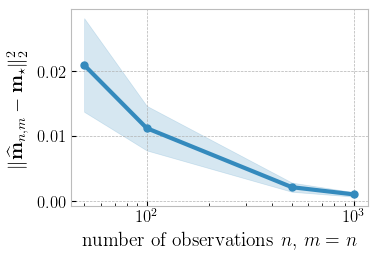}
\includegraphics[width=1.0\textwidth]{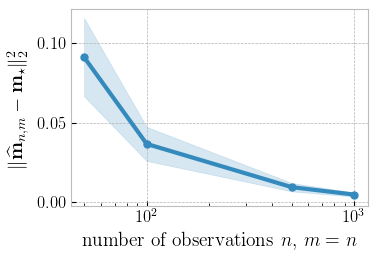}
\caption{\small{MEPRW}}
\end{subfigure}
\begin{subfigure}{0.4\textwidth}
\includegraphics[width=1.0\textwidth]{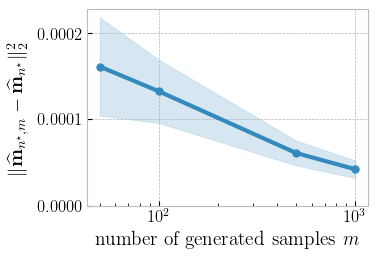}
\includegraphics[width=1.0\textwidth]{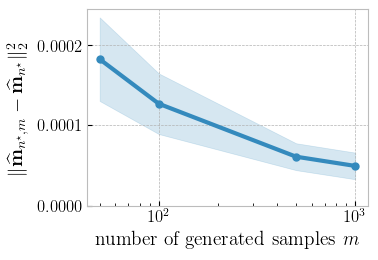}
\includegraphics[width=1.0\textwidth]{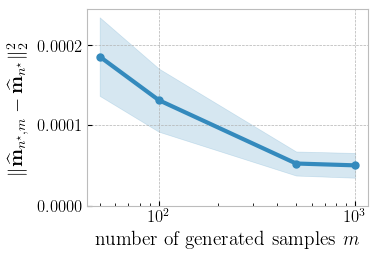}
\caption{\small{MEPRW, $n^\star = 100$}}
\end{subfigure}
\caption{\small{Minimal expected PRW estimations using elliptically contoured stable models and $n$ samples from the mixture of 8 Gaussian distributions (top), 12 Gaussian distributions (middle) and 25 Gaussian distributions (bottom), and $m$ samples generated from the model. Results are averaged over 100 runs and shaded areas represent standard deviation.}}\label{fig:exp2_consistency_ECS_appendix}\vspace*{-1em}
\end{figure}
\paragraph{Model misspecification: Gaussian models.} Figure~\ref{fig:exp2_clt_appendix} shows the distributions centered and rescaled by $\sqrt{n}$ for a range of moderately large $n$, based on the two underlying models including the mixture of 12 Gaussian distributions and the mixture of 25 Gaussian distributions. The left figure supports the convergence rate and the limiting distribution of the estimator as derived in Theorem~\ref{Theorem:general} on the mixture of 12 Gaussian distributions. The right figure suggests that the limiting distribution is not normal when the underlying model is given by the mixture of 25 Gaussian distributions. For the latter case, the result is not as anticipated by  Theorem~\ref{Theorem:general}. This is possibly because we only conduct 5 projected supergradient ascent at each inner loop, which may not be enough to achieve a good approximate projection $u \in \bs_{d, 1}$.

Figure~\ref{fig:exp2_consistent_12_appendix} and~\ref{fig:exp2_consistent_25_appendix} demonstrate the large-sample consistency behavior of MPRW and MEPRW estimators on the mixture of 12 and 25 Gaussian distributions, which are expected since Assumption~\ref{A1}-\ref{A3} are mild. The MEPRW estimator also converges to the MPRW estimator on the mixture of 12 Gaussian distributions, confirming Theorem~\ref{Theorem:convergence-MEPRW-MPRW}. One exception in these experiments is the failure of convergence of MEPRW to MPRW on the mixture of 25 Gaussian distributions. Apparently, the results from Theorem~\ref{Theorem:convergence-MEPRW-MPRW} do not hold in this experiment setting. This is likely due to the violation of Assumption~\ref{A5} that is necessary for Theorem~\ref{Theorem:convergence-MEPRW-MPRW} to hold.

\paragraph{Model misspecification: Elliptically contoured stable models.} Figure~\ref{fig:exp2_consistency_ECS_appendix} (a) illustrates the consistency of the MEPRW estimator $\widehat{\textbf{m}}_{n, m}$, approximated with 5 projected supergradient ascent, the same way as for the Gaussian models. Figure~\ref{fig:exp2_consistency_ECS_appendix} (b) confirms the convergence of $\widehat{\textbf{m}}_{n, m}$ to the MPRW estimator $\widehat{\textbf{m}}_n$, where we fix $n=100$ observations and compute the mean squared error between these two estimators (using 5 projected supergradient ascent) for different values of $m$. Note that the MPRW estimator is approximated with the MEPRW obtained for a large enough value of $m$: $\widehat{\textbf{m}}_n = \widehat{\textbf{m}}_{n, 10^4}$. To this end, our results on elliptically contoured stable models confirm Theorem~\ref{Theorem:consistency-MPRW}, Theorem~\ref{Theorem:consistency-MEPRW} and Theorem~\ref{Theorem:convergence-MEPRW-MPRW} in practice. 

\begin{figure}
\centering\includegraphics[width=0.5\textwidth]{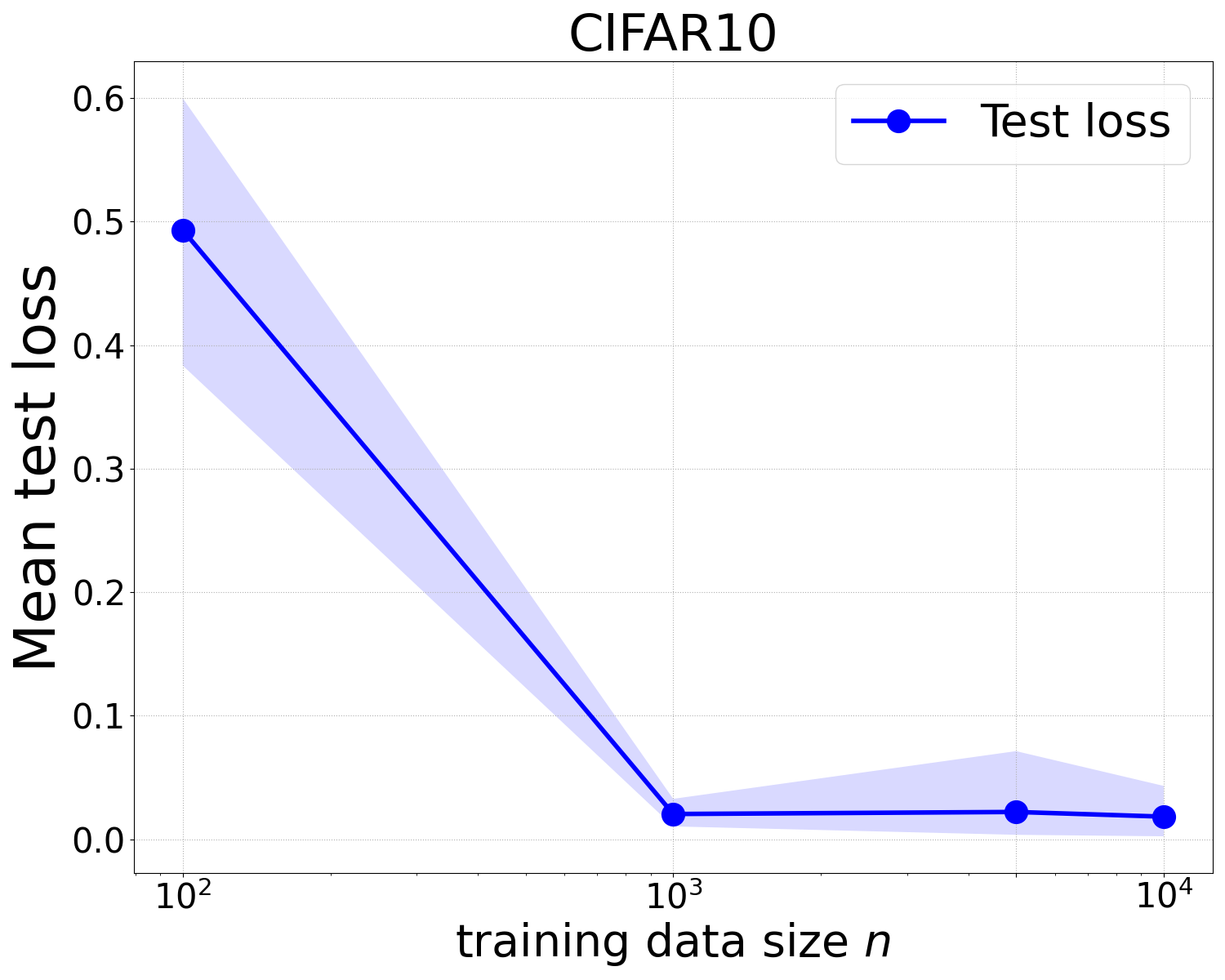}
\caption{\small{Mean test loss for different value of $(n, m)$ on \textsc{CIFAR10}.}}\label{fig:exp3_cifar10}
\end{figure}

\paragraph{Generative modeling.} Figure~\ref{fig:exp3_cifar10} presents the mean test loss on \textsc{CIFAR10} over 10 runs, where the shaded areas show the max-min values over the runs. Here the minimal expected max-SW estimator of order 2 is approximated with $20$ projected gradient ascent steps and $10^{-4}$ learning rate. We trained for 1000 iterations with the ADAM optimizer~\citep{Kingma-2015-ADAM} and $10^{-4}$ learning rate. We also train the NNs with $(n, m) \in \{(100, 20), (1000, 40), (5000, 60), (10000, 100)\}$ where $n$ is the number of training samples and $m$ is the number of generated samples and compute the testing losses using the trained models on the testing dataset ($n=10000$) with $m=250$ generated samples. We compare these testing losses to that of a NN trained using $n = 60000$ (i.e., the entire training dataset) and $m=200$ and present them in Figure~\ref{fig:exp3_cifar10}. Again, our results confirm Theorem~\ref{Theorem:consistency-MEPRW} in practice. 

\end{document}